\documentclass[11pt,leqno]{article}
\normalfont
\renewcommand{\baselinestretch}{1.15}

\usepackage{amsfonts}

\newfont{\eulercursive}{eurm10 at 11pt}
\newcommand{\myl}{\mbox{\eulercursive `}}

\newcommand{\QED}{\raisebox{0.5mm}{\fbox{\rule{0mm}{1.5mm}\ }}}

\newcounter{myfn}[page]

\newcommand{\myfootnote}[1]{\setcounter{footnote}{\value{myfn}}%
    \footnote{#1}\stepcounter{myfn}}

\newcommand{\MarsAttacksTheorem}{Theorem 1.1}
\newcommand{\WildbergerTheorem}{Theorem 1.2}  

\newcommand{\DynkinDiagramFigure}{Figure 1.1}
\newcommand{\BtwoMarsDemo}{Figure 1.2} 
\newcommand{\ECoxeterGraphFigure}{Figure 1.3} 
\newcommand{\EGCMTheorem}{Theorem 1.3} 
\newcommand{\MainTheorems}{Theorems 1.1 and 1.3}

\newcommand{\StrongConvergenceTheorem}{Theorem 2.1}
\newcommand{\StrongConvergenceCorollary}{Lemma 2.2}
\newcommand{\ComparisonTheorem}{Theorem 2.3}
\newcommand{\ComparisonLemma}{Lemma 2.4}
\newcommand{\ComparisonCorollary}{Lemma 2.5}
\newcommand{\ComparisonResults}{Lemmas 2.2 and 2.5}
\newcommand{\NotMarsFriendlyLemma}{Lemma 2.6}
\newcommand{\NotMarsFriendlyCatalog}{Lemma 2.7}
\newcommand{\NotMarsFriendlyFigure}{Figure 2.1}
\newcommand{\DtwittleFiveNodesFigure}{Figure 2.2}
\newcommand{\DtwittleFigure}{Figure 2.3}
\newcommand{\EveryNodeFiredLemma}{Lemma 2.8}
\newcommand{\SubgraphLemma}{Lemma 2.9}
\newcommand{\TwoNodeLemma}{Lemma 2.10}
\newcommand{\FourFamiliesFirstLemma}{Lemma 2.11}
\newcommand{\FourFamiliesSecondLemma}{Lemma 2.12}

\newcommand{\ErikssonWordProposition}{Proposition 3.1}
\newcommand{\ProofCorollary}{Corollary 3.2}
\newcommand{\CorollaryTheoremList}{Theorem 1.1 or 1.2}

\newcommand{\DeodharEquivalenceCorollary}{Corollary 4.1}
\newcommand{\TwoGeneratorAnalysis}{Lemma 4.2}
\newcommand{\HumphreysTheorem}{Proposition 4.3}
\newcommand{\PositiveToNegativeRootsLemma}{Lemma 4.4}
\newcommand{\DeodharProp}{Proposition 4.5}
\newcommand{\OneToOneLemma}{Lemma 4.6}
\newcommand{\ReducedLemma}{Lemma 4.7}
\newcommand{\TFAE}{Proposition 4.8}
\newcommand{\LengthProposition}{Proposition 4.9}
\newcommand{\UnitalExample}{Example 4.10}
\newcommand{\UnitalExampleFigure}{Figure 4.1}
\newcommand{\HRTProofResults}{Propositions 4.5 and 4.9}
\newcommand{\HRTResult}{Proposition 4.11} 
\newcommand{\TitsConeConvergenceResult}{Proposition 4.12}
\newcommand{\RuleOutInfinitePropositions}{Propositions 4.11 and 4.12}
\newcommand{\TitsConeProp}{Proposition 4.13}
\newcommand{\LemmaListForEGCM}{Lemmas 2.2, 2.4, 2.5, 2.6, 2.8, and 
2.9} 
\newcommand{\CycleLemma}{Lemma 4.14}
\newcommand{\ThreeCycleLemma}{Lemma 4.15}
\newcommand{\FourCycleLemma}{Lemma 4.16}
\newcommand{\RuleOutSmallCyclesLemmas}{Lemmas 4.14 and 4.16}

\newcommand{\KeyWJResultLemma}{Lemma 5.1}

\newcommand{\KeyWJResult}{Proposition 5.2}
\newcommand{\LengthUsingLongestWord}{Proposition 5.3} 
\newcommand{\NoRepeatPosRootsProp}{Proposition 5.4}
\newcommand{\PosRootsPropList}{Propositions 5.3 and 5.4}
\newcommand{\AllPositiveRootsProp}{Theorem 5.5}
\newcommand{\AdjacencyFreeLemma}{Lemma 5.6} 
\newcommand{\AdjacencyFreeForFinite}{Proposition 5.7}
\newcommand{\AdjacencyFreeClassification}{Theorem 5.8}

\newcommand{\IntroNum}{1}
\newcommand{\FirstProofNum}{2}
\newcommand{\FirstProofCorollaryNum}{3}
\newcommand{\SecondProofNum}{4}

\newcommand{\ConvergenceNum}{5}
\newcommand{\RemarksNum}{6}

\newcommand{\selt}{\mathbf{s}} \newcommand{\telt}{\mathbf{t}}

\newcommand{\disjointunion}{\setlength{\unitlength}{0.14cm}
\ 
\begin{picture}(2,2) 
\put(0,0){$\cup$}
\put(0.9675,1.5){\circle*{0.5}}
\end{picture}\ }

\newcommand{\CircleInteger}[1]{
\setlength{\unitlength}{0.14cm}
\begin{picture}(3,2) 
\put(1,1){\circle{2}}
\put(0.55,0.6){{\tiny #1}}
\end{picture}
}

\newcommand{\CircleIntegerm}{
\setlength{\unitlength}{0.14cm}
\begin{picture}(3,2) 
\put(1,1){\circle{2}}
\put(0.3,0.6){{\tiny $m$}}
\end{picture}
}

\newcommand{\CircleInfty}{
\setlength{\unitlength}{0.14cm}
\begin{picture}(3,2) 
\put(1,1){\circle{2}}
\put(0.1,0.6){{\tiny $\infty$}}
\end{picture}
}

\newcommand{\BowTie}{\setlength{\unitlength}{1.55in}
\begin{picture}(2.5,1.2)
\put(0,0){
\begin{picture}(0,1.2)
            \put(0.8,0.6){\line(1,0){0.81}}
            \put(1.1,0.625){\CircleInfty}
            \put(0.8,0.6){\vector(1,0){0.2}}
            \put(1.6,0.6){\vector(-1,0){0.2}}
            \put(0.95,0.51){\footnotesize $2$}
            \put(1.4,0.51){\footnotesize $2$}
            \put(0.3,0.1){\line(0,1){1}}
            \put(0.3,0.1){\line(1,1){0.5}}
            \put(0.3,1.1){\line(1,-1){0.5}}
            \put(0.8,0.6){\circle*{0.05}}
            \put(0.3,1.1){\circle*{0.05}}
            \put(0.3,0.1){\circle*{0.05}}
            \put(0.85,0.675){\small $\gamma_{i}$}
            \put(0.125,0.05){\small $\gamma_{k}$}
            \put(0.125,1.1){\small $\gamma_{j}$}
            \put(0.3,0.3){\vector(0,1){0.1}}
            \put(0.195,0.35){\footnotesize $1$}
            \put(0.3,0.9){\vector(0,-1){0.1}}
            \put(0.195,0.8){\footnotesize $3$}
            \put(0.255,0.575){\CircleInteger{6}}
            \put(0.3,0.1){\vector(1,1){0.15}}
            \put(0.8,0.6){\vector(-1,-1){0.15}}
            \put(0.45,0.15){\footnotesize $2$}
            \put(0.675,0.375){\footnotesize $1$}
            \put(0.425,0.375){\CircleInteger{4}}
            \put(0.3,1.1){\vector(1,-1){0.15}}
            \put(0.8,0.6){\vector(-1,1){0.15}}
            \put(0.425,0.755){\CircleInteger{3}}
            \put(0.435,1){\footnotesize $1/5$}
            \put(0.69,0.75){\footnotesize $5$}
\end{picture}
}
\put(1.6,0){
\begin{picture}(0.6,1.2)
            \put(0,0.6){\circle*{0.05}}
            \put(0.5,0.1){\circle*{0.05}}
            \put(0.5,1.1){\circle*{0.05}}
            \put(0,0.6){\line(1,1){0.5}}
            \put(0,0.6){\line(1,-1){0.5}}
            \put(0.5,0.1){\line(0,1){1}}
            \put(0.5,0.3){\vector(0,1){0.1}}
            \put(0.55,0.35){\footnotesize $1+\sqrt{5}$}
            \put(0.5,0.9){\vector(0,-1){0.1}}
            \put(0.55,0.8){\footnotesize $(1+\sqrt{5})/4$}
            \put(0,0.6){\vector(1,1){0.15}}
            \put(0.5,1.1){\vector(-1,-1){0.15}}
            \put(0.075,0.75){\footnotesize $7$}
            \put(0.25,1){\footnotesize $1/7$}
            \put(0,0.6){\vector(1,-1){0.15}}
            \put(0.5,0.1){\vector(-1,1){0.15}}
            \put(0,0.375){\footnotesize $7/2$}
            \put(0.25,0.15){\footnotesize $2/7$} 
            \put(-0.125,0.675){\small $\gamma_{x}$}
            \put(0.58,0.05){\small $\gamma_{z}$}
            \put(0.58,1.1){\small $\gamma_{y}$}
            \put(0.345,0.575){\CircleInteger{5}}
            \put(0.175,0.375){\CircleInteger{3}}
            \put(0.175,0.755){\CircleInteger{3}}
\end{picture}
}
\end{picture}
}

\newcommand{\AnEGraph}{\setlength{\unitlength}{0.75in}
\begin{picture}(6,0.85)
\put(0,0){\begin{picture}(1,0)
            \put(0,0.1){\circle*{0.075}}
            \put(1,0.1){\circle*{0.075}}
            \put(2,0.1){\circle*{0.075}}
            \put(4,0.1){\circle*{0.075}}
            \put(5,0.1){\circle*{0.075}}
            \put(6,0.1){\circle*{0.075}}
            \put(0,0.1){\line(1,0){2}}
            \multiput(2,0.1)(0.4,0){5}{\line(1,0){0.2}}
            \put(4,0.1){\line(1,0){2}}
           \end{picture}}
\end{picture}}

\newcommand{\BnEGraph}{\setlength{\unitlength}{0.75in}
\begin{picture}(6,0.55)
\put(0,0){\begin{picture}(1,0)
            \put(0,0.1){\circle*{0.075}}
            \put(1,0.1){\circle*{0.075}}
            \put(2,0.1){\circle*{0.075}}
            \put(4,0.1){\circle*{0.075}}
            \put(5,0.1){\circle*{0.075}}
            \put(6,0.1){\circle*{0.075}}
            \put(0,0.1){\line(1,0){2}}
            \multiput(2,0.1)(0.4,0){5}{\line(1,0){0.2}}
            \put(4,0.1){\line(1,0){2}}
            \put(5.3,0.15){\CircleInteger{4}}
           \end{picture}}
\end{picture}}

\newcommand{\DnEGraph}{\setlength{\unitlength}{0.75in}
\begin{picture}(6,0.75)
\put(0,-0.25){\begin{picture}(1,0)
            \put(0,0.35){\circle*{0.075}}
            \put(1,0.35){\circle*{0.075}}
            \put(2,0.35){\circle*{0.075}}
            \put(4,0.35){\circle*{0.075}}
            \put(5,0.35){\circle*{0.075}}
            \put(6,0.1){\circle*{0.075}}
            \put(6,0.6){\circle*{0.075}}
            \put(0,0.35){\line(1,0){2}}
            \multiput(2,0.35)(0.4,0){5}{\line(1,0){0.2}}
            \put(4,0.35){\line(1,0){1}}
            \put(5,0.35){\line(4,1){1}}
            \put(5,0.35){\line(4,-1){1}}
           \end{picture}}
\end{picture}}

\newcommand{\EEightEGraph}{\setlength{\unitlength}{0.75in}
\begin{picture}(6,0.75)
\put(0,-0.25){\begin{picture}(1,0)
            \put(0,0.1){\circle*{0.075}}
            \put(1,0.1){\circle*{0.075}}
            \put(2,0.1){\circle*{0.075}}
            \put(2,0.6){\circle*{0.075}}
            \put(3,0.1){\circle*{0.075}}
            \put(4,0.1){\circle*{0.075}}
            \put(5,0.1){\circle*{0.075}}
            \put(6,0.1){\circle*{0.075}}
            \put(0,0.1){\line(1,0){6}}
            \put(2,0.1){\line(0,1){0.5}}
           \end{picture}}
\end{picture}}

\newcommand{\ESevenEGraph}{\setlength{\unitlength}{0.75in}
\begin{picture}(6,0.75)
\put(0,-0.25){\begin{picture}(1,0)
            \put(0,0.1){\circle*{0.075}}
            \put(1,0.1){\circle*{0.075}}
            \put(2,0.1){\circle*{0.075}}
            \put(2,0.6){\circle*{0.075}}
            \put(3,0.1){\circle*{0.075}}
            \put(4,0.1){\circle*{0.075}}
            \put(5,0.1){\circle*{0.075}}
            \put(0,0.1){\line(1,0){5}}
            \put(2,0.1){\line(0,1){0.5}}
            \put(4.95,0.2){\small {\tt *}}
           \end{picture}}
\end{picture}}

\newcommand{\ESixEGraph}{\setlength{\unitlength}{0.75in}
\begin{picture}(6,0.75)
\put(0,-0.25){\begin{picture}(1,0)
            \put(0,0.1){\circle*{0.075}}
            \put(1,0.1){\circle*{0.075}}
            \put(2,0.1){\circle*{0.075}}
            \put(2,0.6){\circle*{0.075}}
            \put(3,0.1){\circle*{0.075}}
            \put(4,0.1){\circle*{0.075}}
            \put(0,0.1){\line(1,0){4}}
            \put(2,0.1){\line(0,1){0.5}}
            \put(-0.05,0.2){\small {\tt *}}
            \put(3.95,0.2){\small {\tt *}}
           \end{picture}}
\end{picture}}

\newcommand{\FFourEGraph}{\setlength{\unitlength}{0.75in}
\begin{picture}(1,0.85)
\put(0,0){\begin{picture}(1,0)
            \put(0,0.1){\circle*{0.075}}
            \put(1,0.1){\circle*{0.075}}
            \put(2,0.1){\circle*{0.075}}
            \put(3,0.1){\circle*{0.075}}
            \put(0,0.1){\line(1,0){3}}
            \put(1.3,0.15){\CircleInteger{4}}
            \end{picture}}
\end{picture}}

\newcommand{\HFourEGraph}{\setlength{\unitlength}{0.75in}
\begin{picture}(1,0.55)
\put(0,0){\begin{picture}(1,0)
            \put(0,0.1){\circle*{0.075}}
            \put(1,0.1){\circle*{0.075}}
            \put(2,0.1){\circle*{0.075}}
            \put(3,0.1){\circle*{0.075}}
            \put(0,0.1){\line(1,0){3}}
            \put(0.3,0.15){\CircleInteger{5}}
            \end{picture}}
\end{picture}}

\newcommand{\HThreeEGraph}{\setlength{\unitlength}{0.75in}
\begin{picture}(1,0.55)
\put(0,0){\begin{picture}(1,0)
            \put(0,0.1){\circle*{0.075}}
            \put(1,0.1){\circle*{0.075}}
            \put(2,0.1){\circle*{0.075}}
            \put(0,0.1){\line(1,0){2}}
            \put(0.3,0.15){\CircleInteger{5}}
            \put(1.95,0.2){\small {\tt *}}
            \end{picture}}
\end{picture}}

\newcommand{\ITwoEGraph}{\setlength{\unitlength}{0.75in}
\begin{picture}(1,0.55)
\put(0,0){\begin{picture}(1,0)
            \put(0,0.1){\circle*{0.075}}
            \put(1,0.1){\circle*{0.075}}
            \put(0,0.1){\line(1,0){1}}
            \put(0.3,0.15){\CircleIntegerm}
            \end{picture}}
\end{picture}}

\newcommand{\BTwoGraphForFigure}{\setlength{\unitlength}{1in}
\begin{picture}(1,0.4)
\put(0,0){\begin{picture}(1,0)
            \put(0,0.1){\circle*{0.05}}
            \put(-0.11,-0.05){\large $\gamma_{1}$}
            \put(1,0.1){\circle*{0.05}}
            \put(1.05,-0.05){\large $\gamma_{2}$}
            \put(0,0.1){\line(1,0){1}}
            \put(0.2,0.1){\vector(1,0){0.1}}
            \put(0.8,0.1){\vector(-1,0){0.1}}
            \put(0.7,0.1){\vector(-1,0){0.1}}
            \end{picture}}
\end{picture}}

\newcommand{\GTwoMarsOK}{\setlength{\unitlength}{0.75in}
\begin{picture}(1,0.55)
\put(0,0){\begin{picture}(1,0)
            \put(0,0.1){\circle*{0.075}}
            \put(1,0.1){\circle*{0.075}}
            \put(0,0.1){\line(1,0){1}}
            \put(0.2,0.1){\vector(1,0){0.1}}
            \put(0.8,0.1){\vector(-1,0){0.1}}
            \put(0.7,0.1){\vector(-1,0){0.1}}
            \put(0.6,0.1){\vector(-1,0){0.1}}
            \end{picture}}
\end{picture}}

\newcommand{\FFourMarsOK}{\setlength{\unitlength}{0.75in}
\begin{picture}(1,0.85)
\put(0,0){\begin{picture}(1,0)
            \put(0,0.1){\circle*{0.075}}
            \put(1,0.1){\circle*{0.075}}
            \put(2,0.1){\circle*{0.075}}
            \put(3,0.1){\circle*{0.075}}
            \put(0,0.1){\line(1,0){3}}
            \put(1.2,0.1){\vector(1,0){0.1}}
            \put(1.3,0.1){\vector(1,0){0.1}}
            \put(1.8,0.1){\vector(-1,0){0.1}}
            \end{picture}}
\end{picture}}

\newcommand{\CnMarsOK}{\setlength{\unitlength}{0.75in}
\begin{picture}(6,0.55)
\put(0,0){\begin{picture}(1,0)
            \put(0,0.1){\circle*{0.075}}
            \put(1,0.1){\circle*{0.075}}
            \put(2,0.1){\circle*{0.075}}
            \put(4,0.1){\circle*{0.075}}
            \put(5,0.1){\circle*{0.075}}
            \put(6,0.1){\circle*{0.075}}
            \put(0,0.1){\line(1,0){2}}
            \multiput(2,0.1)(0.4,0){5}{\line(1,0){0.2}}
            \put(4,0.1){\line(1,0){2}}
            \put(5.8,0.1){\vector(-1,0){0.1}}
            \put(5.7,0.1){\vector(-1,0){0.1}}
            \put(5.2,0.1){\vector(1,0){0.1}}
           \end{picture}}
\end{picture}}

\newcommand{\BnMarsOK}{\setlength{\unitlength}{0.75in}
\begin{picture}(6,0.55)
\put(0,0){\begin{picture}(1,0)
            \put(0,0.1){\circle*{0.075}}
            \put(1,0.1){\circle*{0.075}}
            \put(2,0.1){\circle*{0.075}}
            \put(4,0.1){\circle*{0.075}}
            \put(5,0.1){\circle*{0.075}}
            \put(6,0.1){\circle*{0.075}}
            \put(0,0.1){\line(1,0){2}}
            \multiput(2,0.1)(0.4,0){5}{\line(1,0){0.2}}
            \put(4,0.1){\line(1,0){2}}
            \put(5.8,0.1){\vector(-1,0){0.1}}
            \put(5.2,0.1){\vector(1,0){0.1}}
            \put(5.3,0.1){\vector(1,0){0.1}}
           \end{picture}}
\end{picture}}

\newcommand{\DnMarsOK}{\setlength{\unitlength}{0.75in}
\begin{picture}(6,0.75)
\put(0,-0.25){\begin{picture}(1,0)
            \put(0,0.35){\circle*{0.075}}
            \put(1,0.35){\circle*{0.075}}
            \put(2,0.35){\circle*{0.075}}
            \put(4,0.35){\circle*{0.075}}
            \put(5,0.35){\circle*{0.075}}
            \put(6,0.1){\circle*{0.075}}
            \put(6,0.6){\circle*{0.075}}
            \put(0,0.35){\line(1,0){2}}
            \multiput(2,0.35)(0.4,0){5}{\line(1,0){0.2}}
            \put(4,0.35){\line(1,0){1}}
            \put(5,0.35){\line(4,1){1}}
            \put(5,0.35){\line(4,-1){1}}
           \end{picture}}
\end{picture}}

\newcommand{\AnMarsOK}{\setlength{\unitlength}{0.75in}
\begin{picture}(6,0.85)
\put(0,0){\begin{picture}(1,0)
            \put(0,0.1){\circle*{0.075}}
            \put(1,0.1){\circle*{0.075}}
            \put(2,0.1){\circle*{0.075}}
            \put(4,0.1){\circle*{0.075}}
            \put(5,0.1){\circle*{0.075}}
            \put(6,0.1){\circle*{0.075}}
            \put(0,0.1){\line(1,0){2}}
            \multiput(2,0.1)(0.4,0){5}{\line(1,0){0.2}}
            \put(4,0.1){\line(1,0){2}}
           \end{picture}}
\end{picture}}

\newcommand{\EEightMarsOK}{\setlength{\unitlength}{0.75in}
\begin{picture}(6,0.75)
\put(0,-0.25){\begin{picture}(1,0)
            \put(0,0.1){\circle*{0.075}}
            \put(1,0.1){\circle*{0.075}}
            \put(2,0.1){\circle*{0.075}}
            \put(2,0.6){\circle*{0.075}}
            \put(3,0.1){\circle*{0.075}}
            \put(4,0.1){\circle*{0.075}}
            \put(5,0.1){\circle*{0.075}}
            \put(6,0.1){\circle*{0.075}}
            \put(0,0.1){\line(1,0){6}}
            \put(2,0.1){\line(0,1){0.5}}
           \end{picture}}
\end{picture}}

\newcommand{\ESevenMarsOK}{\setlength{\unitlength}{0.75in}
\begin{picture}(6,0.75)
\put(0,-0.25){\begin{picture}(1,0)
            \put(0,0.1){\circle*{0.075}}
            \put(1,0.1){\circle*{0.075}}
            \put(2,0.1){\circle*{0.075}}
            \put(2,0.6){\circle*{0.075}}
            \put(3,0.1){\circle*{0.075}}
            \put(4,0.1){\circle*{0.075}}
            \put(5,0.1){\circle*{0.075}}
            \put(0,0.1){\line(1,0){5}}
            \put(2,0.1){\line(0,1){0.5}}
           \end{picture}}
\end{picture}}

\newcommand{\ESixMarsOK}{\setlength{\unitlength}{0.75in}
\begin{picture}(6,0.75)
\put(0,-0.25){\begin{picture}(1,0)
            \put(0,0.1){\circle*{0.075}}
            \put(1,0.1){\circle*{0.075}}
            \put(2,0.1){\circle*{0.075}}
            \put(2,0.6){\circle*{0.075}}
            \put(3,0.1){\circle*{0.075}}
            \put(4,0.1){\circle*{0.075}}
            \put(0,0.1){\line(1,0){4}}
            \put(2,0.1){\line(0,1){0.5}}
           \end{picture}}
\end{picture}}

\newcommand{\TwoCitiesGraphWithLabels}{
\setlength{\unitlength}{1in}
\begin{picture}(1,0.25)
\put(0,0){\begin{picture}(1,0)
            \put(0,0.1){\circle*{0.05}}
            \put(-0.15,-0.05){\large $\gamma_{1}$}
            \put(1,0.1){\circle*{0.05}}
            \put(1.05,-0.05){\large $\gamma_{2}$}
            \put(0,0.1){\line(1,0){1}}
            \put(0.2,0.1){\vector(1,0){0.1}}
            \put(0.8,0.1){\vector(-1,0){0.1}}
            \put(0.225,0){\footnotesize $p$}
            \put(0.71,0){\footnotesize $q$}
            \end{picture}}
\end{picture}}

\newcommand{\ATwoGraphWithLabels}{\setlength{\unitlength}{1in}
\begin{picture}(1,0.4)
\put(0.4,0.3){\large $A_{2}$}
\put(0,0){\begin{picture}(1,0)
            \put(0,0.1){\circle*{0.05}}
            \put(-0.15,-0.05){\large $\gamma_{1}$}
            \put(1,0.1){\circle*{0.05}}
            \put(1.05,-0.05){\large $\gamma_{2}$}
            \put(0,0.1){\line(1,0){1}}
            \put(0.2,0.1){\vector(1,0){0.1}}
            \put(0.8,0.1){\vector(-1,0){0.1}}
            \end{picture}}
\end{picture}}

\newcommand{\ATwoGraphNoEdgeLabels}{\setlength{\unitlength}{0.75in}
\begin{picture}(1.65,0.25)
\put(0.25,0){\begin{picture}(1,0)
            \put(0,0.1){\circle*{0.05}}
            \put(-0.20,-0.05){\large $\gamma_{1}$}
            \put(1,0.1){\circle*{0.05}}
            \put(1.05,-0.05){\large $\gamma_{2}$}
            \put(0,0.1){\line(1,0){1}}
            \end{picture}}
\end{picture}}

\newcommand{\EGCMGraphCirclem}{\setlength{\unitlength}{0.75in}
\begin{picture}(1.65,0.35)
\put(0.25,0){\begin{picture}(1,0)
            \put(0,0.1){\circle*{0.05}}
            \put(-0.20,-0.05){\large $\gamma_{1}$}
            \put(1,0.1){\circle*{0.05}}
            \put(1.05,-0.05){\large $\gamma_{2}$}
            \put(0,0.1){\line(1,0){1}}
            \put(0.3,0.15){\CircleIntegerm}
            \end{picture}}
\end{picture}}

\newcommand{\BTwoGraphWithLabels}{\setlength{\unitlength}{1in}
\begin{picture}(1,0.4)
\put(0.4,0.3){\large $B_{2}$}
\put(0,0){\begin{picture}(1,0)
            \put(0,0.1){\circle*{0.05}}
            \put(-0.15,-0.05){\large $\gamma_{1}$}
            \put(1,0.1){\circle*{0.05}}
            \put(1.05,-0.05){\large $\gamma_{2}$}
            \put(0,0.1){\line(1,0){1}}
            \put(0.2,0.1){\vector(1,0){0.1}}
            \put(0.8,0.1){\vector(-1,0){0.1}}
            \put(0.7,0.1){\vector(-1,0){0.1}}
            \end{picture}}
\end{picture}}

\newcommand{\GTwoGraphWithLabels}{\setlength{\unitlength}{1in}
\begin{picture}(1,0.4)
\put(0.4,0.3){\large $G_{2}$}
\put(0,0){\begin{picture}(1,0)
            \put(0,0.1){\circle*{0.05}}
            \put(-0.15,-0.05){\large $\gamma_{1}$}
            \put(1,0.1){\circle*{0.05}}
            \put(1.05,-0.05){\large $\gamma_{2}$}
            \put(0,0.1){\line(1,0){1}}
            \put(0.2,0.1){\vector(1,0){0.1}}
            \put(0.8,0.1){\vector(-1,0){0.1}}
            \put(0.7,0.1){\vector(-1,0){0.1}}
            \put(0.6,0.1){\vector(-1,0){0.1}}
            \end{picture}}
\end{picture}}

\newcommand{\GTwoTwittle}{\setlength{\unitlength}{0.75in}
\begin{picture}(7.5,2)
\put(1.25,1.5){\begin{picture}(1,0)
            \put(0,0.1){\circle*{0.075}}
            \put(1,0.1){\circle*{0.075}}
            \put(2,0.1){\circle*{0.075}}
            \put(0,0.1){\line(1,0){2}}
            \put(0.2,0.1){\vector(1,0){0.1}}
            \put(0.8,0.1){\vector(-1,0){0.1}}
            \put(0.7,0.1){\vector(-1,0){0.1}}
            \put(0.6,0.1){\vector(-1,0){0.1}}
            \put(1.2,0.1){\vector(1,0){0.1}}
            \put(1.8,0.1){\vector(-1,0){0.1}}
            \put(1.7,0.1){\vector(-1,0){0.1}}
            \put(1.6,0.1){\vector(-1,0){0.1}}
            \end{picture}}
\put(4.25,1.5){\begin{picture}(1,0)
            \put(0,0.1){\circle*{0.075}}
            \put(1,0.1){\circle*{0.075}}
            \put(2,0.1){\circle*{0.075}}
            \put(0,0.1){\line(1,0){2}}
            \put(0.2,0.1){\vector(1,0){0.1}}
            \put(0.8,0.1){\vector(-1,0){0.1}}
            \put(0.7,0.1){\vector(-1,0){0.1}}
            \put(0.6,0.1){\vector(-1,0){0.1}}
            \put(1.8,0.1){\vector(-1,0){0.1}}
            \put(1.2,0.1){\vector(1,0){0.1}}
            \put(1.3,0.1){\vector(1,0){0.1}}
            \put(1.4,0.1){\vector(1,0){0.1}}
            \end{picture}}
\put(1.25,0.75){\begin{picture}(1,0)
            \put(0,0.1){\circle*{0.075}}
            \put(1,0.1){\circle*{0.075}}
            \put(2,0.1){\circle*{0.075}}
            \put(0,0.1){\line(1,0){2}}
            \put(0.2,0.1){\vector(1,0){0.1}}
            \put(0.8,0.1){\vector(-1,0){0.1}}
            \put(0.7,0.1){\vector(-1,0){0.1}}
            \put(0.6,0.1){\vector(-1,0){0.1}}
            \put(1.2,0.1){\vector(1,0){0.1}}
            \put(1.8,0.1){\vector(-1,0){0.1}}
            \put(1.7,0.1){\vector(-1,0){0.1}}
            \end{picture}}
\put(4.25,0.75){\begin{picture}(1,0)
            \put(0,0.1){\circle*{0.075}}
            \put(1,0.1){\circle*{0.075}}
            \put(2,0.1){\circle*{0.075}}
            \put(0,0.1){\line(1,0){2}}
            \put(0.2,0.1){\vector(1,0){0.1}}
            \put(0.8,0.1){\vector(-1,0){0.1}}
            \put(0.7,0.1){\vector(-1,0){0.1}}
            \put(0.6,0.1){\vector(-1,0){0.1}}
            \put(1.8,0.1){\vector(-1,0){0.1}}
            \put(1.2,0.1){\vector(1,0){0.1}}
            \put(1.3,0.1){\vector(1,0){0.1}}
            \end{picture}}
\put(1.25,0){\begin{picture}(1,0)
            \put(0,0.1){\circle*{0.075}}
            \put(1,0.1){\circle*{0.075}}
            \put(2,0.1){\circle*{0.075}}
            \put(0,0.1){\line(1,0){2}}
            \put(0.2,0.1){\vector(1,0){0.1}}
            \put(0.8,0.1){\vector(-1,0){0.1}}
            \put(0.7,0.1){\vector(-1,0){0.1}}
            \put(0.6,0.1){\vector(-1,0){0.1}}
            \end{picture}}
\put(4.25,0){\begin{picture}(1,0)
            \put(0,0.1){\circle*{0.075}}
            \put(1,0.1){\circle*{0.075}}
            \put(2,0.1){\circle*{0.075}}
            \put(0,0.1){\line(1,0){2}}
            \put(0.2,0.1){\vector(1,0){0.1}}
            \put(0.3,0.1){\vector(1,0){0.1}}
            \put(0.4,0.1){\vector(1,0){0.1}}
            \put(0.8,0.1){\vector(-1,0){0.1}}
            \end{picture}}
\end{picture}}

\newcommand{\FFourTwittle}{\setlength{\unitlength}{0.75in}
\begin{picture}(7.5,1.25)
\put(1.75,0.75){\begin{picture}(1,0)
            \put(0,0.1){\circle*{0.075}}
            \put(1,0.1){\circle*{0.075}}
            \put(2,0.1){\circle*{0.075}}
            \put(3,0.1){\circle*{0.075}}
            \put(4,0.1){\circle*{0.075}}
            \put(0,0.1){\line(1,0){4}}
            \put(1.2,0.1){\vector(1,0){0.1}}
            \put(1.3,0.1){\vector(1,0){0.1}}
            \put(1.8,0.1){\vector(-1,0){0.1}}
            \end{picture}}
\put(1.75,0){\begin{picture}(1,0)
            \put(0,0.1){\circle*{0.075}}
            \put(1,0.1){\circle*{0.075}}
            \put(2,0.1){\circle*{0.075}}
            \put(3,0.1){\circle*{0.075}}
            \put(4,0.1){\circle*{0.075}}
            \put(0,0.1){\line(1,0){4}}
            \put(1.2,0.1){\vector(1,0){0.1}}
            \put(1.8,0.1){\vector(-1,0){0.1}}
            \put(1.7,0.1){\vector(-1,0){0.1}}
            \end{picture}}
\end{picture}}

\newcommand{\SmallCycles}{\setlength{\unitlength}{0.75in}
\begin{picture}(6.5,3)
\put(-0.5,1.5){\begin{picture}(1,0)
            \put(0,0.6){\circle*{0.075}}
            \put(0.5,0.1){\circle*{0.075}}
            \put(0.5,1.1){\circle*{0.075}}
            \put(0,0.6){\line(1,1){0.5}}
            \put(0,0.6){\line(1,-1){0.5}}
            \put(0.5,0.1){\line(0,1){1}}
            \put(0.5,0.3){\vector(0,1){0.1}}
            \put(0.5,0.9){\vector(0,-1){0.1}}
            \put(0,0.6){\vector(1,1){0.2}}
            \put(0.5,1.1){\vector(-1,-1){0.2}}
            \put(0,0.6){\vector(1,-1){0.2}}
            \put(0.5,0.1){\vector(-1,1){0.2}}
            \put(-0.05,0.8){\footnotesize $q_{1}$}
            \put(0.175,1){\footnotesize $p_{1}$}
            \put(-0.05,0.35){\footnotesize $q_{2}$}
            \put(0.175,0.15){\footnotesize $p_{2}$} 
            \put(0.7,0.6){\footnotesize $p_{1}q_{1} \geq 1, 
            p_{2}q_{2} \geq 1$}
           \end{picture}}
\put(2.3,1.5){\begin{picture}(1,0)
            \put(0,0.6){\circle*{0.075}}
            \put(0.5,0.1){\circle*{0.075}}
            \put(0.5,1.1){\circle*{0.075}}
            \put(0,0.6){\line(1,1){0.5}}
            \put(0,0.6){\line(1,-1){0.5}}
            \put(0.5,0.1){\line(0,1){1}}
            \put(0.5,0.3){\vector(0,1){0.1}}
            \put(0.5,0.9){\vector(0,-1){0.1}}
            \put(0.5,0.8){\vector(0,-1){0.1}}
            \put(0,0.6){\vector(1,1){0.2}}
            \put(0.5,1.1){\vector(-1,-1){0.2}}
            \put(0,0.6){\vector(1,-1){0.2}}
            \put(0.5,0.1){\vector(-1,1){0.2}}
            \put(-0.05,0.8){\footnotesize $q_{1}$}
            \put(0.175,1){\footnotesize $p_{1}$}
            \put(-0.05,0.35){\footnotesize $q_{2}$}
            \put(0.175,0.15){\footnotesize $p_{2}$} 
            \put(0.7,0.6){\footnotesize $p_{1}q_{1} \geq 2, 
            p_{2}q_{2} \geq 2$}
           \end{picture}}
\put(5.15,1.5){\begin{picture}(1,0)
            \put(0,0.6){\circle*{0.075}}
            \put(0.5,0.1){\circle*{0.075}}
            \put(0.5,1.1){\circle*{0.075}}
            \put(0,0.6){\line(1,1){0.5}}
            \put(0,0.6){\line(1,-1){0.5}}
            \put(0.5,0.1){\line(0,1){1}}
            \put(0.5,0.3){\vector(0,1){0.1}}
            \put(0.5,0.9){\vector(0,-1){0.1}}
            \put(0.5,0.8){\vector(0,-1){0.1}}
            \put(0.5,0.7){\vector(0,-1){0.1}}
            \put(0,0.6){\vector(1,1){0.2}}
            \put(0.5,1.1){\vector(-1,-1){0.2}}
            \put(0,0.6){\vector(1,-1){0.2}}
            \put(0.5,0.1){\vector(-1,1){0.2}}
            \put(-0.05,0.8){\footnotesize $q_{1}$}
            \put(0.175,1){\footnotesize $p_{1}$}
            \put(-0.05,0.35){\footnotesize $q_{2}$}
            \put(0.175,0.15){\footnotesize $p_{2}$} 
            \put(0.7,0.6){\footnotesize $p_{1}q_{1} \geq 3$, 
            $p_{2}q_{2} \geq 3$}
           \end{picture}}
\put(0.5,0){\begin{picture}(1,0)
            \put(0,0.6){\circle*{0.075}}
            \put(0.5,0.1){\circle*{0.075}}
            \put(0.5,1.1){\circle*{0.075}}
            \put(1,0.6){\circle*{0.075}}
            \put(0,0.6){\line(1,1){0.5}}
            \put(0,0.6){\line(1,-1){0.5}}
            \put(1,0.6){\line(-1,1){0.5}}
            \put(1,0.6){\line(-1,-1){0.5}}
            \put(0.95,0.65){\vector(-1,1){0.1}}
            \put(0.55,1.05){\vector(1,-1){0.1}}
            \put(0.65,0.95){\vector(1,-1){0.1}}
           \end{picture}}
\put(2,0){\begin{picture}(1,0)
            \put(0,0.6){\circle*{0.075}}
            \put(0.5,0.1){\circle*{0.075}}
            \put(0.5,1.1){\circle*{0.075}}
            \put(1,0.6){\circle*{0.075}}
            \put(0,0.6){\line(1,1){0.5}}
            \put(0,0.6){\line(1,-1){0.5}}
            \put(1,0.6){\line(-1,1){0.5}}
            \put(1,0.6){\line(-1,-1){0.5}}
            \put(0.95,0.65){\vector(-1,1){0.1}}
            \put(0.55,1.05){\vector(1,-1){0.1}}
            \put(0.65,0.95){\vector(1,-1){0.1}}
            \put(0.45,0.15){\vector(-1,1){0.1}}
            \put(0.05,0.55){\vector(1,-1){0.1}}
            \put(0.15,0.45){\vector(1,-1){0.1}}
           \end{picture}}
\put(3.5,0){\begin{picture}(1,0)
            \put(0,0.6){\circle*{0.075}}
            \put(0.5,0.1){\circle*{0.075}}
            \put(0.5,1.1){\circle*{0.075}}
            \put(1,0.6){\circle*{0.075}}
            \put(0,0.6){\line(1,1){0.5}}
            \put(0,0.6){\line(1,-1){0.5}}
            \put(1,0.6){\line(-1,1){0.5}}
            \put(1,0.6){\line(-1,-1){0.5}}
            \put(0.95,0.65){\vector(-1,1){0.1}}
            \put(0.55,1.05){\vector(1,-1){0.1}}
            \put(0.65,0.95){\vector(1,-1){0.1}}
            \put(0.45,0.15){\vector(-1,1){0.1}}
            \put(0.05,0.55){\vector(1,-1){0.1}}
            \put(0.35,0.25){\vector(-1,1){0.1}}
           \end{picture}}
\put(5,0){\begin{picture}(1,0)
            \put(0,0.6){\circle*{0.075}}
            \put(0,0.1){\circle*{0.075}}
            \put(0.5,0.1){\circle*{0.075}}
            \put(0.5,1.1){\circle*{0.075}}
            \put(1,0.6){\circle*{0.075}}
            \put(0,0.6){\line(1,1){0.5}}
            \put(0,0.6){\line(0,-1){0.5}}
            \put(0,0.1){\line(1,0){0.5}}
            \put(1,0.6){\line(-1,1){0.5}}
            \put(1,0.6){\line(-1,-1){0.5}}
            \put(0.95,0.65){\vector(-1,1){0.1}}
            \put(0.55,1.05){\vector(1,-1){0.1}}
            \put(0.65,0.95){\vector(1,-1){0.1}}
           \end{picture}}
\end{picture}}

\newcommand{\CnTwittle}{\setlength{\unitlength}{0.75in}
\begin{picture}(6,2.25)
\put(0,0){\begin{picture}(1,0)
            \put(0,0.6){\circle*{0.075}}
            \put(0,0.1){\circle*{0.075}}
            \put(1,0.35){\circle*{0.075}}
            \put(2,0.35){\circle*{0.075}}
            \put(4,0.35){\circle*{0.075}}
            \put(5,0.35){\circle*{0.075}}
            \put(6,0.35){\circle*{0.075}}
            \put(0,0.1){\line(4,1){1}}
            \put(0,0.6){\line(4,-1){1}}
            \put(1,0.35){\line(1,0){1}}
            \multiput(2,0.35)(0.4,0){5}{\line(1,0){0.2}}
            \put(4,0.35){\line(1,0){2}}
            \put(5.2,0.35){\vector(1,0){0.1}}
            \put(5.8,0.35){\vector(-1,0){0.1}}
            \put(5.7,0.35){\vector(-1,0){0.1}}
           \end{picture}}
\put(0,1){\begin{picture}(1,0)
            \put(0,0.1){\circle*{0.075}}
            \put(1,0.1){\circle*{0.075}}
            \put(2,0.1){\circle*{0.075}}
            \put(4,0.1){\circle*{0.075}}
            \put(5,0.1){\circle*{0.075}}
            \put(6,0.1){\circle*{0.075}}
            \put(0,0.1){\line(1,0){2}}
            \multiput(2,0.1)(0.4,0){5}{\line(1,0){0.2}}
            \put(4,0.1){\line(1,0){2}}
            \put(0.2,0.1){\vector(1,0){0.1}}
            \put(0.8,0.1){\vector(-1,0){0.1}}
            \put(0.7,0.1){\vector(-1,0){0.1}}
            \put(5.8,0.1){\vector(-1,0){0.1}}
            \put(5.7,0.1){\vector(-1,0){0.1}}
            \put(5.2,0.1){\vector(1,0){0.1}}
           \end{picture}}
\put(0,1.75){\begin{picture}(1,0)
            \put(0,0.1){\circle*{0.075}}
            \put(1,0.1){\circle*{0.075}}
            \put(2,0.1){\circle*{0.075}}
            \put(4,0.1){\circle*{0.075}}
            \put(5,0.1){\circle*{0.075}}
            \put(6,0.1){\circle*{0.075}}
            \put(0,0.1){\line(1,0){2}}
            \multiput(2,0.1)(0.4,0){5}{\line(1,0){0.2}}
            \put(4,0.1){\line(1,0){2}}
            \put(0.2,0.1){\vector(1,0){0.1}}
            \put(0.3,0.1){\vector(1,0){0.1}}
            \put(0.8,0.1){\vector(-1,0){0.1}}
            \put(5.8,0.1){\vector(-1,0){0.1}}
            \put(5.7,0.1){\vector(-1,0){0.1}}
            \put(5.2,0.1){\vector(1,0){0.1}}
           \end{picture}}
\end{picture}}

\newcommand{\BnTwittle}{\setlength{\unitlength}{0.75in}
\begin{picture}(6,1.5)
\put(0,0){\begin{picture}(1,0)
            \put(0,0.1){\circle*{0.075}}
            \put(1,0.1){\circle*{0.075}}
            \put(2,0.1){\circle*{0.075}}
            \put(4,0.1){\circle*{0.075}}
            \put(5,0.1){\circle*{0.075}}
            \put(6,0.1){\circle*{0.075}}
            \put(0,0.1){\line(1,0){2}}
            \multiput(2,0.1)(0.4,0){5}{\line(1,0){0.2}}
            \put(4,0.1){\line(1,0){2}}
            \put(0.2,0.1){\vector(1,0){0.1}}
            \put(0.8,0.1){\vector(-1,0){0.1}}
            \put(0.7,0.1){\vector(-1,0){0.1}}
            \put(5.8,0.1){\vector(-1,0){0.1}}
            \put(5.2,0.1){\vector(1,0){0.1}}
            \put(5.3,0.1){\vector(1,0){0.1}}
           \end{picture}}
\put(0,0.5){\begin{picture}(1,0)
            \put(0,0.6){\circle*{0.075}}
            \put(0,0.1){\circle*{0.075}}
            \put(1,0.35){\circle*{0.075}}
            \put(2,0.35){\circle*{0.075}}
            \put(4,0.35){\circle*{0.075}}
            \put(5,0.35){\circle*{0.075}}
            \put(6,0.35){\circle*{0.075}}
            \put(0,0.1){\line(4,1){1}}
            \put(0,0.6){\line(4,-1){1}}
            \put(1,0.35){\line(1,0){1}}
            \multiput(2,0.35)(0.4,0){5}{\line(1,0){0.2}}
            \put(4,0.35){\line(1,0){2}}
            \put(5.8,0.35){\vector(-1,0){0.1}}
            \put(5.2,0.35){\vector(1,0){0.1}}
            \put(5.3,0.35){\vector(1,0){0.1}}
           \end{picture}}
\end{picture}}

\newcommand{\DnTwittle}{\setlength{\unitlength}{0.75in}
\begin{picture}(6,1)
\put(0,0){\begin{picture}(1,0)
            \put(0,0.6){\circle*{0.075}}
            \put(0,0.1){\circle*{0.075}}
            \put(1,0.35){\circle*{0.075}}
            \put(2,0.35){\circle*{0.075}}
            \put(4,0.35){\circle*{0.075}}
            \put(5,0.35){\circle*{0.075}}
            \put(6,0.1){\circle*{0.075}}
            \put(6,0.6){\circle*{0.075}}
            \put(0,0.1){\line(4,1){1}}
            \put(0,0.6){\line(4,-1){1}}
            \put(1,0.35){\line(1,0){1}}
            \multiput(2,0.35)(0.4,0){5}{\line(1,0){0.2}}
            \put(4,0.35){\line(1,0){1}}
            \put(5,0.35){\line(4,1){1}}
            \put(5,0.35){\line(4,-1){1}}
           \end{picture}}
\end{picture}}

\newcommand{\AnTwittle}{\setlength{\unitlength}{0.75in}
\begin{picture}(6,1)
\put(0,0){\begin{picture}(1,0)
            \put(0,0.1){\circle*{0.075}}
            \put(1,0.1){\circle*{0.075}}
            \put(2,0.1){\circle*{0.075}}
            \put(3,0.6){\circle*{0.075}}
            \put(4,0.1){\circle*{0.075}}
            \put(5,0.1){\circle*{0.075}}
            \put(6,0.1){\circle*{0.075}}
            \put(0,0.1){\line(1,0){2}}
            \put(0,0.1){\line(6,1){3}}
            \multiput(2,0.1)(0.4,0){5}{\line(1,0){0.2}}
            \put(4,0.1){\line(1,0){2}}
            \put(6,0.1){\line(-6,1){3}}
           \end{picture}}
\end{picture}}

\newcommand{\EnTwittle}{\setlength{\unitlength}{0.75in}
\begin{picture}(7,3.25)
\put(0,-0.25){\begin{picture}(1,0)
            \put(0,0.1){\circle*{0.075}}
            \put(1,0.1){\circle*{0.075}}
            \put(2,0.1){\circle*{0.075}}
            \put(2,0.6){\circle*{0.075}}
            \put(3,0.1){\circle*{0.075}}
            \put(4,0.1){\circle*{0.075}}
            \put(5,0.1){\circle*{0.075}}
            \put(6,0.1){\circle*{0.075}}
            \put(7,0.1){\circle*{0.075}}
            \put(0,0.1){\line(1,0){7}}
            \put(2,0.1){\line(0,1){0.5}}
           \end{picture}}
\put(0.5,0.75){\begin{picture}(1,0)
            \put(0,0.1){\circle*{0.075}}
            \put(1,0.1){\circle*{0.075}}
            \put(2,0.1){\circle*{0.075}}
            \put(3,0.1){\circle*{0.075}}
            \put(3,0.6){\circle*{0.075}}
            \put(4,0.1){\circle*{0.075}}
            \put(5,0.1){\circle*{0.075}}
            \put(6,0.1){\circle*{0.075}}
            \put(0,0.1){\line(1,0){6}}
            \put(3,0.1){\line(0,1){0.5}}
           \end{picture}}
\put(1.5,1.75){\begin{picture}(1,0)
            \put(0,0.1){\circle*{0.075}}
            \put(1,0.1){\circle*{0.075}}
            \put(2,0.1){\circle*{0.075}}
            \put(2,0.6){\circle*{0.075}}
            \put(2,1.1){\circle*{0.075}}
            \put(3,0.1){\circle*{0.075}}
            \put(4,0.1){\circle*{0.075}}
            \put(0,0.1){\line(1,0){4}}
            \put(2,0.1){\line(0,1){1}}
           \end{picture}}
\end{picture}}

\begin{document}

\newpage
\setcounter{page}{1} 
\renewcommand{\baselinestretch}{1}

\begin{center}
{\large \bf The numbers game, 
geometric representations of 
Coxeter groups,\\ and 
Dynkin diagram 
 classification results}

Robert G.\ Donnelly

Department of Mathematics and Statistics, Murray State
University, Murray, KY 42071
\end{center}

\begin{abstract}
The numbers game is a one-player game played on  
a 
finite simple 
graph with certain ``amplitudes'' assigned to its edges 
and with an initial assignment of real 
numbers to its nodes.  The moves of the game successively transform 
the numbers at the nodes using the amplitudes in a certain way.  
This game has been studied previously by 
Proctor, Mozes, Bj\"{o}rner, Eriksson, and Wildberger.  
We show that those connected 
such graphs for which the numbers game  
meets a certain finiteness requirement are 
precisely the Dynkin diagrams associated with the finite-dimensional 
complex simple Lie algebras.  As a consequence of our 
proof we obtain the classifications of the finite-dimensional 
Kac-Moody 
algebras and of the finite Weyl groups. 
We use Coxeter group theory to 
establish a more general result that applies 
to Eriksson's E-games: 
an E-game 
meets the finiteness requirement if and only if a naturally 
associated Coxeter group is finite. To prove this and some other 
finiteness results we 
further develop Eriksson's theory of a geometric representation of Coxeter 
groups and observe some curious differences of this representation 
from the standard 
geometric representation. 
\begin{center}

\ 


\end{center}
\end{abstract}
\def\abstractname{Contents}
\begin{abstract}
\begin{center}
\parbox{3.75in}{
\IntroNum. Introduction and statements of our main results\\
\FirstProofNum.  Proof of our first main result\\
\FirstProofCorollaryNum. Classifications of finite-dimensional\\ 
\hspace*{0.65in}Kac-Moody 
algebras and finite Weyl groups\\
\SecondProofNum. Quasi-standard geometric representations of\\ 
\hspace*{0.65in}Coxeter groups and a generalization to E-games\\ 
\ConvergenceNum. Some further finiteness aspects of E-game play\\ 
\RemarksNum. Comments 
}
\end{center}
\end{abstract}

\vspace{2ex} 

\noindent
{\Large \bf \IntroNum.\ \ Introduction and 
main results}

\vspace{2ex} 
The numbers game is a one-player game played on a finite simple graph 
with weights (which we call ``amplitudes'') on its edges 
and with an initial assignment of 
real numbers (which we call initial ``populations'') to its nodes.  
At the outset, each of the two edge amplitudes (one for each direction) 
will be negative integers; later
we will relax this integrality requirement. 
The 
move a player can make 
is to ``fire'' one of the nodes with a positive population.  This move 
transforms the population at the fired node 
by changing its sign, and it also 
transforms the population at each adjacent node in a certain way 
using an amplitude   
along the incident edge.  
The player fires the nodes in some sequence of 
the player's choosing, continuing until no node has a positive 
population.  
This numbers game formulated by Mozes \cite{Mozes} has also been 
studied by Proctor 
\cite{PrEur}, \cite{PrDComplete}, Bj\"{o}rner 
\cite{Bjorner}, \cite{BB}, 
Eriksson \cite{ErikssonLinear}, 
\cite{ErikssonDiscrete}, \cite{ErikssonEur}, and 
Wildberger 
\cite{WildbergerAdv}, \cite{WildbergerEur}, \cite{WildbergerPreprint}.  
Wildberger studies a dual version which 
he calls the ``mutation game.''  See Alon {\em et al} \cite{AKP} for a 
brief and readable treatment of the numbers game on ``unweighted'' 
cyclic graphs. Much of the numbers game discussion in Chapter 4 of 
\cite{BB} can be 
found in \cite{ErikssonDiscrete}.  
 Proctor developed this process in \cite{PrEur} 
to compute Weyl group 
orbits of weights with respect to the fundamental weight basis.  For 
this reason we prefer his perspective 
of firing 
nodes with positive, as opposed to negative, populations. 

The motivating question for this paper is: for which such graphs 
does there exist a nontrivial initial assignment of nonnegative 
populations such 
that the numbers game terminates in a finite number of steps?  
For graphs with integer amplitudes, 
our answer to this question (\MarsAttacksTheorem) 
is that the only such connected graphs are the Dynkin 
diagrams of \DynkinDiagramFigure.   
Moreover, from Eriksson's Strong Convergence Theorem (Theorem 3.1 of 
\cite{ErikssonEur}) we are able to conclude 
that for any 
initial assignment of 
populations 
to the nodes of a Dynkin diagram 
and for any legal sequence of node firings, the 
numbers game will terminate 
in the same finite number of steps and ultimately 
yield at each node the same nonpositive terminal 
population. 
Our proof of \MarsAttacksTheorem\ in Section \FirstProofNum\  
requires some Coxeter/Weyl group theory, but  
only implicitly --- in particular the proof of Eriksson's Comparison 
Theorem (Theorem 4.5 of 
\cite{ErikssonDiscrete}).   
As a consequence of our proof of \MarsAttacksTheorem\  
and with the help of another result of Eriksson 
we re-derive in Section \FirstProofCorollaryNum\ 
the classifications of 
the finite-dimensional Kac-Moody  
algebras 
(the finite-dimensional complex semisimple Lie algebras cf.\ 
\cite{Hum}, \cite{Kac}) 
and of the finite Weyl groups (the finite crystallographic 
Coxeter groups of \cite{HumCoxeter}).  

Our second main result (\EGCMTheorem)  
answers our motivating question for a class of 
graphs (the ``E-games'' of \cite{ErikssonDiscrete}) 
whose amplitudes are allowed to be certain real numbers. 
The classification obtained in this theorem 
uses 
the classification of 
finite irreducible Coxeter groups by 
connected 
positive definite Coxeter graphs (cf.\  
\S 2.3-2.7 in \cite{HumCoxeter}).  
The connection to Coxeter groups is made via a particular geometric 
representation studied in \cite{ErikssonDiscrete}, referred to here 
as a ``quasi-standard geometric representation.''  This has many similarities 
to the standard geometric representation, but also some surprising 
differences in regard to some finiteness properties.  Results we develop about 
quasi-standard geometric representations in Section \SecondProofNum\ lead to 
our proof of \EGCMTheorem\ and to some further E-game results in 
Section \ConvergenceNum.   These include a method for computing (in 
certain circumstances) the 
positive roots in the root system for a quasi-standard geometric 
representation of a finite Coxeter group (\AllPositiveRootsProp) 
and a classification (\AdjacencyFreeClassification) of those E-games 
for which the choices of node firings are ``interchangeable'' in some 
sense. 
The latter applies a classification result of Stembridge \cite{StemFC} 
about ``fully commutative'' elements in finite Coxeter groups. 
In Section \RemarksNum, we remark on other 
connections.

We formulate the problem statement and solution precisely as 
follows.  
Fix a positive integer $n$ and a totally ordered set $I_{n}$ with $n$ 
elements (usually $I_{n} := \{1<\ldots<n\}$).  
A {\em generalized Cartan matrix} (or {\em GCM}) is an $n \times n$ 
matrix $M = (M_{ij})_{i,j \in I_{n}}$  
with integer entries satisfying the requirements that each 
main diagonal matrix entry is 2, that all other matrix entries are 
nonpositive, and that if a matrix entry $M_{ij}$ is nonzero then its 
transpose entry $M_{ji}$ is also nonzero.  
Generalized Cartan matrices 
are the starting point for the study of 
Kac-Moody 
algebras: beginning with a GCM, one can 
write down a list of the defining relations for a Kac-Moody 
algebra as well as the associated Weyl group 
(see Section \FirstProofCorollaryNum).  
To an $n \times n$ generalized Cartan matrix 
$M = (M_{ij})_{i,j \in I_{n}}$ we associate a finite 
graph $\Gamma$ (which has undirected edges, 
no loops, and no multiple edges) 
as follows:     
The nodes $(\gamma_{i})_{i \in I_{n}}$ of $\Gamma$ are indexed 
by the set $I_{n}$, 
and   an edge is placed between nodes $\gamma_{i}$ and $\gamma_{j}$ 
if and only if $i \not= j$ 
and the matrix entries $M_{ij}$ and $M_{ji}$ are nonzero.  We call the pair 
$(\Gamma,M)$ a {\em GCM graph}. 
We consider two 
GCM graphs $(\Gamma, M = (M_{ij})_{i,j \in I_{n}})$ 
and $(\Gamma', M' = (M'_{pq})_{p,q \in I'_{n}})$ 
to be the same if under some bijection $\sigma: I_{n} \rightarrow 
I'_{n}$ we have nodes $\gamma_{i}$ and $\gamma_{j}$ in 
$\Gamma$ adjacent if and only if 
$\gamma'_{\sigma(i)}$ and $\gamma'_{\sigma(j)}$ are adjacent in $\Gamma'$ with 
$M_{ij} = M'_{\sigma(i),\sigma(j)}$. 
With $p = -M_{12}$ and $q = -M_{21}$, 
we depict a generic connected two-node GCM graph as follows:

\vspace*{-0.075in}
\noindent
\begin{center}
\TwoCitiesGraphWithLabels
\end{center}

\vspace*{-0.075in}
\noindent 
We use special names and notation to refer to 
two-node GCM 
graphs which have $p = 1$ and $q = 1$, $2$, or $3$ respectively:
\noindent
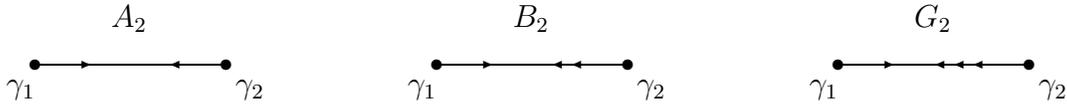
\begin{center}
\ATwoGraphWithLabels
\hspace*{1in}
\BTwoGraphWithLabels
\hspace*{1in}
\GTwoGraphWithLabels
\end{center}
When $p=1$ and $q=1$ it is convenient to use the graph 
\ATwoGraphNoEdgeLabels\ to 
represent the GCM graph $A_{2}$.  
A GCM graph $(\Gamma,M)$ is a {\em Dynkin diagram} if each connected 
component of $(\Gamma,M)$ 
is one of the graphs of \DynkinDiagramFigure.  We number our nodes as 
in \S 11.4 of \cite{Hum}. In these cases the GCMs  are ``Cartan'' 
matrices. 

\begin{figure}[t]
\begin{center}
\DynkinDiagramFigure: Connected Dynkin diagrams. 
\end{center}

\vspace*{-0.45in}
\begin{tabular}{cl}
$A_{n}$ ($n \geq 1$) & \AnMarsOK\\

$B_{n}$ ($n \geq 2$) & \BnMarsOK\\

$C_{n}$ ($n \geq 3$) & \CnMarsOK\\

$D_{n}$ ($n \geq 4$) & \DnMarsOK\\

$E_{6}$ & \ESixMarsOK\\

$E_{7}$ & \ESevenMarsOK\\

$E_{8}$ & \EEightMarsOK\\

$F_{4}$ & \FFourMarsOK\\

$G_{2}$ & \GTwoMarsOK
\end{tabular}
\end{figure}

A {\em position} $\lambda = (\lambda_i)_{i \in 
I_{n}}$ is an assignment of real numbers to the nodes of the GCM graph 
$(\Gamma,M)$; the real number $\lambda_i$ is the {\em population 
at node} $\gamma_{i}$.  The position $\lambda$ is 
{\em dominant} (respectively, {\em strongly dominant}) if 
$\lambda_{i} \geq 0$ 
(resp. $\lambda_i > 0$) for all $i \in I_{n}$; 
$\lambda$ is {\em nonzero} if at least one $\lambda_i \not= 0$. 
For $i \in I_{n}$, the {\em 
fundamental position} $\omega_i$ is the assignment of population 
$1$ at node $\gamma_{i}$ and population $0$ at all other nodes.  
Given a position $\lambda$ for a GCM graph $(\Gamma,M)$, to 
{\em fire} a node $\gamma_{i}$ is to change the population at each node 
$\gamma_{j}$ of $\Gamma$ by the transformation  
\[\lambda_j	 \longmapsto \lambda_j - 
M_{ij}\lambda_i,\] provided the population at node 
$\gamma_{i}$ is positive; otherwise node $\gamma_{i}$ is not allowed 
to be fired. 
Since the generalized Cartan  
matrix $M$ assigns a pair of {\em amplitudes} ($M_{ij}$ and 
$M_{ji}$) to each edge of the 
graph $\Gamma$, we sometimes refer to GCMs as 
{\em amplitude matrices}.  
The {\em numbers game} 
is the one-player game on a GCM graph $(\Gamma,M)$ in which the player 
(1) Assigns an initial position  
to the nodes of $\Gamma$; (2) Chooses a node with a positive 
population and fires the node to obtain a new position; and (3) 
Repeats step (2) for the new position if there is at least one node 
with a positive population.\myfootnote{Mozes 
studied numbers games on GCM graphs for 
which the amplitude matrix $M$ is {\em symmetrizable} 
(i.e.\ there is a nonsingular  
diagonal matrix $D$ such that $D^{-1}M$ is symmetric); 
in \cite{Mozes} he obtained strong convergence results and a 
geometric characterization of the initial positions for which the 
game terminates.}  

Consider now the GCM graph $B_{2}$.  As we can see in \BtwoMarsDemo, 
the numbers game terminates in a finite number of steps for any 
initial position and any legal sequence of node firings, 
if it is understood that the player  
will continue to fire as long as there is at least one 
node with a positive population.  In general, 
given a position $\lambda$, a {\em game sequence 
for} $\lambda$ is the (possibly empty, possibly infinite) sequence 
$(\gamma_{i_{1}}, \gamma_{i_{2}},\ldots)$, where 
$\gamma_{i_{j}}$ 
is the $j$th node that is fired in some 
numbers game with initial position $\lambda$.  
More generally, a {\em firing sequence} from some position $\lambda$ is an 
initial portion of some game sequence played from $\lambda$; the 
phrase {\em legal firing sequence} is used to emphasize that all node 
firings in the sequence are known or assumed to be possible. 
Note that a game sequence 
$(\gamma_{i_{1}}, 
\gamma_{i_{2}},\ldots,\gamma_{i_{l}})$ 
is of finite length $l$ 
(possibly with $l = 0$) if 
the population is nonpositive at each node after the $l$th firing; in 
this case we say the game sequence is {\em convergent} and the 
resulting position is the {\em terminal position}.  
We say a connected GCM graph $(\Gamma,M)$ is {\em admissible} if 
there exists a nonzero dominant initial position with a convergent 
game sequence. 
Our first main result (to be proved in Section 
\FirstProofNum) is: 

\noindent 
{\bf \MarsAttacksTheorem}\ \ {\sl A connected GCM graph $(\Gamma,M)$ is 
admissible if and only if it 
is a connected Dynkin diagram.  
In these cases, for any given initial position every game sequence  
will converge to the same terminal position in the 
same finite number of steps.}

\begin{figure}[htb]
\begin{center}
\BtwoMarsDemo: The numbers game for the GCM graph $B_{2}$. 

\vspace*{0.25in}
\setlength{\unitlength}{1in}
\begin{picture}(1.5,5.6) 
\put(0.25,5.2){\BTwoGraphForFigure}
\put(0.25,5.4){$a$}
\put(1.3,5.4){$b$}
\put(0.3,5){\vector(-4,-3){0.7}}
\put(1.3,5){\vector(4,-3){0.7}}
\put(1.75,3.9){\BTwoGraphForFigure}
\put(1.5,4.1){$a+2b$}
\put(2.7,4.1){$-b$}
\put(2.3,3.7){\vector(0,-1){0.6}}
\put(1.75,2.6){\BTwoGraphForFigure}
\put(1.45,2.8){$-a-2b$}
\put(2.65,2.8){$a+b$}
\put(2.3,2.4){\vector(0,-1){0.6}}
\put(1.75,1.3){\BTwoGraphForFigure}
\put(1.75,1.5){$a$}
\put(2.6,1.5){$-a-b$}
\put(2.1,1.1){\vector(-4,-3){0.7}}
\put(-1.25,3.9){\BTwoGraphForFigure}
\put(-1.35,4.1){$-a$}
\put(-0.35,4.1){$a+b$}
\put(-0.7,3.7){\vector(0,-1){0.6}}
\put(-1.25,2.6){\BTwoGraphForFigure}
\put(-1.5,2.8){$a+2b$}
\put(-0.4,2.8){$-a-b$}
\put(-0.7,2.4){\vector(0,-1){0.6}}
\put(-1.25,1.3){\BTwoGraphForFigure}
\put(-1.55,1.5){$-a-2b$}
\put(-0.2,1.5){$b$}
\put(-0.5,1.1){\vector(4,-3){0.7}}
\put(0.25,0){\BTwoGraphForFigure}
\put(0.15,0.2){$-a$}
\put(1.2,0.2){$-b$}
\end{picture}
\end{center}
\end{figure}

In \cite{ErikssonLinear}, Eriksson proves the following related result 
using combinatorial reasoning and a result from the Perron-Frobenius 
theory for eigenvalues of nonnegative real matrices: For a connected GCM graph 
$(\Gamma,M)$ whose amplitude products are unity, every initial 
position has a convergent game sequence if and only if $(\Gamma,M)$ is 
one of $A_{n}$, $D_{n}$, $E_{6}$, $E_{7}$, or $E_{8}$ from 
\DynkinDiagramFigure\ (and if and only if the ``dominant 
eigenvalue'' of $M$ is less that two).  
Wildberger generalizes this assertion to all 
connected GCM graphs (see \cite{WildbergerPreprint}); the resulting 
GCM graphs are the connected Dynkin diagrams.  His proof also uses the 
Perron-Frobenius theory and, in particular, does not depend on the 
classification of finite Weyl groups.  In the language of this paper, 
say a GCM graph is {\em strongly admissible} if every nonzero 
dominant position has a convergent game sequence.  Then 
\cite{ErikssonLinear} gives an ``ADE'' version of the following 
result: 

\noindent 
{\bf \WildbergerTheorem\ (Wildberger)}\ \ {\sl A connected GCM graph 
is strongly admissible if and only if it is a connected Dynkin diagram.}

\MarsAttacksTheorem\ does not require ``strongly,'' and its proof does 
not refer to eigenvalues.  
We will also demonstrate a more general 
version of \MarsAttacksTheorem\ that 
applies to what Eriksson calls ``E-games'' in 
\cite{ErikssonDiscrete}.\myfootnote{Eriksson uses ``E'' for edge;  
he also allows for ``N-games'' 
where, in addition, nodes can be weighted.}
There Eriksson drops the integrality 
requirement for off-diagonal entries of the amplitude 
matrix 
and asks: For which such graphs will the numbers game be 
``strongly convergent''? 
(This concept is defined below at the beginning of Section \FirstProofNum.) 
His answer (see Theorem 2.2 of \cite{ErikssonDiscrete} or Theorem 3.1 of 
\cite{ErikssonEur}) is that the amplitude matrix 
must be what we will call here an 
{\em E-generalized Cartan matrix}  
or {\em E-GCM}: This is  
an $n \times n$ matrix $M = (M_{ij})_{i,j \in I_{n}}$  
with real entries satisfying the requirements that each 
main diagonal matrix entry is 2, that all other matrix entries are 
nonpositive, that if a matrix entry $M_{ij}$ is nonzero then its 
transpose entry $M_{ji}$ is also nonzero, and that if 
$M_{ij}M_{ji}$ is nonzero then $M_{ij}M_{ji} \geq 4$ or 
$M_{ij}M_{ji} = 4\cos^{2}(\pi/m_{ij})$ for some integer $m_{ij} \geq 
3$.  An {\em E-GCM graph} is the pair $(\Gamma,M)$ for an 
E-generalized Cartan matrix $M$.  
As before, we depict a generic two-node E-GCM graph as follows:  

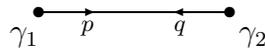
\vspace*{-0.075in}
\noindent
\begin{center}
\TwoCitiesGraphWithLabels
\end{center}

\vspace*{-0.075in}
\noindent
In this graph, $p = -M_{12}$  
and $q = -M_{21}$.  
We use \EGCMGraphCirclem 
for the collection of all two-node E-GCM 
graphs for which $M_{12}M_{21} = pq = 4\cos^{2}(\pi/m)$ 
for an integer $m > 3$;  when $m = 3$ (i.e.\ $pq = 1$), 
we use an unlabelled edge \ATwoGraphNoEdgeLabels\ as before.  
An {\em E-Coxeter graph} will be any E-GCM graph whose 
connected components come from one of the collections of 
\ECoxeterGraphFigure. 
Other terminology of this section used for GCM graphs will also be 
used in reference to E-GCM graphs and their numbers games ({\em 
firing}, {\em position}, {\em admissible}, etc).  
Our second main result generalizes \MarsAttacksTheorem. 

\begin{figure}[t]
\begin{center}
\ECoxeterGraphFigure: Families of connected E-Coxeter graphs. 
\vspace*{0.04in}

{\footnotesize (For adjacent nodes, the notation$\!\!$ 
\CircleIntegerm$\!\!$ means 
that the amplitude product on the edge is $4\cos^{2}(\pi/m)$;\\ for an 
unlabelled edge take $m=3$.  The asterisks for $\mathcal{E}_{6}$, 
$\mathcal{E}_{7}$, and $\mathcal{H}_{3}$ pertain to 
\AdjacencyFreeClassification.)}
\end{center}

\vspace*{-0.45in}
\begin{tabular}{cl}
$\mathcal{A}_{n}$ ($n \geq 1$) & \AnEGraph\\

$\mathcal{B}_{n}$ ($n \geq 3$) & \BnEGraph\\

$\mathcal{D}_{n}$ ($n \geq 4$) & \DnEGraph\\

$\mathcal{E}_{6}$ & \ESixEGraph\\

$\mathcal{E}_{7}$ & \ESevenEGraph\\

$\mathcal{E}_{8}$ & \EEightEGraph\\

$\mathcal{F}_{4}$ & \FFourEGraph\\

$\mathcal{H}_{3}$ & \HThreeEGraph\\

$\mathcal{H}_{4}$ & \HFourEGraph\\

$\mathcal{I}_{2}^{(m)}$ ($m \geq 4$) & \ITwoEGraph
\end{tabular}
\end{figure}

\noindent 
{\bf \EGCMTheorem}\ \ {\sl A connected E-GCM graph $(\Gamma,M)$ is 
admissible if and only if it 
is a connected E-Coxeter graph.  
In these cases, for any given initial position every game sequence  
will converge to the same terminal position in the 
same finite number of steps.} 

\noindent 
{\bf Acknowledgments}\ \ 
We thank  
John Eveland for stimulating discussions during his work on 
an undergraduate research project \cite{Eveland} at Murray 
State University; this helped lead to the motivating question we 
address here.  
We thank 
Norman Wildberger for sharing his perspective on the numbers game 
(including his observation about the appearance of 
``positive root functionals'' in numbers games on GCM graphs); 
this helped us formulate the question as well as 
our proof of our first main result.  
We thank Bob Proctor for pointing us in the 
direction of Eriksson's work and for suggesting that the 
classification results of \ProofCorollary\ follow from the proof of 
\MarsAttacksTheorem. 
We also thank Kimmo Eriksson for his helpful feedback.

\vspace{2ex} 

\noindent
{\Large \bf \FirstProofNum.\ \ Proof of our first main result}

\vspace{2ex} 
Our proof of the ``only if'' direction of the first claim of 
\MarsAttacksTheorem\ uses a series of reductions that are typical in 
Dynkin diagram classification arguments.  These reductions are 
implemented using several results of Eriksson.  This helps us to  
minimize the use of Coxeter group 
theory. 

\noindent 
{\em Proof of the ``only if'' direction of the first claim of 
\MarsAttacksTheorem:} 

\underline{Step 1: Strong convergence.}  
Following \cite{ErikssonEur}, we say 
the numbers game on a GCM graph $(\Gamma,M)$ is {\em strongly 
convergent} if given any initial position, every game sequence 
either diverges or converges to the same terminal position in the 
same number of steps. The next result  
follows from 
Theorem 3.1 of \cite{ErikssonEur}.  

\noindent
{\bf \StrongConvergenceTheorem\ (Eriksson's Strong Convergence Theorem)}\ \ 
{\sl The numbers game on a connected GCM graph 
is strongly 
convergent.}  

For this part of our proof of \MarsAttacksTheorem, 
we only require the 
following weaker result, which also applies when the 
GCM graph is not connected:  

\noindent
{\bf \StrongConvergenceCorollary} \ \ {\sl In any GCM graph, 
if a game sequence for an initial position  
$\lambda$ diverges, then all game  
sequences for $\lambda$ diverge.}
 
\underline{Step 2: Comparison.}  
The next result is an immediate consequence of Theorem 4.5 of 
\cite{ErikssonDiscrete}.  
Eriksson's proof of this result uses some Coxeter group theory. 
 
\noindent 
{\bf \ComparisonTheorem\ (Eriksson's Comparison Theorem)}\ \ 
{\sl Given a GCM graph, suppose that a game sequence 
for an initial position $\lambda = (\lambda_{i})_{i \in 
I_{n}}$ converges.  Suppose that a position $\lambda' := (\lambda'_{i})_{i 
\in I_{n}}$ has the property that $\lambda'_{i} \leq 
\lambda_{i}$ for all $i \in I_{n}$.  Then some game sequence 
for the initial position $\lambda'$ also converges.}

\noindent 
{\bf \ComparisonLemma}\ \ {\sl Let $r$ be a positive real number.  
If 
$(\gamma_{i_{1}},\ldots,\gamma_{i_{l}})$ 
is a convergent 
game sequence for an initial position $\lambda = 
(\lambda_{i})_{i \in I_{n}}$, then 
$(\gamma_{i_{1}},\ldots,\gamma_{i_{l}})$ 
is a convergent 
game sequence for the initial position $r\lambda := 
(r\lambda_{i})_{i \in I_{n}}$.}

{\em Proof.}  For game play from initial position $\lambda$, let 
$c_{j}$ be the positive population at node $\gamma_{i_{j}}$ when that 
node is fired in the game sequence 
$(\gamma_{i_{1}},\ldots,\gamma_{i_{l}})$.  Now $rc_{j}$ will be the 
positive population at the same node when the same game sequence is 
applied to the initial position $r\lambda$.\hfill\QED

The next result follows immediately from  
\ComparisonTheorem\ together with \ComparisonLemma: 

\noindent 
{\bf \ComparisonCorollary}\ \ {\sl 
Suppose $(\Gamma,M)$ is admissible. Let $\lambda = (\lambda_{i})_{i \in 
I_{n}}$ be a dominant initial position such that $\lambda_{j} > 0$ 
for some $j \in I_{n}$. Suppose that a game sequence for $\lambda$ 
converges.  
Then some game   
sequence for the fundamental position $\omega_{j}$ also 
converges.}

\underline{Step 3: A catalog of connected GCM graphs that are not 
admissible.} 
The following immediate consequence of \ComparisonResults\ is useful 
in the proof of \NotMarsFriendlyCatalog: 

\noindent 
{\bf \NotMarsFriendlyLemma}\ \ {\sl A GCM graph is not admissible if and 
only if for each fundamental position there is a divergent game  
sequence.} 

\noindent 
{\bf \NotMarsFriendlyCatalog}\ \ {\sl The connected GCM graphs of 
\NotMarsFriendlyFigure\myfootnote{The infinite ``$\widetilde{A}$'' 
family of GCM graphs of \NotMarsFriendlyFigure\  
is the family of cycles with 
amplitude products of unity on all edges.  
Such cycles were in fact the graphs that motivated Mozes' study of the 
numbers game in \cite{Mozes}.} are not admissible.} 

\begin{figure}[ht]
\begin{center}
\NotMarsFriendlyFigure: Some connected GCM graphs that are not 
admissible. 

\ 

The ``$\widetilde{A}$'' family of GCM graphs\\
\AnTwittle 

\vspace*{0.5in}
The ``$\widetilde{B}$'' family of GCM graphs\\
\BnTwittle 

\vspace*{0.5in}
The ``$\widetilde{C}$'' family of GCM graphs\\
\CnTwittle 

\vspace*{0.5in}
The ``$\widetilde{D}$'' family of GCM graphs\\
\DnTwittle 
\end{center}
\end{figure} 

\begin{figure}[ht]
\begin{center}
\NotMarsFriendlyFigure\ (continued): Some connected GCM graphs that are not 
admissible. 

\ 

The ``$\widetilde{E}$'' family of GCM graphs\\
\EnTwittle 

\vspace*{0.75in}
The ``$\widetilde{F}$'' family of GCM graphs\\
\FFourTwittle 

\vspace*{0.5in}
The ``$\widetilde{G}$'' family of GCM graphs\\
\GTwoTwittle 
\end{center}
\end{figure}

\begin{figure}[ht]
\begin{center}
\NotMarsFriendlyFigure\ (continued): Some connected GCM graphs that are not 
admissible. 

\ 

Families of small cycles\\
\SmallCycles  
\end{center}
\end{figure}

{\em Sketch of proof.} 
By \NotMarsFriendlyLemma, it 
suffices to show that for each graph in 
\NotMarsFriendlyFigure\ one game 
sequence for each fundamental position diverges. 
Our approach is to start with a 
fundamental position assigned to one of these graphs 
and then proceed with a numbers 
game, firing in some predictable manner until a pattern emerges. 
This exercise can be completed by hand in a reasonable amount of time. 
We illustrate the nature of our arguments with two examples.  

First, we show why GCM graphs in the ``$\widetilde{D}$'' family are 
not admissible.  The smallest GCM graph in this family has five 
nodes.  \DtwittleFiveNodesFigure\ 
shows that for each fundamental position there exists a game sequence 
that repeats indefinitely.  In \DtwittleFiveNodesFigure.1, we fired each 
node once.  
In finding sequences of node firings in 
\DtwittleFiveNodesFigure.2, we analyzed the $k=0$ and $k>0$ cases 
separately; for the latter we fired the center node twice.  
In the case of more than five nodes,  
\DtwittleFigure\ 
shows that for each fundamental position there exists a game  
sequence that repeats indefinitely. 
In \DtwittleFigure.1, we started by firing at the 
node with positive population, then fired all nodes to the right in 
succession (and returned), then fired all nodes to the left in 
succession (and returned); we only fired the starting node once. 
In \DtwittleFigure.2, we fired first at the positive population node, 
then proceeded to fire the nodes to its right in succession (and 
returned), then fired the starting node again, and finally fired the 
two leftmost nodes. 
In \DtwittleFigure.3, we traversed the ``isthmus'' of the graph twice. 

\begin{figure}[ht]
\begin{center}
\DtwittleFiveNodesFigure: The five-node  
GCM graph in the ``$\widetilde{D}$'' family is 
not admissible. \\ 
\ \\ 
\vspace*{-0.125in}

{\footnotesize (In each part of this figure, for any $k \geq 0$, 
the position on the 
right can be obtained 
 
from the position on the left by a sequence of 
legal node firings.)}

\ 

\setlength{\unitlength}{0.75in}
\begin{picture}(6,0.85)
\put(-1,0.3){2.2.1}
\put(0,0){\begin{picture}(1,0)
 \put(-0.1,0.675){\footnotesize $-k$}
 \put(-0.1,-0.1){\footnotesize $-k$}
 \put(1.9,0.675){\footnotesize $-k$}
 \put(1.9,-0.1){\footnotesize $-k$}
 \put(0.75,0.45){\footnotesize $2k+1$}
            \put(0,0.6){\circle*{0.075}}
            \put(0,0.1){\circle*{0.075}}
            \put(1,0.35){\circle*{0.075}}
            \put(2,0.1){\circle*{0.075}}
            \put(2,0.6){\circle*{0.075}}
            \put(0,0.1){\line(4,1){1}}
            \put(0,0.6){\line(4,-1){1}}
            \put(1,0.35){\line(4,1){1}}
            \put(1,0.35){\line(4,-1){1}}
           \end{picture}}
\put(2.5,0.35){\vector(1,0){1}}
\put(4,0){\begin{picture}(1,0)
 \put(-0.25,0.675){\footnotesize $-k-1$}
 \put(-0.25,-0.1){\footnotesize $-k-1$}
 \put(1.75,0.675){\footnotesize $-k-1$}
 \put(1.75,-0.1){\footnotesize $-k-1$}
 \put(0.75,0.45){\footnotesize $2k+3$}
            \put(0,0.6){\circle*{0.075}}
            \put(0,0.1){\circle*{0.075}}
            \put(1,0.35){\circle*{0.075}}
            \put(2,0.1){\circle*{0.075}}
            \put(2,0.6){\circle*{0.075}}
            \put(0,0.1){\line(4,1){1}}
            \put(0,0.6){\line(4,-1){1}}
            \put(1,0.35){\line(4,1){1}}
            \put(1,0.35){\line(4,-1){1}}
           \end{picture}}
\end{picture}

\ 

\begin{picture}(6,1)
\put(-1,0.3){2.2.2}
\put(0,0){\begin{picture}(1,0)
 \put(-0.2,0.675){\footnotesize $k+1$}
 \put(-0.05,-0.1){\footnotesize $k$}
 \put(1.95,0.675){\footnotesize $k$}
 \put(1.95,-0.1){\footnotesize $k$}
 \put(0.85,0.45){\footnotesize $-2k$}
            \put(0,0.6){\circle*{0.075}}
            \put(0,0.1){\circle*{0.075}}
            \put(1,0.35){\circle*{0.075}}
            \put(2,0.1){\circle*{0.075}}
            \put(2,0.6){\circle*{0.075}}
            \put(0,0.1){\line(4,1){1}}
            \put(0,0.6){\line(4,-1){1}}
            \put(1,0.35){\line(4,1){1}}
            \put(1,0.35){\line(4,-1){1}}
           \end{picture}}
\put(2.5,0.35){\vector(1,0){1}}
\put(4,0){\begin{picture}(1,0)
 \put(-0.2,0.675){\footnotesize $k+2$}
 \put(-0.2,-0.1){\footnotesize $k+1$}
 \put(1.8,0.675){\footnotesize $k+1$}
 \put(1.8,-0.1){\footnotesize $k+1$}
 \put(0.7,0.45){\footnotesize $-2k-2$}
            \put(0,0.6){\circle*{0.075}}
            \put(0,0.1){\circle*{0.075}}
            \put(1,0.35){\circle*{0.075}}
            \put(2,0.1){\circle*{0.075}}
            \put(2,0.6){\circle*{0.075}}
            \put(0,0.1){\line(4,1){1}}
            \put(0,0.6){\line(4,-1){1}}
            \put(1,0.35){\line(4,1){1}}
            \put(1,0.35){\line(4,-1){1}}
           \end{picture}}
\end{picture}
\end{center}
\end{figure}

\begin{figure}[ht]
\begin{center}
\DtwittleFigure:  
GCM graphs in the ``$\widetilde{D}$'' family with more than five nodes are  
not admissible. \\ 
\ \\ 
\vspace*{-0.125in}

{\footnotesize (In each part of this figure, for any $k \geq 0$, the 
subsequent position 
can be obtained 

from the first position 
by a sequence of legal node firings.)}

\ 

\setlength{\unitlength}{0.75in}
\begin{picture}(7,2.0)
\put(-0.85,1.55){2.3.1}
\put(0,1.25){\begin{picture}(1,0)
 \put(-0.085,0.675){\footnotesize $0$}
 \put(-0.085,-0.1){\footnotesize $0$}
 \put(0.975,0.45){\footnotesize $0$}
 \put(2.4,0.45){\footnotesize $-k$}
 \put(3.25,0.45){\footnotesize $2k+1$}
 \put(4.4,0.45){\footnotesize $-k$}
 \put(5.95,0.45){\footnotesize $0$}
 \put(7.0,0.675){\footnotesize $0$}
 \put(7.0,-0.1){\footnotesize $0$}
            \put(0,0.6){\circle*{0.075}}
            \put(0,0.1){\circle*{0.075}}
            \put(1,0.35){\circle*{0.075}}
            \put(2.5,0.35){\circle*{0.075}}
            \put(3.5,0.35){\circle*{0.075}}
            \put(4.5,0.35){\circle*{0.075}}
            \put(6,0.35){\circle*{0.075}}
            \put(7,0.1){\circle*{0.075}}
            \put(7,0.6){\circle*{0.075}}
            \put(0,0.1){\line(4,1){1}}
            \put(0,0.6){\line(4,-1){1}}
            \put(2.5,0.35){\line(1,0){2}}
            \multiput(1,0.35)(0.4,0){4}{\line(1,0){0.2}}
            \put(6,0.35){\line(4,1){1}}
            \put(6,0.35){\line(4,-1){1}}
            \multiput(4.5,0.35)(0.4,0){4}{\line(1,0){0.2}}
           \end{picture}}
\put(-0.6,0.35){\vector(1,0){1}}
\put(0.65,0){\begin{picture}(1,0)
 \put(-0.085,0.675){\footnotesize $0$}
 \put(-0.085,-0.1){\footnotesize $0$}
 \put(0.975,0.45){\footnotesize $0$}
 \put(2.25,0.45){\footnotesize $-k-1$}
 \put(3.25,0.45){\footnotesize $2k+3$}
 \put(4.25,0.45){\footnotesize $-k-1$}
 \put(5.95,0.45){\footnotesize $0$}
 \put(7.0,0.675){\footnotesize $0$}
 \put(7.0,-0.1){\footnotesize $0$}
            \put(0,0.6){\circle*{0.075}}
            \put(0,0.1){\circle*{0.075}}
            \put(1,0.35){\circle*{0.075}}
            \put(2.5,0.35){\circle*{0.075}}
            \put(3.5,0.35){\circle*{0.075}}
            \put(4.5,0.35){\circle*{0.075}}
            \put(6,0.35){\circle*{0.075}}
            \put(7,0.1){\circle*{0.075}}
            \put(7,0.6){\circle*{0.075}}
            \put(0,0.1){\line(4,1){1}}
            \put(0,0.6){\line(4,-1){1}}
            \put(2.5,0.35){\line(1,0){2}}
            \multiput(1,0.35)(0.4,0){4}{\line(1,0){0.2}}
            \put(6,0.35){\line(4,1){1}}
            \put(6,0.35){\line(4,-1){1}}
            \multiput(4.5,0.35)(0.4,0){4}{\line(1,0){0.2}}
           \end{picture}}
\end{picture}
\end{center}
\end{figure}

Second, we show why GCM graphs of the form \hspace*{0.1in}
\parbox[c]{0.5in}{
\setlength{\unitlength}{0.75in}
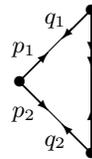
\begin{picture}(0.6,1.2)
            \put(0,0.6){\circle*{0.075}}
            \put(0.5,0.1){\circle*{0.075}}
            \put(0.5,1.1){\circle*{0.075}}
            \put(0,0.6){\line(1,1){0.5}}
            \put(0,0.6){\line(1,-1){0.5}}
            \put(0.5,0.1){\line(0,1){1}}
            \put(0.5,0.3){\vector(0,1){0.1}}
            \put(0.5,0.9){\vector(0,-1){0.1}}
            \put(0.5,0.8){\vector(0,-1){0.1}}
            \put(0,0.6){\vector(1,1){0.2}}
            \put(-0.05,0.8){\footnotesize $p_{1}$}
            \put(0.175,1){\footnotesize $q_{1}$}
            \put(0.5,1.1){\vector(-1,-1){0.2}}
            \put(0,0.6){\vector(1,-1){0.2}}
            \put(0.5,0.1){\vector(-1,1){0.2}}
            \put(-0.05,0.35){\footnotesize $p_{2}$}
            \put(0.175,0.15){\footnotesize $q_{2}$} 
\end{picture}
}
are not admissible.  We assume that the amplitude products 
$p_{1}q_{1}$ and $p_{2}q_{2}$ are at least two 
(at Step 6 below we 
will see that it suffices 
in \NotMarsFriendlyFigure\ 
to assume these products are at most 3). 
Assign populations $a$, $b$, 
and $c$ as follows:  
\hspace*{0.1in}
\parbox[c]{0.5in}{
\setlength{\unitlength}{0.75in}
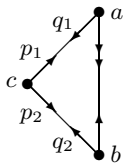
\begin{picture}(0.6,1.2)
            \put(0,0.6){\circle*{0.075}}
            \put(0.5,0.1){\circle*{0.075}}
            \put(0.5,1.1){\circle*{0.075}}
            \put(0,0.6){\line(1,1){0.5}}
            \put(0,0.6){\line(1,-1){0.5}}
            \put(0.5,0.1){\line(0,1){1}}
            \put(0.5,0.3){\vector(0,1){0.1}}
            \put(0.5,0.9){\vector(0,-1){0.1}}
            \put(0.5,0.8){\vector(0,-1){0.1}}
            \put(0,0.6){\vector(1,1){0.2}}
            \put(-0.05,0.8){\footnotesize $p_{1}$}
            \put(0.175,1){\footnotesize $q_{1}$}
            \put(0.5,1.1){\vector(-1,-1){0.2}}
            \put(0,0.6){\vector(1,-1){0.2}}
            \put(0.5,0.1){\vector(-1,1){0.2}}
            \put(-0.05,0.35){\footnotesize $p_{2}$}
            \put(0.175,0.15){\footnotesize $q_{2}$} 
            \put(-0.15,0.58){\footnotesize $c$}
            \put(0.58,0.05){\footnotesize $b$}
            \put(0.58,1.05){\footnotesize $a$}
\end{picture}
}
Set $\kappa := (2p_{1}+2p_{2}-\frac{1}{q_{1}})a + 
(p_{1}+2p_{2}-\frac{1}{q_{2}})b + c$.  
Assume for now that $a \geq 0$, $b \geq 0$, $c \leq 0$, and $\kappa > 0$; 
when these inequalities hold we will say the position $(a,b,c)$ meets 
condition ({\tt *}). Under 
condition ({\tt *}) notice that $a$ and $b$ cannot both be zero. 
Begin by firing only at the two 
rightmost nodes.  When this is no longer possible, fire at the 
leftmost node.  
The resulting corresponding 
populations are $a_{1} = 
q_{1}(\kappa + \frac{1}{q_{2}}b)$, $b_{1} = q_{2}(\kappa + 
\frac{1}{q_{1}}a)$, and $c_{1} = -\kappa-\frac{1}{q_{1}}a-\frac{1}{q_{2}}b$. In 
particular, $a_{1} > 0$,

\clearpage 
\begin{figure}[ht]
\begin{center}
\DtwittleFigure\ (cont'd):  
GCM graphs in the ``$\widetilde{D}$'' family with more than five nodes are  
not admissible. \\ 
\ \\ 
\vspace*{-0.125in}

{\footnotesize (In each part of this figure, for any $k \geq 0$, the 
subsequent position 
can be obtained 

from the first position 
by a sequence of legal node firings.)}

\ 

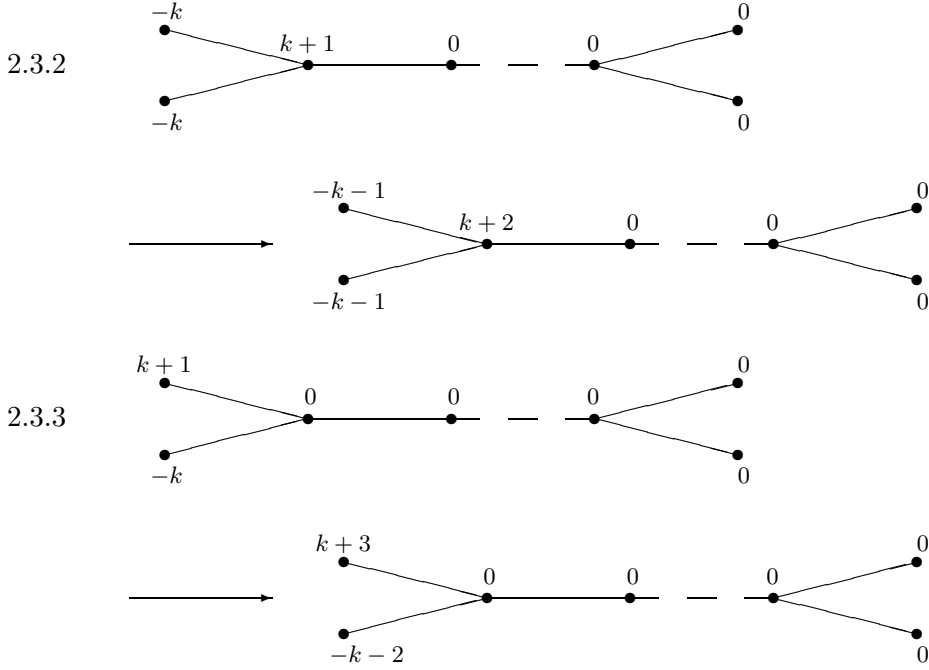
\setlength{\unitlength}{0.75in}
\begin{picture}(5,2.0)
\put(-1.1,1.55){2.3.2}
\put(0,1.25){\begin{picture}(1,0)
 \put(-0.1,0.675){\footnotesize $-k$}
 \put(-0.1,-0.1){\footnotesize $-k$}
 \put(0.8,0.45){\footnotesize $k+1$}
 \put(1.975,0.45){\footnotesize $0$}
 \put(2.95,0.45){\footnotesize $0$}
 \put(4.0,0.675){\footnotesize $0$}
 \put(4.0,-0.1){\footnotesize $0$}
            \put(0,0.6){\circle*{0.075}}
            \put(0,0.1){\circle*{0.075}}
            \put(1,0.35){\circle*{0.075}}
            \put(2,0.35){\circle*{0.075}}
            \put(3,0.35){\circle*{0.075}}
            \put(4,0.1){\circle*{0.075}}
            \put(4,0.6){\circle*{0.075}}
            \put(0,0.1){\line(4,1){1}}
            \put(0,0.6){\line(4,-1){1}}
            \put(1,0.35){\line(1,0){1}}
            \multiput(2,0.35)(0.4,0){3}{\line(1,0){0.2}}
            \put(3,0.35){\line(4,1){1}}
            \put(3,0.35){\line(4,-1){1}}
           \end{picture}}
\put(-0.25,0.35){\vector(1,0){1}}
\put(1.25,0){\begin{picture}(1,0)
 \put(-0.225,0.675){\footnotesize $-k-1$}
 \put(-0.225,-0.1){\footnotesize $-k-1$}
 \put(0.8,0.45){\footnotesize $k+2$}
 \put(1.975,0.45){\footnotesize $0$}
 \put(2.95,0.45){\footnotesize $0$}
 \put(4.0,0.675){\footnotesize $0$}
 \put(4.0,-0.1){\footnotesize $0$}
            \put(0,0.6){\circle*{0.075}}
            \put(0,0.1){\circle*{0.075}}
            \put(1,0.35){\circle*{0.075}}
            \put(2,0.35){\circle*{0.075}}
            \put(3,0.35){\circle*{0.075}}
            \put(4,0.1){\circle*{0.075}}
            \put(4,0.6){\circle*{0.075}}
            \put(0,0.1){\line(4,1){1}}
            \put(0,0.6){\line(4,-1){1}}
            \put(1,0.35){\line(1,0){1}}
            \multiput(2,0.35)(0.4,0){3}{\line(1,0){0.2}}
            \put(3,0.35){\line(4,1){1}}
            \put(3,0.35){\line(4,-1){1}}
           \end{picture}}
\end{picture}

\ 

\begin{picture}(5,2.2)
\put(-1.1,1.55){2.3.3}
\put(0,1.25){\begin{picture}(1,0)
 \put(-0.2,0.675){\footnotesize $k+1$}
 \put(-0.1,-0.1){\footnotesize $-k$}
 \put(0.975,0.45){\footnotesize $0$}
 \put(1.975,0.45){\footnotesize $0$}
 \put(2.95,0.45){\footnotesize $0$}
 \put(4.0,0.675){\footnotesize $0$}
 \put(4.0,-0.1){\footnotesize $0$}
            \put(0,0.6){\circle*{0.075}}
            \put(0,0.1){\circle*{0.075}}
            \put(1,0.35){\circle*{0.075}}
            \put(2,0.35){\circle*{0.075}}
            \put(3,0.35){\circle*{0.075}}
            \put(4,0.1){\circle*{0.075}}
            \put(4,0.6){\circle*{0.075}}
            \put(0,0.1){\line(4,1){1}}
            \put(0,0.6){\line(4,-1){1}}
            \put(1,0.35){\line(1,0){1}}
            \multiput(2,0.35)(0.4,0){3}{\line(1,0){0.2}}
            \put(3,0.35){\line(4,1){1}}
            \put(3,0.35){\line(4,-1){1}}
           \end{picture}}
\put(-0.25,0.35){\vector(1,0){1}}
\put(1.25,0){\begin{picture}(1,0)
 \put(-0.2,0.675){\footnotesize $k+3$}
 \put(-0.1,-0.1){\footnotesize $-k-2$}
 \put(0.975,0.45){\footnotesize $0$}
 \put(1.975,0.45){\footnotesize $0$}
 \put(2.95,0.45){\footnotesize $0$}
 \put(4.0,0.675){\footnotesize $0$}
 \put(4.0,-0.1){\footnotesize $0$}
            \put(0,0.6){\circle*{0.075}}
            \put(0,0.1){\circle*{0.075}}
            \put(1,0.35){\circle*{0.075}}
            \put(2,0.35){\circle*{0.075}}
            \put(3,0.35){\circle*{0.075}}
            \put(4,0.1){\circle*{0.075}}
            \put(4,0.6){\circle*{0.075}}
            \put(0,0.1){\line(4,1){1}}
            \put(0,0.6){\line(4,-1){1}}
            \put(1,0.35){\line(1,0){1}}
            \multiput(2,0.35)(0.4,0){3}{\line(1,0){0.2}}
            \put(3,0.35){\line(4,1){1}}
            \put(3,0.35){\line(4,-1){1}}
           \end{picture}}
\end{picture}
\end{center}
\end{figure}

\noindent
$b_{1} > 0$, and $c_{1} < 0$.  Next we 
check that $\kappa_{1} := (2p_{1}+2p_{2}-\frac{1}{q_{1}})a_{1} + 
(p_{1}+2p_{2}-\frac{1}{q_{2}})b_{1} + c_{1}$ is also positive.  
Now 
\begin{eqnarray*} 
\kappa_{1} & = & 
\left(2p_{1}+2p_{2}-\frac{1}{q_{1}}\right)
\left[q_{1}\left(2p_{1}+2p_{2}-\frac{1}{q_{1}}\right) +  
q_{2}\left(p_{1}+2p_{2}-\frac{1}{q_{2}}\right) - 1\right]\, a\\
 &  & \hspace*{0.3in}+\ \left(p_{1}+2p_{2}-\frac{1}{q_{2}}\right)
\left[q_{1}\left(2p_{1}+2p_{2}-\frac{1}{q_{1}}\right) + 
q_{2}\left(p_{1}+2p_{2}-\frac{1}{q_{2}}\right) - 1\right]\, b\\
 &  & \hspace*{0.3in}+\ \left[q_{1}\left(2p_{1}+2p_{2}-\frac{1}{q_{1}}\right) + 
q_{2}\left(p_{1}+2p_{2}-\frac{1}{q_{2}}\right) - 1\right]\, c\\
 &  & \hspace*{0.3in}+\ \left[\frac{-1}{q_{1}} + 
 \frac{q_{2}}{q_{1}}\left(p_{1}+2p_{2}-\frac{1}{q_{2}}\right)\right]\, a 
 + \left[\frac{-1}{q_{2}} + 
 \frac{q_{1}}{q_{2}}\left(2p_{1}+2p_{2}-\frac{1}{q_{1}}\right)\right]\, b 
\end{eqnarray*} 
Observe that $q_{1}\left(2p_{1}+2p_{2}-\frac{1}{q_{1}}\right) +  
q_{2}\left(p_{1}+2p_{2}-\frac{1}{q_{2}}\right) - 1$ is positive. Since  
$\frac{-1}{q_{1}} + 
 \frac{q_{2}}{q_{1}}\left(p_{1}+2p_{2}-\frac{1}{q_{2}}\right)$, 
$\frac{-1}{q_{2}} + 
 \frac{q_{1}}{q_{2}}\left(2p_{1}+2p_{2}-\frac{1}{q_{1}}\right)$, 
and $\kappa$ are positive as well, then $\kappa_{1} > 0$. 
Thus, $(a_{1},b_{1},c_{1})$ meets condition ({\tt *}), so we 
can repeat the above firing sequence from position $(a_{1},b_{1},c_{1})$ 
to obtain a position $(a_{2},b_{2},c_{2})$ that meets condition 
({\tt *}), etc.  
Since the fundamental 
positions $(a,b,c) = (1,0,0)$ and $(a,b,c) = (0,1,0)$ meet  
condition ({\tt *}), then we see that the indicated 
firing sequence can be repeated indefinitely from these positions.  
For the fundamental 
position $(a,b,c) = (0,0,1)$, begin by firing at the leftmost node to 
obtain the position $(q_{1},q_{2},-1)$.  This latter position meets  
condition ({\tt *}), and so the firing sequence indicated above can 
be repeated indefinitely from this position. 

\underline{Step 4: Every node is fired.} The following is 
proved easily with an induction argument on the number of 
nodes.  

\noindent 
{\bf \EveryNodeFiredLemma}\ \ {\sl Let $\lambda$ be a 
nonzero dominant position 
assigned to the nodes of a connected GCM graph $(\Gamma,M)$.  Then in any 
game sequence for $\lambda$, every node of $\Gamma$ is fired at least 
once.}

\underline{Step 5: Subgraphs.}  
If $I'_{m}$ is a subset of the node set $I_{n}$ of a GCM graph 
$(\Gamma,M)$, 
then let $\Gamma'$ be the subgraph of $\Gamma$ with node set $I'_{m}$ and 
the induced set of edges, and let $M'$ be the corresponding submatrix 
of the amplitude matrix $M$; we call $(\Gamma',M')$ a {\em GCM 
subgraph} of $(\Gamma,M)$.  In light of \EveryNodeFiredLemma, the 
following result amounts to an observation. 

\noindent 
{\bf \SubgraphLemma}\ \ {\sl If a connected GCM graph is 
admissible, then any 
connected GCM subgraph is also admissible.} 

\underline{Step 6: Amplitude products must be 1, 2, or 3.}

\noindent 
{\bf \TwoNodeLemma}\ \ {\sl If $\gamma_{i}$ and $\gamma_{j}$ are adjacent 
nodes in a 
connected admissible GCM graph $(\Gamma,M)$, then the product of the 
amplitudes $M_{ij}M_{ji}$ is 1, 2, or 3.  That 
is, the GCM subgraph of $(\Gamma,M)$ with nodes $i$ and $j$ 
is in this case one of $A_{2}$, $B_{2}$, or $G_{2}$.}

{\em Proof.} By \SubgraphLemma\ we may restrict attention to the admissible  
GCM subgraph $(\Gamma',M')$ 
with node set $\{i,j\}$.  A nonzero dominant position with 
a convergent game sequence might not begin with 
positive populations at 
\underline{both} nodes; nonetheless, by examining the 
proof one sees that Lemma 
3.7 of \cite{ErikssonEur} still applies to show that the product 
$M_{ij}M_{ji}$ of amplitudes in the admissible GCM 
graph $(\Gamma',M')$ is 1, 2, or 3.  

\underline{Conclusion of the ``only if'' part of the first claim of 
\MarsAttacksTheorem.} 
Putting Steps 1 through 6 together, 
we see that the only possible connected 
admissible GCM graphs are the Dynkin diagrams. 

\noindent 
{\em Proof of the remaining claims of 
\MarsAttacksTheorem:} 

Let $(\Gamma,M)$ be a connected Dynkin diagram.  
\StrongConvergenceTheorem\ (Eriksson's Strong Convergence Theorem) 
shows that if a game sequence for some initial position $\lambda$ converges, 
then all game sequences from $\lambda$ converge to the same 
terminal position in the 
same finite number of steps.  
Then in light of 
\ComparisonTheorem\ (Eriksson's Comparison Theorem), it  
suffices to show that for 
any strongly dominant initial position assigned to the nodes of 
$(\Gamma,M)$, 
there is a convergent game sequence.  For the exceptional graphs 
($E_{6}$, $E_{7}$, $E_{8}$, $F_{4}$, and $G_{2}$) this can be checked by 
hand (requiring 36, 63, 120, 24, and 6 firings respectively).  
For the four infinite families of Dynkin diagrams, 
the next result can be proved by induction. 

\noindent 
{\bf \FourFamiliesFirstLemma}\ \ {\sl For any positive integer $n$ 
(respectively, any integer $n \geq 2$, $n \geq 3$, $n \geq 4$) and for 
any strongly dominant position $(a_{1},\ldots,a_{n})$ 
assigned to the nodes of 
$A_{n}$ (respectively, $B_{n}$, $C_{n}$, $D_{n}$), one can obtain the 
position $(a_{1}+\cdots+a_{n},-a_{n},\ldots,-a_{3},-a_{2})$ 
(respectively 
$(a_{1}+2a_{2}+\cdots+2a_{n-1}+a_{n},-a_{2},-a_{3},\ldots,-a_{n})$, 
$(a_{1}+2a_{2}+\cdots+2a_{n-1}+2a_{n},-a_{2},-a_{3},\ldots,-a_{n})$, 
$(a_{1}+2a_{2}+\cdots+2a_{n-2}+a_{n-1}+a_{n},-a_{2},-a_{3},\ldots,
-a_{n-2},-b_{n-1},-b_{n})$ where $b_{n-1} := a_{n-1}$ and $b_{n} := 
a_{n}$ when $n$ is odd 
and where $b_{n-1} := a_{n}$ and $b_{n} := a_{n-1}$ 
when $n$ is even) by a sequence of $\frac{n(n-1)}{2}$ (resp.\ 
$(n-1)^{2}$, $(n-1)^{2}$, $(n-1)(n-2)$) node firings.}

From this we can obtain the following (again by induction): 

\noindent 
{\bf \FourFamiliesSecondLemma}\ \ {\sl For any positive integer $n$ 
(respectively, any integer $n \geq 2$, $n \geq 3$, $n \geq 4$) and for 
any strongly dominant position $(a_{1},\ldots,a_{n})$ 
assigned to the nodes of 
$A_{n}$ (respectively, $B_{n}$, $C_{n}$, $D_{n}$), one can obtain the 
terminal position $(-a_{n},\ldots,-a_{2},-a_{1})$ 
(respectively 
$(-a_{1},-a_{2},\ldots,-a_{n})$, 
$(-a_{1},-a_{2},\ldots,-a_{n})$, 
$(-a_{1},-a_{2},\ldots,
-a_{n-2},-b_{n-1},-b_{n})$ where $b_{n-1} := a_{n-1}$ and $b_{n} := 
a_{n}$ when $n$ is even  
and where $b_{n-1} := a_{n}$ and $b_{n} := a_{n-1}$ 
when $n$ is odd) by a sequence of $\frac{n(n+1)}{2}$ (resp.\ 
$n^{2}$, $n^{2}$, $n(n-1)$) node firings.} 

This completes the proof of \MarsAttacksTheorem.\hfill\QED

\vspace{2ex} 

\noindent
{\Large \bf \FirstProofCorollaryNum.\ \ Classifications of finite-dimensional 
Kac-Moody 
algebras and finite Weyl groups}

\vspace{2ex} 
Since our first proof of \MarsAttacksTheorem\ does not depend on the 
classifications of the finite-dimensional Kac-Moody 
algebras or of 
the finite Weyl groups, we can use \MarsAttacksTheorem\ to obtain 
these results.  
In fact only the weaker result \WildbergerTheorem\ due to Wildberger 
is needed for our proof.  
This is recorded as \ProofCorollary\ below. 
These classifications are obtained 
in \cite{Kac} and \cite{HumCoxeter} respectively by carefully 
studying properties of the generalized Cartan matrix 
(or a closely related 
matrix).   
The definitions we use here basically follow \cite{Kumar} (but see 
also \cite{Kac}).  
The Lie algebra that is constructed next 
does not depend on the specific choices made.  
Given a GCM graph $(\Gamma,M)$ with $n$ nodes, choose a complex vector space 
$\mathfrak{h}$ of dimension $n + \mathrm{corank}(M)$.  Choose $n$ 
linearly independent vectors $\{\beta_{i}^{\vee}\}_{1 \leq i \leq n}$ 
in $\mathfrak{h}$, and find $n$ linearly independent functionals 
$\{\beta_{i}\}_{1 \leq i \leq n}$ in $\mathfrak{h}^{*}$ satisfying 
$\beta_{j}(\beta_{i}^{\vee}) = M_{ij}$.  
The {\sl Kac-Moody 
algebra} $\mathfrak{g} = 
\mathfrak{g}(\Gamma,M)$ is the Lie 
algebra over $\mathbb{C}$ generated by the set $\mathfrak{h} \cup 
\{x_{i},y_{i}\}_{i \in I_{n}}$ with relations 
$[\mathfrak{h},\mathfrak{h}] = 0$; 
$[h,x_{i}] = \beta_{i}(h)x_{i}$ and $[h,y_{i}] = 
-\beta_{i}(h)y_{i}$ for all $h \in \mathfrak{h}$ and $i \in I_{n}$; 
$[x_{i},y_{j}] = \delta_{i,j}\beta_{i}^{\vee}$ for all $i, j \in I_{n}$; 
$(\mathrm{ad} x_{i})^{1-M_{ji}}(x_{j}) = 0$ for $i \not= j$;  
and 
$(\mathrm{ad} y_{i})^{1-M_{ji}}(y_{j}) = 0$ for $i \not= j$, 
where $(\mathrm{ad} z)^{k}(w) = [z,[z,\cdots,[z,w]\cdots]]$.  
It is known (see for example Proposition 1.3.21 of \cite{Kumar}) 
that the associated {\em Weyl group} $W = W(\Gamma,M)$ with identity 
denoted $e$ has the following presentation by generators and 
relations: generators 
$\{s_{i}\}_{i \in I_{n}}$, and relations $s_{i}^{2} = 
e$ for $i \in I_{n}$ and $(s_{i}s_{j})^{m_{ij}} = e$ where the 
integers $m_{ij}$ are determined as follows: 
$m_{ij} = 2$ if $M_{ij}M_{ji} = 0$, 
$m_{ij} = 3$ if $M_{ij}M_{ji} = 1$, 
$m_{ij} = 4$ if $M_{ij}M_{ji} = 2$, 
$m_{ij} = 6$ if $M_{ij}M_{ji} = 3$, and 
$m_{ij} = \infty$ if $M_{ij}M_{ji} \geq 4$. 
This group is the same as the Coxeter group Eriksson associates to 
the GCM graph $(\Gamma,M)$ in \S 3 of \cite{ErikssonDiscrete}. 
If the graph $\Gamma$ has connected components $\Gamma_{1}, 
\ldots, \Gamma_{k}$ with corresponding amplitude matrices 
$M_{1}, \ldots, M_{k}$, then 
$\mathfrak{g}(\Gamma,M) \approx \mathfrak{g}(\Gamma_{1},M_{1}) \oplus \cdots 
\oplus \mathfrak{g}(\Gamma_{k},M_{k})$ and $W(\Gamma,M) \approx 
W(\Gamma_{1},M_{1}) \times \cdots 
\times W(\Gamma_{k},M_{k})$. 
For a firing sequence $(\gamma_{i_{1}}, 
\gamma_{i_{2}}, \ldots, \gamma_{i_{j}})$ from some initial position 
$\lambda$ in a numbers game on $(\Gamma,M)$, 
the corresponding element of $W$ is 
$s_{i_{j}}\cdots{s}_{i_{2}}s_{i_{1}}$. 
The next result follows from Proposition 4.1 of 
\cite{ErikssonDiscrete} and is a key step in our proof of 
\ProofCorollary. 

\noindent 
{\bf \ErikssonWordProposition\ (Eriksson's 
Reduced Word Result)}\ \ {\sl If 
$(\gamma_{i_{1}}, 
\gamma_{i_{2}}, \ldots, \gamma_{i_{j}})$ is a legal sequence of node 
firings in a numbers game played  
from some initial position on a GCM graph $(\Gamma,M)$, 
then $s_{i_{j}}\cdots{s}_{i_{2}}s_{i_{1}}$ is a reduced expression for 
the corresponding element of the Weyl group $W(\Gamma,M)$.}

\noindent 
{\bf \ProofCorollary}\ \ {\sl Given a generalized Cartan matrix,   
the associated Weyl group is finite if and only if the 
associated Kac-Moody  
algebra is finite-dimensional if and only if 
the associated GCM graph is a Dynkin diagram.}

{\em Proof.} 
Suppose a connected GCM graph $(\Gamma,M)$ is not a Dynkin diagram.  
Then by \CorollaryTheoremList, 
we may pick a nonzero dominant position $\lambda$ as an 
initial position and play a nonterminating numbers game. 
By Eriksson's 
Reduced Word Result, 
it follows that the sequence of the first $k$ firings (for any $k > 0$) 
of the corresponding game sequence  
will correspond to a reduced word in the corresponding Weyl group $W$.  
Thus our numbers game 
generates an infinite number of Weyl group elements, and  
hence $|W| = \infty$.  From Proposition 1.4.2 of \cite{Kumar}, 
it follows that the set of roots 
of the Kac-Moody 
algebra $\mathfrak{g}$ 
associated to $(\Gamma,M)$ is infinite.  From the 
root space decomposition of $\mathfrak{g}$ 
(see \S 1.2 of \cite{Kumar}), 
it now follows that $\mathfrak{g}$ is infinite-dimensional.  
Conversely, 
that the Kac-Moody 
algebra $\mathfrak{g}$ associated to a 
connected Dynkin diagram $(\Gamma,M)$ is finite-dimensional follows 
from Serre's Theorem (Theorem 18.3 of \cite{Hum}) in combination with 
the construction in \S 12 of that text of the root system for the 
connected Dynkin diagram; 
finiteness of this root system implies the finiteness of the 
associated Weyl group (again Proposition 1.4.2 of \cite{Kumar}). 
\hfill\QED

It is well known that the Kac-Moody  
algebras associated to the Dynkin 
diagrams of \DynkinDiagramFigure\ are the complex 
finite-dimensional simple Lie algebras (see for example \cite{Hum} \S 
18). 
It is also well known that Lie algebras corresponding to distinct 
Dynkin diagrams of \DynkinDiagramFigure\ are non-isomorphic;   
for the associated Weyl groups, the only 
redundancy is that the 
groups corresponding to the $B_{n}$ and $C_{n}$ graphs for $n \geq 3$  
are the same.

\vspace{2ex} 

\noindent
{\Large \bf \SecondProofNum.\ \ 
Quasi-standard geometric representations of Coxeter groups and a 
generalization to E-games}

\vspace{2ex} 
In this section we study a certain geometric representation of a 
Coxeter group associated to a given E-GCM.   
This geometric 
representation is developed in 
\cite{ErikssonDiscrete} \S 3, 4 and in \cite{BB} \S 4.1, 4.2; it is  
similar to but not in general the same as the standard geometric 
representation (as developed in \cite{HumCoxeter} Ch.\ 5, for 
example).  
Results we derive about geometric aspects of 
this representation extend the standard theory and will 
lead to a proof of the classification result \EGCMTheorem.  
The classification of finite Coxeter groups 
(\cite{HumCoxeter} \S 6.4, 2.7, 2.4)
is used in our proof of 
\EGCMTheorem. 

Now and for the remainder of this section 
fix an E-GCM graph $(\Gamma,M)$.  
Define the associated Coxeter group 
$W = W(\Gamma,M)$ to be the Coxeter group with identity 
denoted $e$, 
generators $\{s_{i}\}_{i \in I_{n}}$, and relations $s_{i}^{2} = 
e$ for $i \in I_{n}$ and 
$(s_{i}s_{j})^{m_{ij}} = e$ for all $i \not= j$, where the 
integers $m_{ij}$ are determined as follows: 

\vspace*{-0.25in}
\[m_{ij} = \left\{\begin{array}{cl}
k & \mbox{\hspace*{0.25in} if 
$M_{ij}M_{ji} = 4\cos^{2}(\pi/k_{ij})$ for some integer $k_{ij} \geq 2$}\\ 
\infty & \mbox{\hspace*{0.25in} if 
$M_{ij}M_{ji} \geq 4$} 
\end{array}\right.\]  
 
Observe that any Coxeter group on a finite set of generators is isomorphic 
to $W(\Gamma,M)$ for some E-GCM graph $(\Gamma,M)$; the Coxeter 
group is {\em irreducible} if $\Gamma$ is connected.  By the 
classification of finite Coxeter groups, an irreducible Coxeter group 
is finite if and only if it is the Coxeter group associated to an 
E-Coxeter graph from \ECoxeterGraphFigure. 
Before proceeding toward our proof of \EGCMTheorem, we note the 
following.  
In Propositions 4.1 and 4.2 of \cite{Deodhar}, Deodhar gives a number of 
statements equivalent to the assertion that a given irreducible 
Coxeter group is finite.  
As an immediate consequence of \EGCMTheorem\ and the classification of 
finite Coxeter groups, we add to that list 
the following equivalence.  

\noindent 
{\bf \DeodharEquivalenceCorollary}\ \ 
{\sl An irreducible Coxeter group $W$ is finite if and only if 
there is an admissible E-GCM graph whose associated Coxeter group is 
$W$ if and only if any E-GCM graph is admissible 
when its associated Coxeter group is $W$.}\hfill\QED 

We say two nodes $\gamma_{i}$ and $\gamma_{j}$ 
in $(\Gamma,M)$ are {\em odd-adjacent} if $m_{ij}$ is odd, {\em 
even-adjacent} if $m_{ij} \geq 4$ is even, and $\infty$-adjacent if 
$m_{ij} = \infty$.  
When 
$m_{ij}$ is odd and 
$M_{ij} \not= M_{ji}$, we say that the adjacent nodes $\gamma_{i}$ 
and $\gamma_{j}$ form an {\em odd asymmetry}. 
Note that $m_{ij} = 2$ (respectively 3, 4, 6) when $M_{ij}M_{ji} = 
0$ (resp.\ 1, 2, 3). 
We let ${\myl}$ denote the length function for $W$. 
For $J \subseteq I_{n}$, let 
$W_{J}$ be the subgroup generated by 
$\{s_{i}\}_{i \in J}$, a {\em parabolic} subgroup, and $W^{J} := \{w 
\in W\, |\, \myl(ws_{j}) > \myl(w) \mbox{ for all } j \in J\}$  
is the set of {\em minimal coset representatives}. 
If $J = \{i,j\}$, then $W_{J}$ is a dihedral group of order 
$2m_{ij}$.  

Let $V$ be a real $n$-dimensional vector space freely generated by 
$(\alpha_{i})_{i \in I_{n}}$ (elements of this ordered basis are {\em 
simple roots}).  
Equip $V$ with a possibly asymmetric bilinear form $B: V \times V 
\rightarrow \mathbb{R}$ defined on the 
basis $(\alpha_{i})_{i \in I_{n}}$ by 
$B(\alpha_{i},\alpha_{j}) := \frac{1}{2}M_{ij}$. 
For each $i \in I_{n}$ define an operator 
$S_{i}: V \rightarrow V$ by the rule $S_{i}(v) := v - 
2B(\alpha_{i},v)\alpha_{i}$ for each $v \in V$.  One can check 
that $S_{i}^{2} = I$ (the identity 
transformation), so $S_{i} \in GL(V)$.  With $V_{i,j} := 
\mathrm{span}_{\mathbb{R}}\{\alpha_{i},\alpha_{j}\}$, 
observe that $S_{k}(V_{i,j}) \subseteq V_{i,j}$ for $k = i,j$.  Let 
$\mathfrak{B}$ be the ordered basis $(\alpha_{i},\alpha_{j})$ for 
$V_{i,j}$, and for any linear mapping $T: V_{i,j} \rightarrow 
V_{i,j}$ let $[T]_{\mathfrak{B}}$ be the matrix for $T$ relative to 
$\mathfrak{B}$.  
Then 
\[ [S_{i}|_{V_{i,j}}]_{\mathfrak{B}} = 
\left(\begin{array}{cc} -1 & -M_{ij}\\ 0 & 1\end{array}\right), 
[S_{j}|_{V_{i,j}}]_{\mathfrak{B}} = 
\left(\begin{array}{cc} 1 & 0\\ -M_{ji} & -1\end{array}\right), 
[S_{i}S_{j}|_{V_{i,j}}]_{\mathfrak{B}} = 
\left(\begin{array}{cc} M_{ij}M_{ji}-1 & M_{ij}\\ -M_{ji} 
& -1\end{array}\right)\] 
Analysis of the eigenvalues for $X_{i,j} := 
[S_{i}S_{j}|_{V_{i,j}}]_{\mathfrak{B}}$ 
as in the proof of Proposition 1.3.21 of \cite{Kumar} shows that 
$X_{i,j}$ has infinite order when $M_{ij}M_{ji} 
\geq 4$, and hence $S_{i}S_{j}$ 
has infinite order as an element of $GL(V)$.  When $0 < M_{ij}M_{ji} < 4$, 
write $M_{ij}M_{ji} = 4\cos^{2}\theta$ for $0 < \theta < \pi/2$.  
In this case check that $X_{i,j}$ has two distinct complex eigenvalues 
($e^{2i\theta}$ and $e^{-2i\theta}$).  It follows that $X_{i,j}$ has 
finite order $m_{ij}$ if {\em and only if} $\theta = \pi/m_{ij}$, which 
coincides with 
Eriksson's constraints on the amplitude products of the 
E-GCM.  
When $M_{ij}M_{ji} = 0$, then $X_{i,j}  = 
\left(\begin{array}{cc} -1 & 0\\ 0 & -1\end{array}\right)$, 
which clearly has order $m_{ij} = 2$. 
With $0 \leq M_{ij}M_{ji} < 4$, one can easily see 
that $V = V_{i,j} \oplus V_{i,j}'$, where $V_{i,j}' := 
\{v \in V\, |\, B(\alpha_{i},v) = 0 = 
B(\alpha_{j},v)\}$.  Since $S_{i}S_{j}$ acts as the identity on 
$V_{i,j}'$, it follows that 
$S_{i}S_{j}$ has order $m_{ij}$ as an element of $GL(V)$.  
Then  there is a unique homomorphism 
$\sigma_{M}: W \rightarrow GL(V)$ for which $\sigma_{M}(s_{i}) = 
S_{i}$; we call $\sigma_{M}$ a {\em quasi-standard geometric representation} 
of $W$. 
We now have $W$ acting on $V$, and for all $w \in W$ and $v \in V$ we 
write $w.v$ for $\sigma_{M}(w)(v)$. 
Define $\Phi_{M} := \{\alpha \in V\, |\, \alpha = w.\alpha_{i} \mbox{ for 
some } i \in I_{n} \mbox{ and } w \in W\}$.  Elements of $\Phi_{M}$ are 
{\em roots}; if $\alpha = \sum c_{i}\alpha_{i}$ is a root with all 
$c_{i}$ nonnegative (respectively nonpositive), then say $\alpha$ is a {\em 
positive} (resp.\ {\em negative}) root, and write $\alpha >_{_{M}} 0$ 
(resp.\ $\alpha <_{_{M}} 0$).  Let $\Phi_{M}^{+}$ and $\Phi_{M}^{-}$ 
denote the collections of positive and negative roots respectively; 
it is a consequence of \HumphreysTheorem\ below that $\Phi_{M}$ is 
partitioned by $\Phi_{M}^{+}$ and $\Phi_{M}^{-}$.  
The possible asymmetry of the bilinear form is a crucial difference 
between this quasi-standard 
geometric realization and the standard geometric 
realization; 
for example, 
$\sigma_{M}(W)$ preserves the form $B$ if and only if $M$ is 
symmetric, in which case the representation studied here is the same 
as the standard geometric representation. 
Under this quasi-standard action of $W$ on $V$, 
sometimes $K\alpha_{x}$ is 
a root for $K \not= \pm{1}$, as \UnitalExample\ 
shows (see also Exercise 4.9 of 
\cite{BB}).\myfootnote{Just prior to the statement of Proposition 4.4 in  
\cite{ErikssonDiscrete}, it is mistakenly asserted that the only 
multiples in $\Phi_{M}$ of a simple root $\alpha_{x}$ are 
$\pm\alpha_{x}$.  This only affects Proposition 4.4 of that 
paper.}   
Our proofs of the main results of this section require us to 
understand how this $W$-action generates scalar multiples in $\Phi$.  
First we 
analyze how 
$s_{i}$ and $s_{j}$ act in tandem on $V_{i,j}$. Our next result 
strengthens Lemma 4.2.4 of \cite{BB}. 

\noindent 
{\bf \TwoGeneratorAnalysis}\ \ 
{\sl Fix $i \not= j$ in $I_{n}$, and let $k$ be a positive integer.   
If $m_{ij} = \infty$, then 
$(s_{i}s_{j})^{k}.\alpha_{i} = 
a\alpha_{i} + b\alpha_{j}$ and 
$s_{j}(s_{i}s_{j})^{k}.\alpha_{i} = c\alpha_{i} + d\alpha_{j}$,  
for positive coefficients $a$, $b$, $c$, and $d$.  
Now suppose $m_{ij} < \infty$. 
If $2k<m_{ij}$, then $(s_{i}s_{j})^{k}.\alpha_{i} = a\alpha_{i} + 
b\alpha_{j}$ with $a \geq 0$ and $b > 0$.  Moreover, 
this is a multiple of $\alpha_{j}$ if and only if $m_{ij}$ is odd 
and $k = (m_{ij}-1)/2$, in which case 
$(s_{i}s_{j})^{k}.\alpha_{i} = 
\frac{-M_{ji}}{2\cos(\pi/m_{ij})}\alpha_{j}$.  
Similarly, if $2k<m_{ij}-1$, then 
$s_{j}(s_{i}s_{j})^{k}.\alpha_{i} = a\alpha_{i} + b\alpha_{j}$ with $a > 
0$ and $b \geq 0$.  Moreover, 
this is a multiple of $\alpha_{i}$ if and only if $m_{ij}$ is even 
and $k = (m_{ij}-2)/2$, in which case 
$s_{j}(s_{i}s_{j})^{k}.\alpha_{i} = \alpha_{i}$.}

{\em Proof.} The details are somewhat tedious but routine.  For 
convenience set $p := -M_{ij}$ and 
$q := -M_{ji}$. 
Start with $m_{ij} < \infty$, and set $\theta := \pi/m_{ij}$. 
We will work with the ordered basis $\mathfrak{B}$ for $V_{i,j}$ and 
operators $X_{i} := [S_{i}|_{V_{i,j}}]_{\mathfrak{B}}$ and $X_{j} := 
[S_{j}|_{V_{i,j}}]_{\mathfrak{B}}$ as above.  
Then to understand $(s_{i}s_{j})^{k}.\alpha_{i}$ and 
$s_{j}(s_{i}s_{j})^{k}.\alpha_{i}$ we compute 
$X_{i,j}^{k}$ and $X_{j}X_{i,j}^{k}$.  In the case that $pq = 
4\cos^{2}(\theta)$, then $X_{i,j}$ can be written as 
$X_{i,j} = PDP^{-1}$ for a 
nonsingular matrix $P$ and diagonal matrix $D$ in the following way: 
\[\frac{1}{q(e^{2i\theta} - e^{-2i\theta})}
\left(\begin{array}{cc}e^{2i\theta}+1 & e^{-2i\theta}+1\\ 
q & q\end{array}\right)
\left(\begin{array}{cc} e^{2i\theta} & 0\\ 0 & 
e^{-2i\theta}\end{array}\right) 
\left(\begin{array}{cc} q & -e^{-2i\theta}-1\\ -q 
& e^{2i\theta}+1\end{array}\right).\]
Then for any positive integer $k$ we have  
\[X_{i,j}^{k} = PD^{k}P^{-1} = \frac{1}{\sin(2\theta)}
\left(\begin{array}{cc}
\sin(2(k+1)\theta) + \sin(2k\theta) & 
-p\sin(2k\theta)\\
q\sin(2k\theta) & -\sin(2k\theta) - \sin(2(k-1)\theta)
\end{array}\right)\]
and 
\[X_{j}X_{i,j}^{k} = \frac{1}{\sin(2\theta)}
\left(\begin{array}{cc}
\sin(2(k+1)\theta) + \sin(2k\theta) & 
-p\sin(2k\theta)\\
q\sin(2(k+1)\theta) & (1-pq)\sin(2k\theta) + \sin(2(k-1)\theta)
\end{array}\right)
\] 
Use the first column of $X_{i,j}^{k}$ and 
$X_{j}X_{i,j}^{k}$ 
to see that 
$(s_{i}s_{j})^{k}.\alpha_{i} = \frac{1}{\sin(2\theta)}
[\sin(2(k+1)\theta) + 
\sin(2k\theta)]\alpha_{i} +  
\frac{q}{\sin(2\theta)}\sin(2k\theta)\alpha_{j}$ and 
$s_{j}(s_{i}s_{j})^{k}.\alpha_{i} =  
\frac{1}{\sin(2\theta)}[\sin(2(k+1)\theta) + 
\sin(2k\theta)]\alpha_{i} +  
\frac{q}{\sin(2\theta)}\sin(2(k+1)\theta)\alpha_{j}$.  As long as 
$2(k+1) < m_{ij}$, then all the coefficients of these linear 
combinations will be positive.  So now suppose $2(k+1) \geq m_{ij}$.  
First we consider $(s_{i}s_{j})^{k} = a\alpha_{i} + b\alpha_{j}$ 
for some positive $k$ with $2k < 
m_{ij}$.  There are two possibilities now: $2(k+1) = m_{ij}$ or 
$2(k+1) = m_{ij}+1$.  In the former case 
both $a$ and $b$ are positive.  In the latter case we have $m_{ij}$ 
odd, $a = \frac{1}{\sin(2\theta)}
[\sin(2(k+1)\theta) + 
\sin(2k\theta)] = 0$, and $b = \frac{q\sin\theta}{\sin(2\theta)} = 
\frac{q}{2\cos\theta}$.  Second we consider 
$s_{j}(s_{i}s_{j})^{k} = a\alpha_{i} + b\alpha_{j}$ 
for some positive $k$ with $2k < m_{ij}-1$. Now the fact that $2(k+1) 
\geq m_{ij}$ implies we have $2(k+1) = m_{ij}$.  In particular, 
$m_{ij}$ is even.  With $k = (m_{ij}-2)/2$ now, one can check that 
$b = 0$ and $a = 1$. 

For $m_{ij} = \infty$, first take $pq = 4$.  We can write $X_{i,j} = 
PYP^{-1}$ for nonsingular $P$ and upper triangular $Y$ as follows: 
\[X_{i,j} = \frac{1}{p}
\left(\begin{array}{cc}p & p\\ 
2 & 1\end{array}\right)
\left(\begin{array}{cc} 1 & 1\\ 0 & 
1\end{array}\right) 
\left(\begin{array}{cc} -1 & p\\ 2 
& -p\end{array}\right).\] 
Then for any positive integer $k$ we obtain 
$X_{i,j}^{k} = 
\left(\begin{array}{cc}2k+1 & -kp\\ 
kq & -2k+1\end{array}\right)$.  It follows that 
$(s_{i}s_{j})^{k}.\alpha_{i} = (2k+1)\alpha_{i} +  
kq\alpha_{j}$, with both coefficients of the linear combination 
positive. From the first column of the matrix $X_{j}X_{i,j}^{k}$ we 
see that  
$s_{j}(s_{i}s_{j})^{k}.\alpha_{i} = (2k+1)\alpha_{i} +  
(2k+1)q\alpha_{j}$, with both coefficients of the linear combination 
positive.    
Next take $pq > 4$.  In this case we get distinct eigenvalues 
$\lambda = \frac{1}{2}(pq-2+\sqrt{pq(pq-4)}) > 1$ and  
$\mu = \frac{1}{2}(pq-2-\sqrt{pq(pq-4)}) < 1$ for $X_{i,j}$ (here we 
have $\lambda\mu = 1$).  As before, write
$X_{i,j} = PDP^{-1}$ for the diagonal 
matrix $D = \left(\begin{array}{cc}\lambda & 0\\ 
0 & \mu\end{array}\right)$ 
and a nonsingular matrix $P$ to obtain 
\[X_{i,j}^{k} = \frac{1}{p(\lambda-\mu)}
\left(\begin{array}{cc}p & p\\ 
\lambda' & \mu'\end{array}\right)
\left(\begin{array}{cc} \lambda^{k} & 0\\ 0 & 
\mu^{k}\end{array}\right) 
\left(\begin{array}{cc} \mu' & -p\\ -\lambda' 
& p\end{array}\right),\]
for any positive integer $k$, with $\lambda' := \mu+1$  
and $\mu' := \lambda+1$. This (eventually) simplifies to 
\[X_{i,j}^{k} =  
\frac{1}{\lambda-\mu}
\left(\begin{array}{cc}\mu'\lambda^{k} - \lambda'\mu^{k} & 
-p(\lambda^{k}-\mu^{k})\\ 
q(\lambda^{k}-\mu^{k}) & \mu'\mu^{k} - 
\lambda'\lambda^{k}\end{array}\right).\] From this we also get 
\[X_{j}X_{i,j}^{k} =  
\frac{1}{\lambda-\mu}
\left(\begin{array}{cc}\mu'\lambda^{k} - \lambda'\mu^{k} & 
-p(\lambda^{k}-\mu^{k})\\ 
q(\lambda^{k+1}-\mu^{k+1}) & \lambda'\mu^{k} - 
\mu'\lambda^{k}\end{array}\right).\] 
The factor $\frac{1}{\lambda-\mu}$ is positive, 
and for both matrices $X_{i,j}^{k}$ and 
$X_{j}X_{i,j}^{k}$, 
the first column entries are positive.  
So,  
$(s_{i}s_{j})^{k}.\alpha_{i} = a\alpha_{i} + b\alpha_{j}$ with both 
$a$ and $b$ positive, and 
$s_{j}(s_{i}s_{j})^{k}.\alpha_{i} = c\alpha_{i} + d\alpha_{j}$ with $c$ 
and $d$ both positive.\hfill\QED

In view of this result, for odd $m_{ij}$, 
let $v_{ji}$ be the element 
$(s_{i}s_{j})^{(m_{ij}-1)/2}$, and 
set $K_{ji} := \frac{-M_{ji}}{2\cos(\pi/m_{ij})}$, which is positive. 
Then 
$v_{ji}.\alpha_{i} = K_{ji}\alpha_{j}$.  
Observe that $K_{ij}K_{ji} = 1$ and moreover that 
$v_{ij} = v_{ji}^{-1}$.  
A {\em path with odd adjacencies} (or {\em OA-path}, 
for short) 
in $(\Gamma,M)$ is a sequence  
$\mathcal{P} := 
[\gamma_{i_{0}},\gamma_{i_{1}},\ldots,\gamma_{i_{p}}]$ of pairwise 
odd-adjacent nodes of $\Gamma$; this OA-path has length $p$, and we allow 
OA-paths to have length zero.  
We say $\gamma_{i_{0}}$ and 
$\gamma_{i_{p}}$ are the {\em start} and {\em end} nodes of the OA-path, 
respectively.  
If OA-path $\mathcal{Q} = 
[\gamma_{j_{0}},\gamma_{j_{1}},\ldots,\gamma_{j_{q}}]$ 
has the same start node as the end node of  
$\mathcal{P}$, then their {\em concatenation} 
$\mathcal{P}\natural\mathcal{Q}$ is the OA-path 
$[\gamma_{i_{0}},\gamma_{j_{1}},\ldots,
\gamma_{i_{p}=j_{0}},\ldots,\gamma_{j_{q}}]$. 
Let $w_{_{\mathcal{P}}} \in W$ be the Coxeter group element 
$v_{i_{p}i_{p-1}}\cdots{v}_{i_{2}i_{1}}v_{i_{1}i_{0}}$, and 
let $\Pi_{_{\mathcal{P}}} := 
K_{i_{p}i_{p-1}}\cdots{K}_{i_{2}i_{1}}K_{i_{1}i_{0}}$, where 
$w_{_{\mathcal{P}}} = e$ with $\Pi_{_{\mathcal{P}}} = 1$ 
when $\mathcal{P}$ has length zero. 
Note that $w_{_{\mathcal{P}}}.\alpha_{i_{0}} = 
\Pi_{_{\mathcal{P}}}\alpha_{i_{p}}$ and that 
$w_{_{\mathcal{P}\natural\mathcal{Q}}} = 
w_{_{\mathcal{Q}}}w_{_{\mathcal{P}}}$. 

\noindent 
{\bf \HumphreysTheorem}\ \ 
{\sl Let $w \in W$ and $i \in I_{n}$.  If $\myl(ws_{i}) > 
\myl(w)$, then $w.\alpha_{i} >_{_{M}} 0$, and in this case 
$w.\alpha_{i} = 
K\alpha_{x}$ ($x \in I_{n}, K > 0$) 
if and only if $w.\alpha_{i} = 
w_{_{\mathcal{P}}}.\alpha_{i}$ for some OA-path $\mathcal{P} = 
[\gamma_{i_{0}=i},
\gamma_{i_{1}},\ldots,\gamma_{i_{p-1}},\gamma_{i_{p}=x}]$, so 
$K = \Pi_{_{\mathcal{P}}}$.  
Similarly, if $\myl(ws_{i}) < 
\myl(w)$, then $w.\alpha_{i} <_{_{M}} 0$, and in this case $w.\alpha_{i} = 
K\alpha_{x}$ ($x \in I_{n}, K<0$) 
if and only if $w.\alpha_{i} = (w_{_{\mathcal{P}}}s_{i}).\alpha_{i}$ 
for some OA-path $\mathcal{P} = 
[\gamma_{i_{0}=i},
\gamma_{i_{1}},\ldots,\gamma_{i_{p-1}},\gamma_{i_{p}=x}]$, so 
$K = -\Pi_{_{\mathcal{P}}}$.} 

The assertions that $\myl(ws_{i}) > 
\myl(w) \Rightarrow w.\alpha_{i} >_{_{M}} 0$ and 
$\myl(ws_{i}) < 
\myl(w) \Rightarrow w.\alpha_{i} <_{_{M}} 0$ are Proposition 4.2.5 of 
\cite{BB}.  Our similar proof (below) of these assertions follows 
the proof of  
Theorem 5.4 of \cite{HumCoxeter}, whose set up we require for our 
analysis of roots which are scalar multiples of simple roots.    

{\em Proof of \HumphreysTheorem.} 
The first part of our argument follows the proof of 
Theorem 5.4 from 
\cite{HumCoxeter}; however, the argument here is easier since we may 
use the fact that for $J \subseteq I_{n}$, 
the length function $\myl_{_{J}}$ on $W_{J}$ agrees with $\myl$. 
Note that the second assertion of the theorem follows from the first.  For the 
first assertion  
of the theorem, induct on $\myl(w)$.  When $\myl(w) = 0$, there 
is nothing to prove.  Now suppose $\myl(w) > 0$.  Take any $j \in I_{n}$ 
for which $\myl(ws_{j}) = \myl(w) - 1$; since $\myl(ws_{i}) > \myl(w)$, then $i 
\not= j$.  Let $J := \{i,j\}$, and let $v$ be the 
unique element in $W^{J}$ and $v_{_{J}}$ the unique element in $W_{J}$ 
for which $w = vv_{_{J}}$.  Then $\myl(w) = \myl(v) + \myl(v_{_{J}})$.  Observe that 
$\myl(v_{_{J}}) > 0$ since $\myl(v) < \myl(w)$.  From 
Humphreys' proof, we can see that $\myl(vs_{i}) > \myl(v)$, $\myl(vs_{j}) > 
\myl(v)$, and $\myl(v_{_{J}}s_{i}) > \myl(v_{_{J}})$.  
Apply the induction hypothesis to conclude that $v.\alpha_{i} >_{_{M}} 0$  
and $v.\alpha_{j} >_{_{M}} 0$.  
It is possible that $\myl(w) = \myl(v_{_{J}})$, so the induction hypothesis 
might not apply to $v_{_{J}}$.  
But since 
$\myl(v_{_{J}}s_{i}) > \myl(v_{_{J}})$, it follows that any reduced expression 
for $v_{_{J}}$ (necessarily an alternating product of 
$s_{i}\, '${\footnotesize s} and 
$s_{j}\, '${\footnotesize s}) 
must end in $s_{j}$.  Then we may apply \TwoGeneratorAnalysis\ to 
conclude that $v_{_{J}}.\alpha_{i} >_{_{M}} 0$.  
Together these facts imply that $w.\alpha_{i} >_{_{M}} 0$.  

Now we address the issue of scalar multiples.   Suppose $w.\alpha_{i} = 
K\alpha_{x}$ for some $x \in I_{n}$ and real number $K > 0$.  
Write $v_{_{J}}.\alpha_{i} = 
a\alpha_{i} + b\alpha_{j}$, and suppose $a > 0$ and $b > 0$.  Note 
that $v.\alpha_{i}$ and $v.\alpha_{j}$ cannot both be multiples of 
the same $\alpha_{x}$ (otherwise $v^{-1}.\alpha_{x}$ is a multiple of 
both $\alpha_{i}$ and $\alpha_{j}$).  But now a 
simple calculation shows that in this case $w.\alpha_{i}$ will not be 
a multiple of $\alpha_{x'}$ for any $x' \in I_{n}$.  We conclude that 
$v_{_{J}}.\alpha_{i}$ must be a 
scalar multiple of $\alpha_{i}$ or $\alpha_{j}$. 
If $v_{_{J}} = s_{j}$, then $w.\alpha_{i} = v.\alpha_{i}$, and the 
induction hypothesis now applies to $v$ to obtain the desired 
result.  
So now suppose $\myl(v_{_{J}})>1$.  Then $v_{_{J}} = (s_{i}s_{j})^{k}$ or 
$s_{j}(s_{i}s_{j})^{k}$ for some positive integer $k$. 
From \TwoGeneratorAnalysis, it 
follows that if $m_{ij} = \infty$, then $v_{_{J}}.\alpha_{i} = 
a\alpha_{i} + b\alpha_{j}$ with $a$ and $b$ both positive.  
Therefore $m_{ij}$ is finite.  In this case the longest element in 
$W_{J}$ has length $m_{ij}$ and can be written in two ways, one ending 
in $s_{i}$.  Therefore $\myl(v_{_{J}}) < m_{ij}$.  
We are again in the situation of 
\TwoGeneratorAnalysis.  If $v_{_{J}}.\alpha_{i} = \alpha_{i}$, then 
$w.\alpha_{i} = v.\alpha_{i}$, and we can apply the  
induction hypothesis to $v$.  If $v_{_{J}}.\alpha_{i} \not= \alpha_{i}$, then 
we see that $m_{ij}$ is odd and $v_{_{J}}.\alpha_{i} = 
v_{ji}.\alpha_{i} = K_{ji}\alpha_{j}$.  Again apply the 
induction hypothesis to $v$, where now $v.\alpha_{j} = 
\frac{K}{K_{ji}}\alpha_{x} = w_{_{\mathcal{Q}}}.\alpha_{j}$ for some 
OA-path $\mathcal{Q} = 
[\gamma_{j_{0}=j},\gamma_{j_{1}},\ldots,\gamma_{j_{q}=x}]$. 
Let $\mathcal{P} := [\gamma_{i},\gamma_{j}]\natural\mathcal{Q}$. 
Then $w.\alpha_{i} = K\alpha_{x} = 
w_{_{\mathcal{P}}}.\alpha_{i}$.\hfill\QED 

As with Corollary 5.4 of \cite{HumCoxeter} it follows 
that the representation $\sigma_{M}$ is faithful.  
It also follows that $\Phi_{M}$ is 
partitioned by the sets of positive and negative roots.  
For any $w \in W$, set $N_{M}(w) := 
\{\alpha \in \Phi_{M}^{+}\, |\, w.\alpha \in \Phi_{M}^{-}\}$. 

\noindent
{\bf \PositiveToNegativeRootsLemma}\ \ {\sl For any $i \in I_{n}$, 
$s_{i}(\Phi_{M}^{+}\setminus\{K\alpha_{i}|K\in\mathbb{R}\}) = 
\Phi_{M}^{+}\setminus\{K\alpha_{i}|K\in\mathbb{R}\}$.  Now let $w \in W$.  
If $w.\alpha_{i} >_{_{M}} 0$, then $N_{M}(ws_{i}) = s_{i}(N_{M}(w)) 
\disjointunion 
(\{K\alpha_{i}|K\in\mathbb{R}\} \cap \Phi_{M}^{+})$, a disjoint union. 
If $w.\alpha_{i} <_{_{M}} 0$, then $N_{M}(ws_{i}) = 
s_{i}(N_{M}(w)\setminus\{K\alpha_{i}|K\in\mathbb{R}\})$.} 

{\em Proof.}   The proof of Proposition 5.6.(a) from \cite{HumCoxeter} 
is easily adjusted to prove the first claim.  Proofs for 
the remaining claims involve routine set inclusion 
arguments.\hfill\QED

For $J \subseteq I_{n}$, 
let $(\Phi_{M})^{J} := \{\alpha \in \Phi_{M}^{+}\, |\, 
\alpha \not\in \mathrm{span}_{\mathbb{R}}\{\alpha_{j}\}_{j \in J}\}$. 
Our next result 
observes that an assertion from the proof of 
Proposition 4.2 of \cite{Deodhar} also holds here:  

\noindent 
{\bf \DeodharProp}\ \ {\sl If $(\Gamma,M)$ is connected and $\Phi_{M}$ is 
infinite, then $\Phi_{M}^{J}$ is infinite when $J$ is a 
proper subset of $I_{n}$.}  

{\em Proof:} In the ``({\em ix}) $\Rightarrow$ ({\em ii})'' part of 
the proof of Proposition 4.2 in \cite{Deodhar}, 
begin reading at the assumption ``$|\Phi^{J}| < \infty$,'' 
replacing $\Phi^{J}$ with $\Phi_{M}^{J}$.\hfill\QED

Before we analyze the sets 
$\{K\alpha_{x}|K\in\mathbb{R}\}\cap\Phi_{M}^{+}$, we need some 
further notation.  
An OA-path 
$\mathcal{P} = [\gamma_{i_{0}},\ldots,\gamma_{i_{p}}]$ 
is an {\em OA-cycle} if $\gamma_{i_{p}} = \gamma_{i_{0}}$; 
it is a {\em unital} OA-cycle if 
$\Pi_{\mathcal{P}} = 1$.  For OA-paths $\mathcal{P}$ and $\mathcal{Q}$, 
write $\mathcal{P} \sim 
\mathcal{Q}$ and say $\mathcal{P}$ and $\mathcal{Q}$ are {\em 
equivalent}  
if these OA-paths have 
the same start and end nodes and 
$\Pi_{_{\mathcal{P}}} = \Pi_{_{\mathcal{Q}}}$; this is an 
equivalence relation on the set of all OA-paths. 
An OA-path $\mathcal{P}$ 
is {\em simple} if it has no repeated nodes with the possible 
exception that the start and end nodes may coincide.  We say the E-GCM graph 
$(\Gamma,M)$ is {\em unital OA-cyclic} if and only if 
$\Pi_{_{\mathcal{C}}} = 1$ for all OA-cycles $\mathcal{C}$.  
Note that $(\Gamma,M)$ is unital OA-cyclic if and only if 
$\mathcal{P} \sim \mathcal{Q}$ whenever 
$\mathcal{P}$ and $\mathcal{Q}$ are OA-paths 
with the same start and end nodes.  
If $\Gamma$ is a tree, 
then $(\Gamma,M)$ is unital OA-cyclic (vacuously so).  
From the definitions it follows that  
$(\Gamma,M)$ is unital OA-cyclic if it has no odd asymmetries. 
If $M$ is 
symmetrizable, then by Exercise 2.1 of \cite{Kac}, $(\Gamma,M)$ is 
unital OA-cyclic. 
However, a unital OA-cyclic E-GCM graph need not have a symmetrizable 
amplitude matrix $M$, as \UnitalExample\ shows.  
To check if an E-GCM graph is unital OA-cyclic, it is enough to check 
that each simple OA-cycle is unital. 
An E-GCM graph is {\em OA-connected} if 
any two nodes can be joined by an OA-path.  An {\em OA-connected 
component} of an E-GCM graph $(\Gamma,M)$ is an E-GCM subgraph 
$(\Gamma',M')$ whose nodes 
form a maximal collection of nodes in $(\Gamma,M)$ which 
can be pairwise joined by OA-paths.  

\noindent 
{\bf \OneToOneLemma}\ \ {\sl Let $\gamma_{i}$ and $\gamma_{j}$ be nodes in the 
same OA-connected component of the E-GCM graph $(\Gamma,M)$.  
Then there is a a one-to-one 
correspondence between the sets 
$\{K\alpha_{i}|K\in\mathbb{R}\}\cap\Phi_{M}^{+}$ and 
$\{K'\alpha_{j}|K'\in\mathbb{R}\}\cap\Phi_{M}^{+}$.} 

{\em Proof.} Let 
$\mathfrak{S}_{i} := \{K\alpha_{i}|K\in\mathbb{R}\}\cap\Phi_{M}^{+}$ and 
$\mathfrak{S}_{j} := 
\{K'\alpha_{j}|K'\in\mathbb{R}\}\cap\Phi_{M}^{+}$. 
For $\phi: \mathfrak{S}_{i} \rightarrow 
\mathfrak{S}_{j}$ define $K\alpha_{i} \stackrel{\phi}{\longmapsto} 
K_{ji}K\alpha_{j}$.  This is well-defined since $K\alpha_{i} \in 
\mathfrak{S}_{i}$ means $K = \Pi_{_{\mathcal{P}}}$ for some OA-path 
with some start node $\gamma_{x}$ and end node $\gamma_{i}$ (cf.\ 
\HumphreysTheorem).  Then 
$K_{ji}K = \Pi_{_{\mathcal{P}\natural[\gamma_{i},\gamma_{j}]}}$, 
and hence $K_{ji}K\alpha_{j} \in \mathfrak{S}_{j}$.  A similar 
argument shows that $\psi: \mathfrak{S}_{j} \rightarrow 
\mathfrak{S}_{i}$ given by $K'\alpha_{j} \stackrel{\psi}{\longmapsto} 
K_{ij}K'\alpha_{i}$ is well-defined.  That $\phi$ and $\psi$ are 
inverses follows from the fact that $K_{ij}K_{ji} = 1$.\hfill\QED. 

\noindent 
{\bf \ReducedLemma}\ \ 
{\sl Suppose $(\Gamma,M)$ is unital OA-cyclic.  Then 
for any OA-path $\mathcal{P}$ there is a simple OA-path  
which is equivalent to $\mathcal{P}$.}

{\em Proof.} If $\mathcal{P} = [\gamma_{i_{0}},\ldots,\gamma_{i_{p}}]$ 
is not simple, then 
let $\gamma_{i_{s}}$ be the first repeated node, appearing again 
(say) as $\gamma_{i_{t}}$ in position $t > s$ of the sequence.  Let 
$\mathcal{P}_{1} = [\gamma_{i_{0}},\ldots,\gamma_{i_{s}}]$, 
$\mathcal{Q} = [\gamma_{i_{s}},\ldots,\gamma_{i_{t}}]$ (an OA-cycle), and 
$\mathcal{P}_{2} = [\gamma_{i_{t}},\ldots,\gamma_{i_{p}}]$. 
Clearly $\mathcal{P} = 
\mathcal{P}_{1}\natural\mathcal{Q}\natural\mathcal{P}_{2}$.  
Since $\Pi_{_{\mathcal{Q}}} = 1$, then $\Pi_{_{\mathcal{P}}} = 
\Pi_{_{\mathcal{P}_{1}\natural\mathcal{P}_{2}}}$.  So we have 
$\mathcal{P}_{1}\natural\mathcal{P}_{2} \sim \mathcal{P}$, and the 
former is shorter than the latter.  
Continuing this process we arrive at a simple 
OA-path equivalent to $\mathcal{P}$.\hfill\QED

\noindent 
{\bf \TFAE}\ \ 
{\sl Suppose $(\Gamma',M')$ is an OA-connected component of 
$(\Gamma,M)$ with nodes corresponding to some subset $J \subseteq 
I_{n}$. Then the following are equivalent:} 

\hspace{0.35in} {\sl (1) $(\Gamma',M')$ is unital OA-cyclic;} 

\hspace{0.35in} {\sl (2) $|\{K\alpha_{x}|K\in\mathbb{R}\}\cap\Phi_{M}^{+}| < 
\infty$ for some $x \in J$;} 

\hspace{0.35in} {\sl (3) $|\{K\alpha_{x}|K\in\mathbb{R}\}\cap\Phi_{M}^{+}| < 
\infty$ for all $x \in J$.}   

{\sl In these cases we have $|\{K\alpha_{x}|K\in\mathbb{R}\}\cap\Phi_{M}^{+}| = 
|\{K\alpha_{y}|K\in\mathbb{R}\}\cap\Phi_{M}^{+}|$ for all $x, y \in J$.}

{\em Proof.}  We show (2) $\Rightarrow$ (1) $\Rightarrow$ (3), the 
implication (3) $\Rightarrow$ (2) being obvious.  
For (1) 
$\Rightarrow$ (3), let $x \in J$.  
Observe that if $K\alpha_{x} \in \Phi_{M}^{+}$, then 
by \HumphreysTheorem\ we must have $K = \Pi_{_{\mathcal{P}}}$ for 
some OA-path $\mathcal{P}$ with end node $\gamma_{x}$.  Therefore 
$\mathcal{P}$ is in $(\Gamma',M')$.  
By \ReducedLemma, we may take a simple  
OA-path $\mathcal{Q}$ equivalent to $\mathcal{P}$ 
(all OA-paths equivalent to 
$\mathcal{P}$ must be in $(\Gamma',M')$), so that $K = 
\Pi_{_{\mathcal{Q}}}$.  Since there can be at most a finite number of 
simple OA-paths, then there can be at most finitely many positive roots that 
are scalar multiples of a given $\alpha_{x}$.  For (2) $\Rightarrow$ 
(1), we show the contrapositive.  Let $\mathcal{C} = 
[\gamma_{x},\ldots,\gamma_{x}]$ be a non-unital OA-cycle with 
start/end node $\gamma_{x}$ for an $x \in J$. So 
necessarily $\mathcal{C}$ has nonzero length.  Note that 
$w_{_{\mathcal{C}}}.\alpha_{x} = \Pi_{_{\mathcal{C}}}\alpha_{x}$. 
Next, for $y \in J$ (and possibly $y = 
x$) take any OA-path $\mathcal{P}$ 
with start node $\gamma_{x}$ 
and end node $\gamma_{y}$.  Since 
$w_{_{\mathcal{P}}}.\alpha_{x} = \Pi_{_{\mathcal{P}}}\alpha_{y}$, it 
follows that $w_{_{\mathcal{P}}}w_{_{\mathcal{C}}}^{k}.\alpha_{x} = 
\Pi_{_{\mathcal{P}}}\Pi_{_{\mathcal{C}}}^{k}\alpha_{y}$ for any integer $k$.  
In particular, for all $y \in I_{n}$, we have 
$|\{K\alpha_{y}|K\in\mathbb{R}\}\cap\Phi_{M}^{+}| = \infty$. The final 
claim of the proposition statement follows from \OneToOneLemma.\hfill\QED

When the E-GCM $(\Gamma,M)$ is OA-connected and unital OA-cyclic, let 
$f_{\Gamma,M} := |\{K\alpha_{x}|K\in\mathbb{R}\}\cap\Phi_{M}^{+}|$ for 
any fixed $x \in I_{n}$; then in this case   
\HumphreysTheorem, \PositiveToNegativeRootsLemma, and \TFAE\ 
allow us to modify the proof of Proposition 5.6 of \cite{HumCoxeter} 
to obtain the result that for all $w \in W$, 
$|N_{M}(w)| = f_{\Gamma,M}\, \myl(w)$. 
\LengthProposition\ below generalizes this statement.    
When $W$ is infinite, the length function must take arbitrarily large 
values; from \LengthProposition\ it will follow that 
$\Phi_{M}$ is infinite as well.  

\noindent 
{\bf \LengthProposition}\ \ 
{\sl If the E-GCM graph $(\Gamma,M)$ is unital OA-cyclic then for all $w \in 
W$ we have}
\[f_{1}\, \myl(w) \leq |N_{M}(w)| \leq f_{2}\, \myl(w),\]
{\sl where $f_{1}$ is the min and $f_{2}$ is the max of all integers 
in the set} 
\[\{f_{\Gamma',M'}\, |\, (\Gamma',M') \mbox{ is an 
OA-connected component of } (\Gamma,M)\}.\]

{\em Proof.} Induct on $\myl(w)$.  For $\myl(w) = 0$, the result is 
obvious.  Now take $w$ with $w = w's_{i}$, $\myl(w) = \myl(w')+1$, 
and $\gamma_{i}$ in an OA-connected component $(\Gamma',M')$ 
of $(\Gamma,M)$.  Then by 
\PositiveToNegativeRootsLemma, $|N_{M}(w)| = |N_{M}(w')| + 
f_{\Gamma',M'}$.  Since $f_{1}\, \myl(w') \leq |N_{M}(w')| \leq 
f_{2}\, \myl(w')$, the result follows.\hfill\QED

\noindent 
{\bf \UnitalExample}\ \ In \UnitalExampleFigure\ is depicted a connected, 
unital OA-cyclic E-GCM graph $(\Gamma,M)$ 
with two OA-connected components: $(\Gamma_{1},M_{1})$ is the E-GCM 
subgraph with nodes $\gamma_{i}$ and $\gamma_{j}$, and 
$(\Gamma_{2},M_{2})$ has nodes $\gamma_{x}$, $\gamma_{y}$, and 
$\gamma_{z}$.  
The amplitude matrix $M$ is not symmetrizable by Exercise 2.1 of 
\cite{Kac}. 
Pertaining to the pair $(\gamma_{y},\gamma_{z})$, we have 
$4\cos^{2}(\pi/5) = \frac{3+\sqrt{5}}{2}$ and $2\cos(\pi/5) = 
\frac{1+\sqrt{5}}{2}$.  Since $M_{yz} = -\frac{1+\sqrt{5}}{4}$ and 
$M_{zy} = -(1+\sqrt{5})$, 
then $K_{yz} = \frac{-M_{yz}}{2\cos(\pi/5)} = 
\frac{1}{2}$ and $K_{zy} = \frac{-M_{zy}}{2\cos(\pi/5)} = 2$.  
For all other odd adjacencies 
$(\gamma_{p},\gamma_{q})$ in this graph, $m_{pq} = 3$, so $K_{pq} = 
-M_{pq}$ and $K_{qp} = -M_{qp}$.  
Use \HumphreysTheorem\ to see that 
$f_{\Gamma_{1},M_{1}} = 2$ and $f_{\Gamma_{2},M_{2}} = 3$.  For 
example, 
$\{K\alpha_{i}|K\in\mathbb{R}\}\cap\Phi_{M}^{+} = 
\{\alpha_{i}, \frac{1}{5}\alpha_{i}\} = N_{M}(s_{i})$ and 
$\{K\alpha_{x}|K\in\mathbb{R}\}\cap\Phi_{M}^{+} = 
\{\alpha_{x}, \frac{1}{7}\alpha_{x}, \frac{2}{7}\alpha_{x}\} = 
N_{M}(s_{x})$.  
By \LengthProposition, we can see that 
\[f_{\Gamma_{1},M_{1}}\myl(s_{x}s_{i}) = 4 \leq |N_{M}(s_{x}s_{i})| 
\leq 6 = f_{\Gamma_{2},M_{2}}\myl(s_{x}s_{i}).\] 
More precisely, from \PositiveToNegativeRootsLemma, we get 
\[N_{M}(s_{x}s_{i}) = s_{i}(N_{M}(s_{x}))
\disjointunion 
(\{K\alpha_{i}|K\in\mathbb{R}\}\cap\Phi_{M}^{+}),\]
from which we see that $|N_{M}(s_{x}s_{i})| = 5$.\hfill\QED 

\begin{figure}[ht]
\begin{center}
\UnitalExampleFigure: A unital OA-cyclic E-GCM graph for 
\UnitalExample. \\
{\footnotesize (In this figure, when the amplitude product on an edge is unity, 
we place$\!\!\!$ \CircleInteger{3}$\!\!\!$ beside the edge for emphasis.)}

\ 

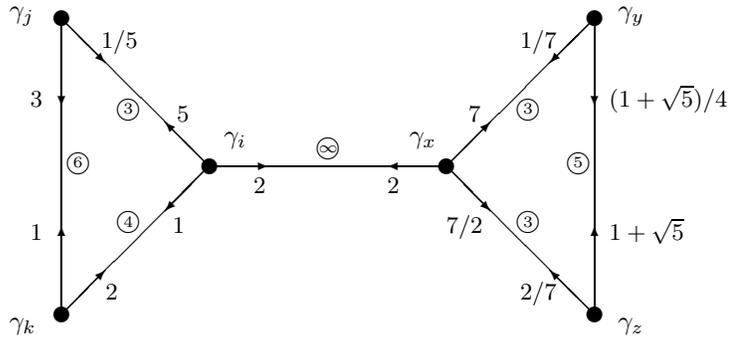
\BowTie   
\end{center}
\end{figure}

We have the natural pairing $\langle \lambda, v \rangle := 
\lambda(v)$ for elements $\lambda$ in the dual space $V^{*}$ and 
vectors $v$ in $V$. 
We think of $V^{*}$ as the space of positions for 
numbers games played on $(\Gamma,M)$:   
For $\lambda \in V^{*}$, the populations for the corresponding 
position are $(\lambda_{i})_{i \in I_{n}}$ where for each $i \in 
I_{n}$ we have $\lambda_{i} := \langle \lambda, \alpha_{i} 
\rangle$.  
Regard the fundamental positions $(\omega_{i})_{i \in I_{n}}$ to be the 
basis for $V^{*}$ dual to the basis $(\alpha_{j})_{j \in I_{n}}$ for 
$V$ relative to the natural pairing $\langle \cdot,\cdot \rangle$, so $\langle 
\omega_{i}, \alpha_{j} \rangle = \delta_{ij}$. 
Given $\sigma_{M}: W \rightarrow GL(V)$ the 
contragredient representation 
$\sigma_{M}^{*}: W \rightarrow GL(V^{*})$ 
is determined by $\langle \sigma_{M}^{*}(w)(\lambda), v \rangle = 
\langle \lambda, \sigma_{M}(w^{-1})(v) \rangle$.  From here on, when $w 
\in W$ and $\lambda \in V^{*}$, 
write $w.\lambda$ for $\sigma_{M}^{*}(w)(\lambda)$.  
Then $s_{i}.\lambda$ is the result of firing node $\gamma_{i}$ when 
the E-GCM graph is assigned position $\lambda$, whether the firing is 
legal or not. 
Set $P_{i} := \{\lambda \in V^{*}\, |\, \langle 
\lambda, \alpha_{i} \rangle > 0\}$, the set of positions with 
positive population at node $\gamma_{i}$.  
Similarly define 
$P_{i}' := \{\lambda \in V^{*}\, |\, \langle 
\lambda, \alpha_{i} \rangle < 0\}$ and 
$Z_{i} := \{\lambda \in V^{*}\, |\, \langle 
\lambda, \alpha_{i} \rangle = 0\}$. 
Let $C := \cap_{i \in I_{n}}P_{i}$, 
the set of strongly dominant positions, and 
let $D := \overline{C} = 
\cap_{i \in I_{n}}(P_{i}\cup{Z}_{i})$, the set of dominant positions.  
The {\em Tits 
cone} is $U_{M} := \cup_{w \in W}wD$.  
In view of 
\HumphreysTheorem, the results of \cite{HumCoxeter} \S 5.13 hold 
here.  And in view of \HRTProofResults, we can use the proof of 
Proposition 3.2 of \cite{HRT} verbatim to get the following 
generalization of their result;  their proof requires that all 
$N_{M}(w)$ be finite, hence our hypothesis that 
$(\Gamma,M)$ is unital OA-cyclic.

\noindent 
{\bf \HRTResult}\ \ {\sl Suppose $(\Gamma,M)$ is connected and unital 
OA-cyclic.  If the Coxeter group $W = W(\Gamma,M)$ 
is infinite, then 
$U_{M} \cap (-U_{M}) = \{0\}$.}\hfill\QED
 
Lemma 5.13 of \cite{HumCoxeter} is the basis for the argument in 
\S 4 of \cite{ErikssonDiscrete} 
characterizing the set of initial positions for which the game 
converges. In contrast to \cite{ErikssonDiscrete}, 
here we fire at nodes with positive rather than negative 
populations, so we have $-U_{M}$ instead of $U_{M}$ in the following 
proposition 
statement.  

\noindent 
{\bf \TitsConeConvergenceResult\ (Eriksson)}\ \ 
{\sl The set of initial positions for which the numbers game on the 
E-GCM graph $(\Gamma,M)$ converges is precisely $-U_{M}$.} 

Although the next proposition is not used in any subsequent proofs, 
it is closely related to the results of this section (cf.\ Exercise 
5.13 of \cite{HumCoxeter}). 

\noindent 
{\bf \TitsConeProp}\ \ 
{\sl Suppose $W = W(\Gamma,M)$ is finite. Then $U_{M} = -U_{M} = V^{*}$.} 

{\em Proof.} Since $W$ is finite, then by Eriksson's Reduced Word 
Result for E-GCM graphs (Proposition 4.1 of 
\cite{ErikssonDiscrete}) it follows that the set of initial positions 
for which the numbers game on $(\Gamma,M)$ converges is all of 
$V^{*}$.  
\TitsConeConvergenceResult\ now implies that $-U_{M} = 
V^{*}$, which is therefore $U_{M}$.\hfill\QED

To understand admissibility for E-GCM graphs that are not 
unital OA-cyclic, we will 
revisit parts of the proof of \MarsAttacksTheorem\ from Section 
\FirstProofNum.  
We note that  
Eriksson's Strong Convergence Theorem  (Theorem 3.1 of 
\cite{ErikssonEur}, stated here for GCM graphs as 
\StrongConvergenceTheorem), 
Eriksson's Comparison Theorem (Theorem 4.5 of 
\cite{ErikssonDiscrete}, stated here for GCM graphs as 
\ComparisonTheorem), and 
Eriksson's Reduced Word Result (Proposition 4.1 of 
\cite{ErikssonDiscrete}, stated above for GCM graphs as 
\ErikssonWordProposition) all hold for E-GCM graphs.  
Moreover, \LemmaListForEGCM\ all hold for E-GCM 
graphs since their proofs nowhere depend on the assumption of integral 
amplitude products. 
We say the $n$-node graph $\Gamma$ is a {\em loop} if the nodes 
can be numbered $\gamma_{1},\ldots,\gamma_{n}$  in such a way that 
for all $1 \leq i \leq n$, 
$\gamma_{i}$ is adjacent precisely to $\gamma_{i+1}$ 
and $\gamma_{i-1}$, understanding that $\gamma_{0} = \gamma_{n}$ 
and $\gamma_{1} = \gamma_{n+1}$. 

\noindent 
{\bf \CycleLemma}\ \ 
{\sl Suppose that the underlying 
graph $\Gamma$ of the E-GCM graph $(\Gamma,M)$ is a loop 
and that for any pair of adjacent nodes the amplitude product is one.  
Then $(\Gamma,M)$ is not admissible.} 

{\em Proof.} We find a divergent 
game sequence starting from 
the fundamental position $\omega_{1}$. Then by renumbering the nodes, 
we see that every fundamental position will have a divergent game 
sequence, and by \NotMarsFriendlyLemma\ it then follows that 
$(\Gamma,M)$ is not admissible. 
Let the OA-cycle $\mathcal{C}$ 
be $[\gamma_{1},\gamma_{2},\ldots,\gamma_{n},\gamma_{1}]$. 
From initial position $\omega_{1}$ 
we propose starting with the firing 
sequence 
$(\gamma_{1},\ldots,\gamma_{n-1},
\gamma_{n},\gamma_{n-1},\ldots\gamma_{2})$.  One can check that all 
of these node firings are legal and that the resulting populations are zero 
at all nodes 
other than $\gamma_{1}$, $\gamma_{2}$, and $\gamma_{n}$.  The 
populations at the latter nodes are, respectively, 
$1+\Pi_{_{\mathcal{C}}}+\Pi_{_{\mathcal{C}}}^{-1}$, 
$M_{12}(\Pi_{_{\mathcal{C}}}^{-1})$, and 
$M_{1n}(\Pi_{_{\mathcal{C}}})$. 
By repeating the proposed firing sequence 
$(\gamma_{1},\ldots,\gamma_{n-1},
\gamma_{n},\gamma_{n-1},\ldots\gamma_{2})$ from this position 
we obtain zero 
populations everywhere except at $\gamma_{1}$, $\gamma_{2}$, and 
$\gamma_{n}$, which are now 
$1+\Pi_{_{\mathcal{C}}}+\Pi_{_{\mathcal{C}}}^{-1}+
\Pi_{_{\mathcal{C}}}^{2}+\Pi_{_{\mathcal{C}}}^{-2}$, 
$M_{12}(\Pi_{_{\mathcal{C}}}^{-1}+\Pi_{_{\mathcal{C}}}^{-2})$, and 
$M_{1n}(\Pi_{_{\mathcal{C}}}+\Pi_{_{\mathcal{C}}}^{2})$ 
respectively.  After $k$ applications of the proposed firing sequence 
we have populations 
$1+\sum_{j=1}^{k}\Pi_{_{\mathcal{C}}}^{j}+\Pi_{_{\mathcal{C}}}^{-j}$, 
$M_{12}(\sum_{j=1}^{k}\Pi_{_{\mathcal{C}}}^{-j})$, and 
$M_{1n}(\sum_{j=1}^{k}\Pi_{_{\mathcal{C}}}^{j})$ 
at nodes 
$\gamma_{1}$, $\gamma_{2}$, and 
$\gamma_{n}$, and zero populations elsewhere.  Thus we have exhibited 
a divergent game sequence.\hfill\QED 

The proof of the next lemma can be adjusted to account for loops on 
three nodes whose E-GCM graphs have even adjacencies. 

\noindent 
{\bf \ThreeCycleLemma}\ \ 
{\sl Suppose $(\Gamma,M)$ is the following three-node E-GCM graph:}  
\hspace*{0.1in}
\parbox[c]{0.5in}{
\setlength{\unitlength}{0.75in}
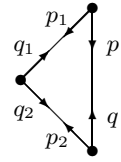
\begin{picture}(0.6,1.2)
            \put(0,0.6){\circle*{0.075}}
            \put(0.5,0.1){\circle*{0.075}}
            \put(0.5,1.1){\circle*{0.075}}
            \put(0,0.6){\line(1,1){0.5}}
            \put(0,0.6){\line(1,-1){0.5}}
            \put(0.5,0.1){\line(0,1){1}}
            \put(0.5,0.3){\vector(0,1){0.1}}
            \put(0.6,0.3){\footnotesize $q$}
            \put(0.5,0.9){\vector(0,-1){0.1}}
            \put(0.6,0.8){\footnotesize $p$}
            \put(0,0.6){\vector(1,1){0.2}}
            \put(0.5,1.1){\vector(-1,-1){0.2}}
            \put(0,0.6){\vector(1,-1){0.2}}
            \put(0.5,0.1){\vector(-1,1){0.2}}
            \put(-0.05,0.8){\footnotesize $q_{1}$}
            \put(0.175,1){\footnotesize $p_{1}$}
            \put(-0.05,0.35){\footnotesize $q_{2}$}
            \put(0.175,0.15){\footnotesize $p_{2}$} 
\end{picture}
} 
\hspace*{0.25in}{\sl Assume that all adjacencies are odd.  
Then $(\Gamma,M)$ is not admissible.} 

{\em Proof.}  
Below \NotMarsFriendlyCatalog\ in Section \FirstProofNum\ we showed that a 
certain three-node GCM graph is not admissible.  The proof here is 
tedious but 
follows the pattern of that argument. With amplitudes as depicted in 
the lemma statement, assign populations $a$, $b$, 
and $c$ as follows:  
\hspace*{0.1in}
\parbox[c]{0.5in}{
\setlength{\unitlength}{0.75in}
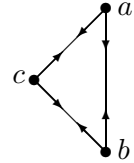
\begin{picture}(0.6,1.2)
            \put(0,0.6){\circle*{0.075}}
            \put(0.5,0.1){\circle*{0.075}}
            \put(0.5,1.1){\circle*{0.075}}
            \put(0,0.6){\line(1,1){0.5}}
            \put(0,0.6){\line(1,-1){0.5}}
            \put(0.5,0.1){\line(0,1){1}}
            \put(0.5,0.3){\vector(0,1){0.1}}
            \put(0.5,0.9){\vector(0,-1){0.1}}
            \put(0,0.6){\vector(1,1){0.2}}
            \put(0.5,1.1){\vector(-1,-1){0.2}}
            \put(0,0.6){\vector(1,-1){0.2}}
            \put(0.5,0.1){\vector(-1,1){0.2}}
            \put(-0.15,0.58){$c$}
            \put(0.58,0.05){$b$}
            \put(0.58,1.05){$a$}
\end{picture}
} 
\hspace*{0.2in} 
Call this position $\lambda = (a,b,c)$, so population $a$ is at node 
$\gamma_{1}$, population $b$ is at node $\gamma_{2}$, and population 
$c$ is at node $\gamma_{3}$.  Without loss of generality, assume that 
$pq \leq p_{1}q_{1}$ and that $pq \leq p_{2}q_{2}$.  Set 
\[
\kappa_{1} := \frac{p p_{2} + p_{1}\sqrt{pq}}{\sqrt{pq}(2-\sqrt{pq})} 
\hspace*{0.25in}\mbox{and}\hspace*{0.25in}
\kappa_{2} := \frac{q p_{1} + p_{2}\sqrt{pq}}{\sqrt{pq}(2-\sqrt{pq})}.
\]
Assume that 
$a \geq 0$, $b \geq 0$, $c \leq 0$, 
and that 
$\displaystyle \left(\kappa_{1}-\frac{p}{q_{2}\sqrt{pq}}\right)a + 
\left(\kappa_{2}-\frac{q}{q_{1}\sqrt{pq}}\right)b + c > 0$; these 
hypotheses will be referred to as 
condition ({\tt *}).  Notice that $a$ and $b$ cannot both be zero 
under condition ({\tt *}).  A justification of the following claim 
will be given at the end of the proof: 

\noindent 
{\bf Claim:}\ \ {\sl 
Under condition} ({\tt *}) 
{\sl there is a sequence of legal node firings from initial 
position $\lambda = (a,b,c)$ which results 
in the position $\lambda' = (a', b', c') = 
(\frac{-q}{\sqrt{pq}}b, \frac{-p}{\sqrt{pq}}a, \kappa_{1}a + 
\kappa_{2}b + c)$.}  

In this case, observe that $a' \leq 0$, $b' \leq 0$, and $c' > 
0$.  Now fire at node $\gamma_{3}$ to obtain the position $\lambda^{(1)} = 
(a_{1}, b_{1}, c_{1})$ with $a_{1} = 
q_{1}[\kappa_{1}a + 
(\kappa_{2}-\frac{q}{q_{1}\sqrt{pq}})b + c]$, 
$b_{1} = q_{2}[(\kappa_{1}-\frac{p}{q_{2}\sqrt{pq}})a + 
\kappa_{2}b + c]$, and $c_{1} = -(\kappa_{1}a + 
\kappa_{2}b + c)$. Now condition ({\tt *}) implies that $a_{1} > 0$, 
$b_{1} > 0$, and $c_{1} < 0$.  
At this point to see that $\lambda^{(1)} = (a_{1}, b_{1}, c_{1})$ itself meets 
condition ({\tt *}), we only need to show that 
$\displaystyle \left(\kappa_{1}-\frac{p}{q_{2}\sqrt{pq}}\right)a_{1} + 
\left(\kappa_{2}-\frac{q}{q_{1}\sqrt{pq}}\right)b_{1} + c_{1} > 0$.  
As a first 
step, we argue that ({\em i}) $q_{1}(\kappa_{1}-\frac{p}{q_{2}\sqrt{pq}}) \geq 
1$ and that ({\em ii}) $q_{2}(\kappa_{2}-\frac{q}{q_{1}\sqrt{pq}}) \geq 1$.  
We only show ({\em i}) since ({\em ii}) follows by similar 
reasoning.  
(From the inequalities ({\em i}) and ({\em ii}), a third inequality 
({\em iii}) follows 
immediately:  
$q_{1}(\kappa_{1}-\frac{p}{q_{2}\sqrt{pq}}) + 
q_{2}(\kappa_{2}-\frac{q}{q_{1}\sqrt{pq}}) - 1 > 0$.)  
For the first of the inequalities ({\em i}), note that since $1 \leq pq$, 
then $2-\sqrt{pq} \leq pq$.  Since $pq \leq p_{2}q_{2}$, then 
$2-\sqrt{pq} \leq p_{2}q_{2}$.  (Similarly $2-\sqrt{pq} \leq 
p_{1}q_{1}$.) Thus 
$\frac{p_{2}q_{2}}{2-\sqrt{pq}}-1 \geq 0$, and hence 
$\frac{p_{2}}{2-\sqrt{pq}} - \frac{1}{q_{2}} \geq 0$.  Therefore, 
$\frac{q_{1}pp_{2}}{\sqrt{pq}(2-\sqrt{pq})} - 
\frac{q_{1}p}{q_{2}\sqrt{pq}} \geq 0$.  Since 
$\frac{p_{1}q_{1}\sqrt{pq}}{\sqrt{pq}(2-\sqrt{pq})} \geq 1$, then 
$\frac{q_{1}pp_{2}}{\sqrt{pq}(2-\sqrt{pq})} + 
\frac{p_{1}q_{1}\sqrt{pq}}{\sqrt{pq}(2-\sqrt{pq})} - 
\frac{q_{1}p}{q_{2}\sqrt{pq}} \geq 1$. From this we get 
$q_{1}(\kappa_{1}-\frac{p}{q_{2}\sqrt{pq}}) \geq 1$, which is ({\em i}). 
The following identity is easy to verify: 
\begin{eqnarray*}
\left(\kappa_{1}-\frac{p}{q_{2}\sqrt{pq}}\right)a_{1} + 
\left(\kappa_{2}-\frac{q}{q_{1}\sqrt{pq}}\right)b_{1} + c_{1} &  & 
\hspace*{3in}
\end{eqnarray*}
\begin{eqnarray*} 
\hspace*{0.75in} & = & 
\left(\kappa_{1}-\frac{p}{q_{2}\sqrt{pq}}\right)\left[
q_{1}\left(\kappa_{1}-\frac{p}{q_{2}\sqrt{pq}}\right) + 
q_{2}\left(\kappa_{2}-\frac{q}{q_{1}\sqrt{pq}}\right) - 1\right]\,\! 
a_{1}\\
\hspace*{0.75in} &  & +\   
\left(\kappa_{2}-\frac{q}{q_{1}\sqrt{pq}}\right)\left[
q_{1}\left(\kappa_{1}-\frac{p}{q_{2}\sqrt{pq}}\right) + 
q_{2}\left(\kappa_{2}-\frac{q}{q_{1}\sqrt{pq}}\right) - 1\right]\,\! 
b_{1}\\ 
\hspace*{0.75in} &  & +\ 
\left[q_{1}\left(\kappa_{1}-\frac{p}{q_{2}\sqrt{pq}}\right) + 
q_{2}\left(\kappa_{2}-\frac{q}{q_{1}\sqrt{pq}}\right) - 1\right]\,\! 
c_{1}\\ 
\hspace*{0.75in} &  & +\  
\frac{p}{q_{2}\sqrt{pq}}\left[
q_{1}\left(\kappa_{1}-\frac{p}{q_{2}\sqrt{pq}}\right) - 1\right]\,\! 
a_{1} + 
\frac{q}{q_{1}\sqrt{pq}}\left[
q_{2}\left(\kappa_{2}-\frac{q}{q_{1}\sqrt{pq}}\right) - 1\right]\,\! 
b_{1}
\end{eqnarray*}

Now the inequalities ({\em i}), ({\em ii}), and ({\em iii}) 
of the previous paragraph together with the 
inequality $\displaystyle \left(\kappa_{1}-\frac{p}{q_{2}\sqrt{pq}}\right)a + 
\left(\kappa_{2}-\frac{q}{q_{1}\sqrt{pq}}\right)b + c > 0$ 
from condition ({\tt *}) imply that 
$\displaystyle \left(\kappa_{1}-\frac{p}{q_{2}\sqrt{pq}}\right)a_{1} + 
\left(\kappa_{2}-\frac{q}{q_{1}\sqrt{pq}}\right)b_{1} + c_{1} > 0$, as 
desired.  This means that position 
$\lambda^{(1)} = (a_{1}, b_{1}, c_{1})$ meets condition ({\tt *}) and 
none of its populations are zero.  
In view of our {\bf Claim}, we may apply to position $\lambda^{(1)}$ a 
legal sequence of node 
firings followed by firing node $\gamma_{3}$ as before to obtain a 
position $\lambda^{(2)} = (a_{2}, b_{2}, c_{2})$ that meets 
condition ({\tt *}) with none of its populations zero, etc.  
So from any such $\lambda = (a,b,c)$ we have a divergent game 
sequence.  In view of inequalities ({\em i}) and ({\em ii}), 
the fundamental positions 
$\omega_{1} = (1,0,0)$ and $\omega_{2} = (0,1,0)$ meet condition 
({\tt *}).  The fundamental position $\omega_{3} = (0,0,1)$ does not 
meet condition ({\tt *}); however, by firing at node $\gamma_{3}$ we 
obtain the position $(q_{1}, q_{2}, -1)$, which meets  
condition ({\tt *}) by inequality ({\em iii}).  
Thus from any fundamental position there is a 
divergent game sequence, and so by \NotMarsFriendlyLemma\ the 
three-node E-GCM graph we started with is not admissible. 

It still remains to justify our {\bf Claim}.   Beginning with 
position $\lambda = (a,b,c)$ under condition ({\tt *}), we propose to 
fire at nodes 
$\gamma_{1}$ and $\gamma_{2}$ in alternating order until this is no 
longer possible.  
We assert that the resulting population will be 
$\lambda' = (a', b', c') = 
(\frac{-q}{\sqrt{pq}}b, \frac{-p}{\sqrt{pq}}a, \kappa_{1}a + 
\kappa_{2}b + c)$. 
There are three cases to consider: (I), $a$ and $b$ are both 
positive, (II), $a>0$ and $b=0$, and (III), $a=0$ and $b>0$.  
For (I), 
we wish to show that $(\gamma_{1},\gamma_{2},\ldots,\gamma_{1})$ of 
length $m_{12}$ is a sequence of legal node firings.  That is, we 
must check that 
\begin{eqnarray} 
\langle (s_{2}s_{1})^{k}.\lambda , \alpha_{1} \rangle & = & 
\langle \lambda , (s_{1}s_{2})^{k}.\alpha_{1} \rangle \ > \ 0 
\hspace*{0.25in} \mbox{for}\ \ 0 \leq k \leq (m_{12}-1)/2, \mbox{ and}\\
\langle s_{1}(s_{2}s_{1})^{k}.\lambda , \alpha_{2} \rangle & = & 
\langle \lambda , s_{1}(s_{2}s_{1})^{k}.\alpha_{2} \rangle \ > \ 0 
\hspace*{0.25in} \mbox{for}\ \ 0 \leq k < (m_{12}-1)/2 
\end{eqnarray}
For (II), 
we wish to show that 
$(\gamma_{1},\gamma_{2},\ldots,\gamma_{1},\gamma_{2})$ of 
length $m_{12}-1$ is a sequence of legal node firings.  That is, we 
must check that 
\begin{eqnarray} 
\langle (s_{2}s_{1})^{k}.\lambda , \alpha_{1} \rangle & = & 
\langle \lambda , (s_{1}s_{2})^{k}.\alpha_{1} \rangle \ > \ 0 
\hspace*{0.25in} \mbox{for}\ \ 0 \leq k < (m_{12}-1)/2, \mbox{ and}\\
\langle s_{1}(s_{2}s_{1})^{k}.\lambda , \alpha_{2} \rangle & = & 
\langle \lambda , s_{1}(s_{2}s_{1})^{k}.\alpha_{2} \rangle \ > \ 0 
\hspace*{0.25in} \mbox{for}\ \ 0 \leq k < (m_{12}-1)/2 
\end{eqnarray}
For (III), 
we wish to show that 
$(\gamma_{2},\gamma_{1},\ldots,\gamma_{2},\gamma_{1})$ of 
length $m_{12}-1$ is a sequence of legal node firings.  That is, we 
must check that 
\begin{eqnarray} 
\langle (s_{1}s_{2})^{k}.\lambda , \alpha_{2} \rangle & = & 
\langle \lambda , (s_{2}s_{1})^{k}.\alpha_{2} \rangle \ > \ 0 
\hspace*{0.25in} \mbox{for}\ \ 0 \leq k < (m_{12}-1)/2, \mbox{ and}\\
\langle s_{2}(s_{1}s_{2})^{k}.\lambda , \alpha_{1} \rangle & = & 
\langle \lambda , s_{2}(s_{1}s_{2})^{k}.\alpha_{1} \rangle \ > \ 0 
\hspace*{0.25in} \mbox{for}\ \ 0 \leq k < (m_{12}-1)/2 
\end{eqnarray}

Our justification of (1) through (6) has 
similarities to the proof of \TwoGeneratorAnalysis.  Under the 
representation $\sigma_{M}$ we have $S_{i} = \sigma_{M}(s_{i})$ for 
$i = 1,2,3$.  With respect to the ordered basis $\mathfrak{B} = 
(\alpha_{1}, \alpha_{2}, \alpha_{3})$ for $V$ we have 
$X_{1} := [S_{1}]_{\mathfrak{B}} = 
\left(\begin{array}{ccc}-1 & p & p_{1}\\ 0 & 1 & 0\\ 0 & 0 & 
1\end{array}\right)$ and 
$X_{2} := [S_{2}]_{\mathfrak{B}} = 
\left(\begin{array}{ccc}1 & 0 & 0\\ 
q & -1 & p_{2}\\
0 & 0 & 1\end{array}\right)$, and so 
\[X_{1,2} := [S_{1}S_{2}]_{\mathfrak{B}} = X_{1}X_{2} = 
\left(\begin{array}{ccc}
pq-1 & -p & p_{2}p + p_{1}\\
q & -1 & p_{2}\\
0 & 0 & 1
\end{array}\right)\]
and 
\[X_{2,1} := [S_{2}S_{1}]_{\mathfrak{B}} = X_{2}X_{1} = 
\left(\begin{array}{ccc}
-1 & p & p_{1}\\
-q & pq-1 & p_{1}q + p_{2}\\
0 & 0 & 1
\end{array}\right).\]
For (1) through (6) above, we need to understand 
$X_{1,2}^{k}${\footnotesize 
$\left(\begin{array}{c}1\\ 0\\ 0\end{array}\right)$},  
$X_{2}X_{1,2}^{k}${\footnotesize 
$\left(\begin{array}{c}1\\ 0\\ 0\end{array}\right)$},  
$X_{2,1}^{k}${\footnotesize 
$\left(\begin{array}{c}0\\ 1\\ 0\end{array}\right)$}, and  
$X_{1}X_{2,1}^{k}${\footnotesize 
$\left(\begin{array}{c}0\\ 1\\ 0\end{array}\right)$}.   
Set $\theta := \pi/m_{12}$.  
Then we can write $X_{1,2} = PDP^{-1}$ for nonsingular $P$ and 
diagonal matrix $D$ as in 
\[
\frac{1}{q(e^{2i\theta}-e^{-2i\theta})}
\left(\begin{array}{ccc}
e^{2i\theta}+1 & e^{-2i\theta}+1 & p_{2}p + 2p_{1}\\
q & q & p_{1}q+2p_{2}\\
0 & 0 & 4-pq
\end{array}\right)
\left(\begin{array}{ccc}
e^{2i\theta} & 0 & 0\\
0 & e^{-2i\theta} & 0\\
0 & 0 & 1
\end{array}\right)
\left(\begin{array}{ccc}
q & -e^{-2i\theta}-1 & C_{1}\\
-q & e^{2i\theta}+1 & C_{2}\\
0 & 0 & C_{3}
\end{array}\right)
\]
where 
$C_{1} = [-q(p_{2}p+2p_{1})+(e^{-2i\theta}+1)(p_{1}q+2p_{2})]/(4-pq)$, 
$C_{2} = [q(p_{2}p+2p_{1})-(e^{2i\theta}+1)(p_{1}q+2p_{2})]/(4-pq)$, and 
$C_{3} = q(e^{2i\theta}-e^{-2i\theta})/(4-pq)$.  With some work we can 
calculate $X_{1,2}^{k}$, which results in 
\[X_{1,2}^{k} = \left(\begin{array}{ccc}
\frac{\sin(2(k+1)\theta)+\sin(2k\theta)}{\sin(2\theta)} & 
\frac{-p\sin(2k\theta)}{\sin(2\theta)} & C_{1}'\\
\frac{q\sin(2k\theta)}{\sin(2\theta)} & 
\frac{-\sin(2k\theta)-\sin(2(k-1)\theta)}{\sin(2\theta)} & 
C_{2}'\\
0 & 0 & 1
\end{array}\right),\]
with \[C_{1}' = -\frac{p_{2}p+2p_{1}}
{4-pq}\left[\frac{\sin(2(k+1)\theta)+\sin(2k\theta)}
{\sin(2\theta)} 
- 1\right] + 
\frac{p(p_{1}q+2p_{2})\sin(2k\theta)}{(4-pq)\sin(2\theta)}\]  
and \[C_{2}' = -\frac{q(p_{2}p+2p_{1})\sin(2k\theta)}{(4-pq)\sin(2\theta)} + 
\frac{p_{1}q+2p_{2}}
{4-pq}\left[\frac{\sin(2k\theta)+\sin(2(k-1)\theta)}
{\sin(2\theta)} + 1\right].\] 
Similar reasoning (or simply interchanging the roles of $\alpha_{1}$ 
and $\alpha_{2}$ in the preceding calculations, or noting that 
$X_{2,1}^{k} = (X_{1,2}^{-1})^{k} = X_{1,2}^{-k}$\, ) shows that 
\[X_{2,1}^{k} = \left(\begin{array}{ccc}
\frac{-\sin(2k\theta)-\sin(2(k-1)\theta)}{\sin(2\theta)} & 
\frac{p\sin(2k\theta)}{\sin(2\theta)} & C_{1}''\\
\frac{-q\sin(2k\theta)}{\sin(2\theta)} & 
\frac{\sin(2(k+1)\theta)+\sin(2k\theta)}{\sin(2\theta)} & 
C_{2}''\\
0 & 0 & 1
\end{array}\right),\]
with 
\[C_{1}'' =   
\frac{p_{2}p+2p_{1}}
{4-pq}\left[\frac{\sin(2k\theta)+\sin(2(k-1)\theta)}
{\sin(2\theta)} + 1\right] - 
\frac{p(p_{1}q+2p_{2})\sin(2k\theta)}{(4-pq)\sin(2\theta)}\] and 
\[C_{2}'' = \frac{q(p_{2}p+2p_{1})\sin(2k\theta)}{(4-pq)\sin(2\theta)} - 
\frac{p_{1}q+2p_{2}}{4-pq}\left[\frac{\sin(2(k+1)\theta)+\sin(2k\theta)}
{\sin(2\theta)} 
- 1\right].\]  
Then 
\[X_{2}X_{1,2}^{k} = \left(\begin{array}{ccc}
\frac{\sin(2(k+1)\theta)+\sin(2k\theta)}{\sin(2\theta)} & 
\frac{-p\sin(2k\theta)}{\sin(2\theta)} & C_{1}'\\
\frac{q\sin(2(k+1)\theta)}{\sin(2\theta)} &
\frac{(1-pq)\sin(2k\theta)+\sin(2(k-1)\theta)}{\sin(2\theta)} & 
qC_{1}' - C_{2}' + p_{2}\\
0 & 0 & 1
\end{array}\right)\]
and 
\[X_{1}X_{2,1}^{k} = \left(\begin{array}{ccc}
\frac{(1-pq)\sin(2k\theta)+\sin(2(k-1)\theta)}{\sin(2\theta)} & 
\frac{p\sin(2(k+1)\theta)}{\sin(2\theta)} & -C_{1}'' + pC_{2}'' + p_{1}\\
\frac{-q\sin(2k\theta)}{\sin(2\theta)} & 
\frac{\sin(2(k+1)\theta)+\sin(2k\theta)}{\sin(2\theta)} & 
C_{2}''\\
0 & 0 & 1
\end{array}\right).\] 
Now we can justify (1) through (6).  For example, for (4) we see that 
since $X_{1}X_{2,1}^{k}${\footnotesize 
$\left(\begin{array}{c}0\\ 1\\ 0\end{array}\right)$} is the second 
column of the matrix $X_{1}X_{2,1}^{k}$, 
then $\langle \lambda , s_{1}(s_{2}s_{1})^{k}.\alpha_{2} \rangle = 
a\frac{p\sin(2(k+1)\theta)}{\sin(2\theta)}$, which is positive 
since $a > 0$, $p > 0$, and (recalling that $m_{12}$ is odd) 
$2(k+1) < m_{12}$.  

Then the proposed firing sequence for each of cases (I), (II), and 
(III) is legal.  To see in case (I) 
that the resulting position is the claimed 
$\lambda' = (a', b', c') = 
(\frac{-q}{\sqrt{pq}}b, \frac{-p}{\sqrt{pq}}a, \kappa_{1}a + 
\kappa_{2}b + c)$, 
we need to calculate 
$\langle s_{1}(s_{2}s_{1})^{k}.\lambda , \alpha_{i} \rangle = 
\langle \lambda , s_{1}(s_{2}s_{1})^{k}.\alpha_{i} \rangle$ for each 
of $i = 1,2,3$, where $k$ is now $(m_{12}-1)/2$.  With patience one 
can confirm that 
\[X_{1}X_{2,1}^{k} = 
\left(\begin{array}{ccc}
0 & 
-p/\sqrt{pq} & \kappa_{1}\\
-q/\sqrt{pq} & 0 &  
\kappa_{2}\\
0 & 0 & 1
\end{array}\right),\] from which the claim follows.  Similar 
computations confirm the claim for cases (II) and (III).\hfill\QED

We have one more loop to rule out; since the details are by now 
routine, we omit the proof.

\noindent 
{\bf \FourCycleLemma}\ \ 
{\sl An E-GCM graph in the family} 
\hspace*{0.1in}
\parbox[c]{0.75in}{
\setlength{\unitlength}{0.5in} 
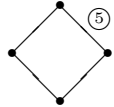
\begin{picture}(1,1.35)
            \put(0,0.6){\circle*{0.075}}
            \put(0.5,0.1){\circle*{0.075}}
            \put(0.5,1.1){\circle*{0.075}}
            \put(1,0.6){\circle*{0.075}}
            \put(0,0.6){\line(1,1){0.5}}
            \put(0,0.6){\line(1,-1){0.5}}
            \put(1,0.6){\line(-1,1){0.5}}
            \put(1,0.6){\line(-1,-1){0.5}}
            \put(0.6,0.85){\CircleInteger{5}}
           \end{picture}
}{\sl is not admissible.}\hfill\QED

We can now prove \EGCMTheorem. 

{\em Proof of \EGCMTheorem.}  
First we use induction on $n$, the number of nodes, to show that any 
connected admissible E-GCM graph $(\Gamma,M)$ must be from one of the 
families of \ECoxeterGraphFigure.  
Clearly a one-node E-GCM graph is admissible.  
For some $n \geq 2$, 
suppose the result is true for all connected admissible E-GCM 
graphs with fewer than $n$ nodes.  Let $(\Gamma,M)$ be a connected, 
admissible, $n$-node E-GCM graph.  
Suppose $(\Gamma,M)$ is 
unital OA-cyclic.  Then 
by \RuleOutInfinitePropositions, we must have $W(\Gamma,M)$ finite.  
Then by the classification of finite irreducible Coxeter groups, 
$(\Gamma,M)$ must be in one of the families of graphs in  
\ECoxeterGraphFigure.  
Now suppose $(\Gamma,M)$ is not unital OA-cyclic.  
First we show that any cycle (OA or otherwise) in $(\Gamma,M)$ must 
use all $n$ nodes.  Indeed, the (connected) 
E-GCM subgraph $(\Gamma',M')$ whose nodes are the nodes of 
a cycle must be admissible by  
\SubgraphLemma.  If $(\Gamma',M')$ has fewer than $n$ nodes, then 
the induction hypothesis applies; but E-Coxeter graphs have no 
cycles (OA or otherwise), so $(\Gamma',M')$ must be all  
of $(\Gamma,M)$. 
Second, $(\Gamma,M)$ 
has an OA-cycle $\mathcal{C}$ for which $\Pi_{_{\mathcal{C}}} \not= 
1$. We can make the following choice for $\mathcal{C}$: Choose 
$\mathcal{C}$ to be a simple OA-cycle with $\Pi_{_{\mathcal{C}}} 
\not= 1$ whose length is as small as possible.  This smallest length 
must therefore be $n$.  
We wish to show that the underlying graph $\Gamma$ 
is a loop.  Let the numbering of the nodes of $\Gamma$ follow 
$\mathcal{C}$, so 
$\mathcal{C} = [\gamma_{1}, \gamma_{2}, \ldots, \gamma_{n}, 
\gamma_{1}]$.  If $\Gamma$ is not a loop, then there are adjacencies 
amongst the $\gamma_{i}$$'${\footnotesize s} besides those of 
consecutive elements of 
$\mathcal{C}$.  But this in turn means that $(\Gamma,M)$ has a 
cycle that uses fewer than $n$ nodes.   So $\Gamma$ is a loop.  Of 
course we must have $n \geq 3$. 
\ThreeCycleLemma\ rules out the possibility that $n=3$.  
Any 
E-GCM subgraph $(\Gamma',M')$ obtained from $(\Gamma,M)$ by removing 
a single node must now be a ``branchless'' E-Coxeter graph from 
\ECoxeterGraphFigure\ whose adjacencies are all odd.  So if $n=4$, 
$(\Gamma,M)$ must be in one of the families 
\hspace*{0.1in}
\parbox[c]{0.75in}{
\setlength{\unitlength}{0.5in} 
\begin{picture}(1,1.35)
            \put(0,0.6){\circle*{0.075}}
            \put(0.5,0.1){\circle*{0.075}}
            \put(0.5,1.1){\circle*{0.075}}
            \put(1,0.6){\circle*{0.075}}
            \put(0,0.6){\line(1,1){0.5}}
            \put(0,0.6){\line(1,-1){0.5}}
            \put(1,0.6){\line(-1,1){0.5}}
            \put(1,0.6){\line(-1,-1){0.5}}
           \end{picture}
} 
or 
\hspace*{0.1in}
\parbox[c]{0.75in}{
\setlength{\unitlength}{0.5in} 
\begin{picture}(1,1.35)
            \put(0,0.6){\circle*{0.075}}
            \put(0.5,0.1){\circle*{0.075}}
            \put(0.5,1.1){\circle*{0.075}}
            \put(1,0.6){\circle*{0.075}}
            \put(0,0.6){\line(1,1){0.5}}
            \put(0,0.6){\line(1,-1){0.5}}
            \put(1,0.6){\line(-1,1){0.5}}
            \put(1,0.6){\line(-1,-1){0.5}}
            \put(0.6,0.85){\CircleInteger{5}}
           \end{picture}
}\hspace*{-0.1in}, which are ruled out by \RuleOutSmallCyclesLemmas\ 
respectively.  
If $n \geq 5$, the only possibility is that $(\Gamma,M)$ meets 
the hypotheses of \CycleLemma\ and therefore 
is not admissible.  In all cases, 
we see that if $(\Gamma,M)$ is not unital OA-cyclic, then it is 
not admissible.  This completes the induction step, so we have shown 
that a connected admissible E-GCM graph must be in one of the 
families of \ECoxeterGraphFigure. 

On the other hand, if $(\Gamma,M)$ is 
from \ECoxeterGraphFigure, then the Coxeter group $W$ 
is finite (again by the classification), 
so there is an upper bound on the length of any element in $W$. 
So by Eriksson's Reduced Word Result for E-GCM graphs, the numbers game 
converges for any initial position. 
The remaining claims of \EGCMTheorem\ now follow from 
Eriksson's Strong Convergence Theorem for E-GCM graphs.\hfill\QED

\vspace{2ex} 

\noindent
{\Large \bf \ConvergenceNum.\ \ 
Some further finiteness aspects of E-game play}

\vspace{2ex} 
The results of this section are independent of the classifications  
obtained in \MainTheorems.  Suppose the Coxeter group $W$ 
associated to an E-GCM graph $(\Gamma,M)$ is finite.  
In this case we may choose the (unique) longest 
element $w_{0}$ in $W$.  
Since we must have $\myl(w_{0}s_{i}) < 
\myl(w_{0})$ for all $i \in I_{n}$, it follows that $w_{0}.\alpha_{i} 
<_{_{M}} 0$ for all $i$.  So if $\alpha = \sum c_{i}\alpha_{i} >_{_{M}} 0$, 
then $w_{0}.\alpha <_{_{M}} 0$, i.e.\   
$N_{M}(w_{0}) = \Phi_{M}^{+}$. 

\noindent 
{\bf \KeyWJResultLemma}\ \ 
{\sl  Let $(\Gamma,M)$ be an E-GCM graph with associated Coxeter 
group $W = W(\Gamma,M)$.  
Let $J \subseteq I_{n}$, 
and suppose the parabolic subgroup 
$W_{J}$ is finite.  
Suppose $\alpha = \sum_{j \in 
J}c_{j}\alpha_{j}$ is a root in $\Phi_{M}^{+}$.  Then for some $w_{_{J}}$ in 
$W_{J}$, we have $w_{_{J}}.\alpha <_{_{M}} 0$.} 

{\em Proof.} Since $W_{J}$ is finite, consider the longest element 
$(w_{0})_{_{J}}$ in $W_{J}$.  Note that any element of $W_{J}$ preserves 
the subspace $V_{J} := \mathrm{span}_{\mathbb{R}}\{\alpha_{j}\}_{j \in 
J}$.  As seen just above, $(w_{0})_{_{J}}$ will send each simple 
root $\alpha_{j}$ for $j \in J$ to some root in $\Phi_{M}^{-}$.  So 
apply $(w_{0})_{_{J}}$ to the given $\alpha$ to see that 
$(w_{0})_{_{J}}.\alpha <_{_{M}} 0$. \hfill\QED

In what follows, for any subset $J$ of 
$I_{n}$, let $C_{J}$ be the set of all dominant positions $\lambda$ 
for which $\lambda = \sum_{i \in 
I_{n}\setminus{J}}\lambda_{i}\omega_{i}$ with $\lambda_{i} > 0$ for 
all $i \in I_{n}\setminus{J}$, that is, $C_{J} := (\cap_{i \in 
I_{n}\setminus{J}} P_{i}) \cap 
(\cap_{j \in J} Z_{j})$.  
The next result generalizes Proposition 4.2 of \cite{ErikssonDiscrete}. 

\noindent 
{\bf \KeyWJResult}\ \ 
{\sl For an E-GCM graph $(\Gamma,M)$ with Coxeter group $W = 
W(\Gamma,M)$, let $J$ be any subset of $I_{n}$ such that $W_{J}$ is finite.  Let 
$s_{i_{p}}\cdots{s}_{i_{2}}s_{i_{1}}$ be any reduced 
expression for an element $w^{^{J}}$ of the set $W^{J}$ of minimal coset 
representatives.  Let $\lambda$ be in $C_{J}$.  
Then for $1 \leq q \leq p$, 
$(\gamma_{i_{1}},\ldots,\gamma_{i_{q}})$ is a legal sequence of node 
firings for a numbers game played from initial position $\lambda$.} 

{\em Proof.}  For $1 \leq q \leq p$, we must show that $\langle 
s_{i_{q-1}}\cdots{s}_{i_{2}}s_{i_{1}}.\lambda , \alpha_{i_{q}} 
\rangle > 0$.  But $\langle 
s_{i_{q-1}}\cdots{s}_{i_{2}}s_{i_{1}}.\lambda , \alpha_{i_{q}} 
\rangle = \langle 
\lambda , s_{i_{1}}{s}_{i_{2}}\cdots{s}_{i_{q-1}}.\alpha_{i_{q}} 
\rangle$.  Now 
$s_{i_{1}}{s}_{i_{2}}\cdots{s}_{i_{q-1}}.\alpha_{i_{q}} >_{_{M}} 0$ 
since $\myl(s_{i_{1}}{s}_{i_{2}}\cdots{s}_{i_{q-1}}s_{i_{q}}) > 
\myl(s_{i_{1}}{s}_{i_{2}}\cdots{s}_{i_{q-1}})$   
(this is because $\myl(s_{i_{q}}s_{i_{q-1}}\cdots{s}_{i_{2}}s_{i_{1}}) > 
\myl(s_{i_{q-1}}\cdots{s}_{i_{2}}s_{i_{1}})$ for these two reduced 
expressions).  Moreover, suppose 
$w_{_{J}}^{-1} = s_{j_{1}}\cdots{s}_{j_{r}}$ is a reduced expression for some 
$w_{_{J}}^{-1} \in W_{J}$.  Then since $w^{^{J}}w_{_{J}}^{-1} = 
s_{i_{p}}\cdots{s}_{i_{2}}s_{i_{1}}s_{j_{1}}\cdots{s}_{j_{r}}$ is 
reduced (cf.\ Proposition 2.4.4 of \cite{BB}), it follows that 
both expressions  
$s_{j_{r}}\cdots{s}_{j_{1}}s_{i_{1}}\cdots{s}_{i_{q-2}}s_{i_{q-1}}$ 
and $s_{j_{r}}\cdots{s}_{j_{1}}s_{i_{1}}\cdots{s}_{i_{q-1}}s_{i_{q}}$ 
are reduced, with the latter longer than the former. 
In particular 
$s_{j_{r}}\cdots{s}_{j_{1}}s_{i_{1}}
\cdots{s}_{i_{q-2}}s_{i_{q-1}}.\alpha_{i_{q}} = w_{_{J}}s_{i_{1}}
\cdots{s}_{i_{q-2}}s_{i_{q-1}}.\alpha_{i_{q}} >_{_{M}} 0$ for all $w_{_{J}} 
\in W_{J}$.  
We wish to show that 
$s_{i_{1}}{s}_{i_{2}}\cdots{s}_{i_{q-1}}.\alpha_{i_{q}}$ cannot be 
contained in $\mathrm{span}_{\mathbb{R}}\{\alpha_{j}\}_{j \in J}$.  
Suppose otherwise, so  
$s_{i_{1}}{s}_{i_{2}}\cdots{s}_{i_{q-1}}.\alpha_{i_{q}} = 
\sum_{j \in J}c_{j}\alpha_{j}$.  We now invoke the finiteness 
assumption for $W_{J}$ and \KeyWJResultLemma\ as follows:  
Let $(w_{0})_{_{J}}$ be 
the longest element of $W_{J}$.  Then 
$(w_{0})_{_{J}}s_{i_{1}}{s}_{i_{2}}\cdots{s}_{i_{q-1}}.\alpha_{i_{q}} <_{_{M}} 
0$, a contradiction.  
Then it must be the case that 
$s_{i_{1}}{s}_{i_{2}}\cdots{s}_{i_{q-1}}.\alpha_{i_{q}} = 
\sum_{i \in I_{n}}c_{i}\alpha_{i}$ with $c_{k} > 0$ for some $k \in 
I_{n}\setminus{J}$.  So 
$\langle 
s_{i_{q-1}}\cdots{s}_{i_{2}}s_{i_{1}}.\lambda , \alpha_{i_{q}} 
\rangle = \langle 
\lambda , s_{i_{1}}{s}_{i_{2}}\cdots{s}_{i_{q-1}}.\alpha_{i_{q}} 
\rangle 
= \langle \lambda , \sum_{i \in I_{n}}c_{i}\alpha_{i} \rangle 
= \sum_{i \in I_{n}}c_{i}\lambda_{i}$, which is positive since all 
$c_{i}$'s are nonnegative, $\lambda_{k} > 0$, and $c_{k} > 0$.  Then 
firing at node $\gamma_{i_{q}}$ from game position 
$s_{i_{q-1}}\cdots{s}_{i_{2}}s_{i_{1}}.\lambda$ is legal.\hfill\QED

From the previous proof we see that 
\KeyWJResult\ could be generalized to other kinds of subsets 
$J \subseteq I_{n}$ 
by extending \KeyWJResultLemma.  Doing so for all subsets $J 
\subseteq I_{n}$ would yield a simple proof of \EGCMTheorem: An E-GCM 
graph $(\Gamma,M)$ would have a convergent game sequence for some 
$\lambda \in C_{J}$ if and only if $W^{J}$ is finite if and only if 
$W$ if finite (by Proposition 4.2 of \cite{Deodhar}).  For an 
arbitrary E-GCM graph $(\Gamma,M)$, let 
$\mathfrak{P}(\lambda)$ denote the set of all positions obtainable 
from legal firing sequences in numbers games with initial position 
$\lambda$.  Clearly $\mathfrak{P}(\lambda) \subseteq W\lambda$, where 
the latter is the orbit of $\lambda$ under the $W$-action on 
$V^{*}$.  
Since the 
statement of  
Theorem 5.13 of \cite{HumCoxeter} holds for quasi-standard geometric 
representations, then $W_{J}$ is the full stabilizer of any $\lambda 
\in C_{J}$, so $W\lambda$ and $W^{J}$ can be identified. 
So from \KeyWJResult\ we see that for $\lambda \in C_{J}$ with  
$W_{J}$ finite, then $\mathfrak{P}(\lambda) = W\lambda$.  

\noindent 
{\bf \LengthUsingLongestWord}\ \ 
{\sl Suppose $W = W(\Gamma,M)$ is finite for some E-GCM graph 
$(\Gamma,M)$.  Let $J \subseteq I_{n}$ be any subset, and let 
$(w_{0})_{_{J}}$ be the longest word in the parabolic subgroup $W_{J}$.  
Then all game sequences for any $\lambda \in C_{J}$ have length 
$\myl(w_{0}) - \myl((w_{0})_{_{J}})$.} 

{\em Proof.} Write $w_{0} = (w_{0})^{^{J}}(w_{0})_{_{J}}$, with 
$(w_{0})^{^{J}}$ the minimal coset representative for $w_{0}W_{J}$.  
\KeyWJResult\ implies that there is a game sequence for $\lambda$ with 
length $\myl((w_{0})^{^{J}}) = \myl(w_{0}) - \myl((w_{0})_{_{J}})$.  
By Eriksson's 
Strong Convergence Theorem, this must be the length of any game 
sequence for $\lambda$.\hfill\QED

Our next two results expand on Remark 4.6 of \cite{ErikssonDiscrete}. 
Let $\lambda$ be a strongly dominant position thought of as an 
initial position for a numbers game on an E-GCM graph.  Consider a 
legal sequence of node firings 
$(\gamma_{i_{1}},\ldots,\gamma_{i_{p}})$.  
Let $\alpha = \sum_{i \in I_{n}}c_{i}\alpha_{i}$ be the root 
$s_{i_{1}}\cdots{s}_{i_{p}}.\alpha_{k}$.  
Then  the population at node 
$\gamma_{k}$ can be computed as $\langle 
s_{i_{p}}\cdots{s}_{i_{1}}.\lambda , \alpha_{k} \rangle = 
\langle \lambda , s_{i_{1}}\cdots{s}_{i_{p}}.\alpha_{k} \rangle = 
\langle \lambda , \alpha \rangle = 
\sum_{i \in I_{n}}c_{i}\lambda_{i}$.  
For indeterminates 
$x_{1},\ldots,x_{n}$, let $\phi_{\alpha} := 
\phi_{\alpha}(x_{1},\ldots,x_{n}) := \sum_{i \in I_{n}}c_{i}x_{i}$, 
and call $\phi_{\alpha}$ the {\em root functional for} $\alpha$; a 
root functional is {\em positive} or {\em negative} depending on 
whether $\alpha$ is positive or negative.  
By \KeyWJResult\ (or Proposition 4.2 of \cite{ErikssonDiscrete}), 
this game sequence is valid for all strongly dominant $\lambda$, so 
we will say that 
$\phi_{\alpha}$ is the {\em root functional at node} $\gamma_{k}$ 
{\em for 
the firing sequence} $(\gamma_{i_{1}},\ldots,\gamma_{i_{p}})$.  If $k = 
i_{p}$ and the firing sequence is understood, then we just say 
$\phi_{\alpha}$ is the {\em root functional at node} 
$\gamma_{i_{p}}$.  
Note that $\phi_{\alpha} = \phi_{\beta}$ for roots $\alpha$ and 
$\beta$ if and only if $\alpha = \beta$.  
The same positive root functional can appear at different nodes 
during a given numbers game (e.g.\ play the numbers game on the GCM 
graph $A_{3}$ with a generic strongly dominant initial position 
$\lambda = (a,b,c)$; fire nodes $\gamma_{3}$, $\gamma_{2}$, and 
$\gamma_{1}$ in that order to see that the root functional 
$\phi_{\alpha_{2}+\alpha_{3}}$ appears twice). However,   
our next result shows that positive root 
functionals are not repeated among the \underline{fired} 
nodes in any numbers game. 

\noindent 
{\bf \NoRepeatPosRootsProp}\ \ 
{\sl Let $\lambda$ be a strongly dominant position and let 
$(\gamma_{i_{1}}, \gamma_{i_{2}}, \ldots)$ be any game sequence for a 
numbers game played on an E-GCM graph $(\Gamma,M)$ from initial 
position $\lambda$.  Then for all $k \geq 1$, the positive 
root functional at node 
$\gamma_{i_{k}}$ is not the same as the positive root functional at 
any $\gamma_{i_{j}}$ for $1 \leq j < k$.} 

{\em Proof.} Let $\beta_{j} := 
s_{i_{1}}\cdots{s}_{i_{j-1}}.\alpha_{i_{j}}$, which is necessarily 
positive.  Then $\phi_{\beta_{j}}$ is the root functional at node 
$\gamma_{i_{j}}$ for the given game sequence.  Suppose that for some 
$k > j$, we have $\phi_{\beta_{k}} = \phi_{\beta_{j}}$.  
Since  $\beta_{j} = \beta_{k}$, then one can see that 
$s_{i_{j}}\cdots{s}_{i_{k-1}}.\alpha_{i_{k}} = \alpha_{i_{j}}$, and 
so $s_{i_{j+1}}\cdots{s}_{i_{k-1}}.\alpha_{i_{k}} = -\alpha_{i_{j}} 
<_{_{M}} 0$.  But $s_{i_{j+1}}\cdots{s}_{i_{k-1}}s_{i_{k}}$ is reduced 
and longer than $s_{i_{j+1}}\cdots{s}_{i_{k-1}}$, which means we must 
have $s_{i_{j+1}}\cdots{s}_{i_{k-1}}.\alpha_{i_{k}} >_{_{M}} 0$.  
From this contradiction we conclude that we cannot have 
$\phi_{\beta_{k}} = \phi_{\beta_{j}}$ for any $k > j$.\hfill\QED

This gives us a strategy for generating all of the positive roots 
when $\Phi_{M}^{+}$ is finite.  However, if $(\Gamma,M)$ has odd asymmetries, 
then not every positive root will be encountered 
as a positive root functional in a single game sequence, as the 
following result shows. 

\noindent 
{\bf \AllPositiveRootsProp}\ \ {\sl For an E-GCM 
graph $(\Gamma,M)$, suppose $W = W(\Gamma,M)$ is finite.  
Let $(\gamma_{i_{1}},\ldots,\gamma_{i_{l}})$ be any game 
sequence for a strongly dominant 
position $\lambda$.  Then the following are equivalent:} 

\hspace{0.05in} {\sl (1) For each $\alpha \in \Phi_{M}^{+}$, $\phi_{\alpha}$ is 
the positive root functional at some node $\gamma_{i_{j}}$ for the game  
sequence;} 

\hspace{0.05in} {\sl (2) Each  
OA-connected component $(\Gamma',M')$ of $(\Gamma,M)$ is unital 
OA-cyclic with $f_{\Gamma',M'} = 1$;} 

\hspace{0.05in} {\sl (3) $(\Gamma,M)$ 
has no odd asymmetries;} 
 
\hspace{0.05in} {\sl (4) $\myl(w_{0}) = |\Phi_{M}^{+}| = l$.}

{\em Proof.} To show (1) $\Rightarrow$ (2), 
choose an OA-connected component  $(\Gamma',M')$.  The proof of 
\TFAE\ shows that $(\Gamma',M')$ must be unital OA-cyclic, else $W$ will 
be infinite. 
Let $J$ be the subset of 
$I_{n}$ corresponding to the nodes of the subgraph $\Gamma'$.  
As in the proof of \LengthUsingLongestWord, write 
$w = w_{0} = (w_{0})^{^{J}}(w_{0})_{_{J}} = w^{^{J}}w_{_{J}}$, where $w_{_{J}} 
= (w_{0})_{_{J}}$ is the 
longest word in $W_{J}$ and $w^{^{J}} = (w_{0})^{^{J}}$ is the minimal coset 
representative for $w_{0}W_{J}$.  
Set $w_{_{J}} = s_{i_{k}}\cdots{s}_{i_{2}}s_{i_{1}}$, a reduced expression.  
Using \PositiveToNegativeRootsLemma, we see that 
\begin{eqnarray*} 
|N_{M}(ws_{i_{1}})| & = & |N_{M}(w)| - f_{\Gamma',M'},\\
|N_{M}(ws_{i_{1}}s_{i_{2}})| & = & |N_{M}(ws_{i_{1}})| - 
f_{\Gamma',M'}\ \, =\ \, |N_{M}(w)| - 2f_{\Gamma',M'},
\end{eqnarray*}
so that eventually $|N_{M}(w)| = |N_{M}(w^{^{J}})| + 
\myl(w_{_{J}})f_{\Gamma',M'}$.  Now by hypothesis 
each positive root functional 
appears once and therefore, by \NoRepeatPosRootsProp, exactly once.  
Thus any game sequence for $\lambda$ has length $|\Phi_{M}^{+}|$.  
So  by \LengthUsingLongestWord\ 
we see that $\myl(w) = |\Phi_{M}^{+}|$.  By 
\LengthProposition, $|N_{M}(w^{^{J}})| \geq \myl(w^{^{J}})$.  
Summarizing, 
\[\myl(w^{^{J}}) + \myl(w_{_{J}}) = \myl(w) = |N_{M}(w)| = 
|N_{M}(w^{^{J}})| + \myl(w_{_{J}})f_{\Gamma',M'} \geq 
\myl(w^{^{J}}) + \myl(w_{_{J}})f_{\Gamma',M'},\]
from which $f_{\Gamma',M'} = 1$.  For (2) $\Leftrightarrow$ (3), 
note that by \HumphreysTheorem\ we get a nontrivial positive multiple 
of some simple root if and only if there are odd asymmetries.  
For (2) $\Rightarrow$ (4), 
\LengthProposition\ and 
the fact that $f_{\Gamma',M'} = 1$ for each OA-connected component 
tell us that $\myl(w) = |N_{M}(w)| = |\Phi_{M}^{+}|$, which by 
\LengthUsingLongestWord\ is $l$.  For (4) 
$\Rightarrow$ (1), see that when $\myl(w) = |N_{M}(w)| = 
|\Phi_{M}^{+}|$, then 
\PosRootsPropList\ imply that 
each positive root functional 
must appear at least once in any game sequence for 
$\lambda$.\hfill\QED 

At the end of our proof of \MarsAttacksTheorem, we saw that for any 
connected Dynkin diagram, the lengths of the convergent game  
sequences for any two strongly dominant initial positions are  the 
same.  In view of the previous result, this common value is 
the number of positive roots in the standard 
root system associated to the Weyl group (cf.\ \cite{HumCoxeter} \S 
2.10, 2.11).  
For $\mathcal{H}_{3}$ and $\mathcal{H}_{4}$ the lengths of the longest 
words are $15$ and $60$ respectively.  To obtain these values consult 
\cite{HumCoxeter} \S 2.13, or in light of \AllPositiveRootsProp\ 
just play the numbers game on the appropriate E-GCM graphs with no odd 
asymmetries. 

In some related work with Norman Wildberger \cite{DW}, we will 
take an interest in 
what we call here ``adjacency-free positions.''  
For a firing sequence $(\gamma_{i_{1}}, \gamma_{i_{2}}, \ldots)$ from 
a position $\lambda$, then any position 
$s_{i_{j}}\cdots{s}_{i_{1}}.\lambda$ (including $\lambda$ itself) 
is an {\em intermediate 
position} for the sequence.  
A position $\lambda$ is {\em adjacency-free} if there  
exists a game sequence played from $\lambda$ such that no intermediate 
position has a pair of adjacent nodes with 
positive populations.  
We will see that this notion is related to the notion of 
``full commutativity'' of Coxeter group 
elements studied by Stembridge in \cite{StemFC}.  
In the discussion that follows, we view a Coxeter group $W$ 
as $W(\Gamma,M)$ for some E-GCM 
graph $(\Gamma,M = (M_{ij})_{i,j \in I_{n}})$.  
Following \S 1.1 of \cite{StemFC} and \S 8.1 of \cite{HumCoxeter},   
we let $\mathcal{W} = I_{n}^{*}$ 
be the free monoid on the set $I_{n}$.  
Elements of $\mathcal{W}$ are {\em words} and will be viewed as 
finite sequences of elements from $I_{n}$; the binary operation is 
concatenation, and the identity $\varepsilon$ is the 
empty word. Fix a word 
$\mathbf{s} := ({i_{1}},\ldots,{i_{r}})$. 
Then $\myl_{\mathcal{W}}(\mathbf{s}) := r$ is the {\em length} of 
$\selt$.  
A {\em subword} of  
$\mathbf{s}$ 
is any subsequence 
$(i_{p},i_{p+1},\ldots,i_{q})$ of consecutive elements of $\selt$.  
For a nonnegative integer $m$ and $x,y \in I_{n}$, 
let $\langle x,y \rangle_{m}$ denote the sequence 
$(x,y,x,y,\ldots) \in \mathcal{W}$ so that 
$\myl_{\mathcal{W}}(\langle x,y \rangle_{m}) = m$. 
We employ several types of ``elementary simplifications'' in $\mathcal{W}$.  
An {\em elementary simplification of 
braid type} replaces a subword $\langle x,y \rangle_{m_{xy}}$ with 
the subword $\langle y,x \rangle_{m_{xy}}$ if $2 \leq m_{xy} < \infty$.   
An {\em elementary simplification of length-reducing type} 
replaces a subword $(x,x)$ with 
the empty subword. 
We let $\mathcal{S}(\mathbf{s})$ be the set of all words that can be 
obtained from $\mathbf{s}$ by some sequence of elementary 
simplifications of braid or length-reducing type.  
Since $s_{i}$ in $W$ is its own inverse for each $i \in I_{n}$, there is an 
induced mapping $\mathcal{W} \rightarrow W$.  We compose this with 
the mapping $W \rightarrow W$ for which $w \mapsto w^{-1}$ to get 
$\psi: \mathcal{W} \rightarrow W$ given by $\psi(\mathbf{s}) 
= s_{i_{r}}\cdots{s}_{i_{1}}$.  
Tits' Theorem for the word problem on Coxeter 
groups (cf.\ Theorem 8.1 of \cite{HumCoxeter}) 
implies that: {\sl For words $\mathbf{s}$ and $\mathbf{t}$ in 
$\mathcal{W}$, $\psi(\mathbf{s}) = \psi(\mathbf{t})$ if and only if 
$\mathcal{S}(\mathbf{s}) \cap \mathcal{S}(\mathbf{t}) \not= 
\emptyset$.} (This theorem is the basis for Eriksson's Reduced Word 
Result.) 
We say $\mathbf{s}$ is a {\em reduced} word for $w = \psi(\mathbf{s})$ if 
$\myl_{\mathcal{W}}(\mathbf{s}) = \myl(w)$ (assume this is the case 
for the remainder of the paragraph); let $\mathcal{R}(w) 
\subseteq \mathcal{W}$ denote  
the set of all reduced words for $w$.  Suppose that $\mathbf{t} \in 
\mathcal{R}(w)$.  By Tits' Theorem, $\mathcal{S}(\mathbf{s}) \cap 
\mathcal{S}(\mathbf{t}) \not= \emptyset$, 
so that $\mathbf{t}$ can be obtained from 
$\mathbf{s}$ by a sequence of elementary simplifications of braid or 
length-reducing type.  Since $\myl_{\mathcal{W}}(\mathbf{s}) = \myl(w) = 
\myl_{\mathcal{W}}(\mathbf{t})$, then no elementary simplifications 
of length-reducing type can be used to obtain $\mathbf{t}$ from 
$\mathbf{s}$.  Then any member of $\mathcal{R}(w)$ can be obtained 
from any other member by a sequence of elementary simplifications of 
braid type.  
An {\em elementary simplification of commuting type} 
replaces a subword $(x,y)$ 
with the subword 
$(y,x)$ 
if $m_{xy} = 2$. 
The {\em commutativity class} $\mathcal{C}(\mathbf{s})$ of 
the word $\mathbf{s}$ is the set of all words that can be obtained 
from $\mathbf{s}$ by a sequence of elementary simplifications of 
commuting type.  
Clearly 
$\mathcal{C}(\mathbf{s}) \subseteq \mathcal{R}(w)$.  In fact there is a 
decomposition of $\mathcal{R}(w)$ into commutativity classes: 
$\mathcal{R}(w) = \mathcal{C}_{1} 
\disjointunion 
\cdots 
\disjointunion 
\mathcal{C}_{k}$, a disjoint union.  If 
$\mathcal{R}(w)$ has just one commutativity class, then $w$ is {\em 
fully commutative}.  Proposition 1.1 
of \cite{StemFC} states: {\sl An element $w \in W$ is fully 
commutative if and only if for all $x, y \in I_{n}$ such that $3 \leq 
m_{xy} < \infty$, there is no member of $\mathcal{R}(w)$ that 
contains $\langle x,y \rangle_{m_{xy}}$ as a subword}. 
The following is a variation of this result. 

\noindent 
{\bf \AdjacencyFreeLemma}\ \ 
{\sl For an E-GCM graph $(\Gamma,M)$, 
an element $w \in W(\Gamma,M)$ is fully commutative if and only if 
there is a 
commutativity class $\mathcal{C}$ of $\mathcal{R}(w)$ such that for 
all $x, y \in I_{n}$ with $3 \leq m_{xy} < \infty$, no 
member of $\mathcal{C}$ contains $\langle {x},{y} 
\rangle_{m_{xy}}$ as a subword.} 

{\em Proof.} By Proposition 1.1 
of \cite{StemFC} the ``$\Rightarrow$'' direction is clear, so we will 
show the ``$\Leftarrow$'' direction.  Say $\mathcal{C} = 
\mathcal{C}(\mathbf{s})$ for some $\mathbf{s} \in \mathcal{R}(w)$.  
Now take any commutativity class 
$\mathcal{C}(\mathbf{t})$ for $\mathcal{R}(w)$.   Since 
$\psi(\mathbf{s}) = \psi(\mathbf{t})$, then $\mathcal{S}(\mathbf{s})$ 
meets $\mathcal{S}(\telt)$ by Tit's Theorem.  Since 
$\myl_{\mathcal{W}}(\mathbf{s}) = \myl(w) = 
\myl_{\mathcal{W}}(\mathbf{t})$, then $\mathbf{t}$ can be obtained 
from $\mathbf{s}$ using only elementary simplifications of braid 
type.  By hypothesis, only elementary simplifications of 
length-reducing or commuting type can be applied to any 
member of the commutativity class 
$\mathcal{C}$.  In particular, it must be the case that $\mathbf{t} 
\in \mathcal{C}$, so $\mathcal{R}(w)$ has only one commutativity 
class.  Then $w$ is fully commutative.\hfill\QED

The significance of 
the next two results is discussed in 
the first two paragraphs of Section \RemarksNum. 

\noindent 
{\bf \AdjacencyFreeForFinite}\ \ 
{\sl For an E-GCM graph $(\Gamma,M)$, suppose the Coxeter group $W = 
W(\Gamma,M)$ is finite.  Let $J \subseteq I_{n}$. 
(1) Suppose an adjacency-free position $\lambda$ is in $C_{J}$.  Then 
every element $w^{^{J}}$ of $W^{J}$ is fully commutative, and for any 
reduced expression $w^{^{J}} = s_{i_{k}}\cdots{s}_{i_{1}}$, no 
intermediate position for the firing sequence 
$(\gamma_{i_{1}},\ldots,\gamma_{i_{k}})$ from $\lambda$ has positive 
populations on adjacent nodes. (2) Suppose each element of $W^{J}$ is 
fully commutative.  Then any position in $C_{J}$ is adjacency-free.}

{\em Proof.} Write $w_{0} = (w_{0})^{^{J}}(w_{0})_{_{J}}$, with 
$(w_{0})^{^{J}} \in W^{J}$ 
and $(w_{0})_{_{J}}$ longest in $W_{J}$.  Let $L := \myl(w_{0}) - 
\myl(\, (w_{0})_{_{J}})$.  
Our proof of (1) is by induction on the lengths of elements in 
$W^{J}$.  It is clear that the identity element is fully commutative. 
Now suppose that for all $v^{^{J}}$ in $W^{J}$ with $\myl(v^{^{J}}) < 
k$, it is the case that $v^{^{J}}$ is fully commutative and that for 
any reduced expression $v^{^{J}} = s_{i_{p}}\cdots{s}_{i_{1}}$, no 
intermediate position for the firing sequence 
$(\gamma_{i_{1}},\ldots,\gamma_{i_{p}})$ starting at $\lambda$ 
has positive populations on 
adjacent nodes.  Now consider $w^{^{J}}$ in $W^{J}$ such that 
$\myl(w^{^{J}}) = k$.  Suppose that for some adjacent $\gamma_{x}$ 
and $\gamma_{y}$ in $\Gamma$ with $3 \leq m_{xy} < \infty$, we have 
$\langle x,y \rangle_{m_{xy}}$ as a 
subword of some reduced word $\mathbf{s} = (i_{1},\ldots,i_{k}) \in 
\mathcal{R}(w^{^{J}})$.  Since $(i_{1},\ldots,i_{k-1})$ is a reduced 
word and $s_{i_{k-1}}\cdots{s}_{i_{1}}$ is in $W^{J}$, then 
$\langle x,y \rangle_{m_{xy}}$ cannot be a subword of 
$(i_{1},\ldots,i_{k-1})$.  Therefore 
it must be 
the case that $\mathbf{s} = 
(i_{1},\ldots,i_{p},\langle x,y \rangle_{m_{xy}})$ for 
$p = k-m_{xy}$.  
But by \KeyWJResult, the corresponding firing sequence 
is legal, so in 
particular after the first $p$ firings of the sequence there must be 
positive populations on adjacent nodes $\gamma_{x}$ and $\gamma_{y}$. 
This contradicts the induction hypothesis for element 
$v^{^{J}} = s_{i_{p}}\cdots{s}_{i_{1}}$.  So $w^{^{J}}$ is fully 
commutative.  Now suppose that the firing sequence 
$(\gamma_{j_{1}},\ldots,\gamma_{j_{k}})$ corresponding to some reduced 
word $\mathbf{t} = (j_{1},\ldots,j_{k})$ in $\mathcal{R}(w^{^{J}})$ 
results in positive 
populations at two adjacent nodes, say $\gamma_{x}$ and $\gamma_{y}$, 
for which $3 \leq m_{xy} < \infty$.  This means that for some 
$j_{k+m_{xy}+1},\ldots,j_{L}$ the word $\selt := 
(j_{1},\ldots,j_{k},\langle x,y \rangle_{m_{xy}},j_{k+m_{xy}+1},\ldots,j_{L})$ 
corresponds to a game sequence played from $\lambda$ and is a reduced 
word for some $u$ in $W$ of length $L$.  Write $u = u^{^{J}}u_{_{J}}$ 
with $u^{^{J}} \in W^{J}$.  By \KeyWJResult, any reduced expression 
for $u^{J}$ corresponds to a legal firing sequence from $\lambda$; 
since $u^{^{J}}.\lambda = u.\lambda$ is the terminal position for any 
game sequence played from $\lambda$, then this legal firing sequence 
must be a game sequence, and hence by \LengthUsingLongestWord\ we 
have $\myl(u^{^{J}}) = L$.  Then $u_{_{J}} = e$, so $u = u^{^{J}} \in 
W^{J}$.  Since $\myl(u(w_{0})_{_{J}}) = L+\myl(\, 
(w_{0})_{_{J}}) = \myl(w_{0})$, then $u(w_{0})_{_{J}} = w_{0} = 
(w_{0})^{^{J}}(w_{0})_{_{J}}$, so $u = (w_{0})^{^{J}}$. So $\selt$ is 
a reduced word for $(w_{0})^{^{J}}$.  
By Proposition 1.1 of \cite{StemFC}, we see then that $(w_{0})^{^{J}}$ is 
not fully commutative.  Therefore by \AdjacencyFreeLemma, every 
commutativity class of $\mathcal{R}((w_{0})^{^{J}})$ has a member 
containing such a subword.  But this means that every firing sequence 
corresponding to a reduced expression for $(w_{0})^{^{J}}$ must 
have intermediate positions with 
positive populations at adjacent nodes, which contradicts 
the hypothesis that $\lambda$ is an adjacency-free position.  
Therefore no reduced 
word $\mathbf{t} = (j_{1},\ldots,j_{k})$ in $\mathcal{R}(w^{^{J}})$ 
has firing sequence 
$(\gamma_{j_{1}},\ldots,\gamma_{j_{k}})$ which 
results in positive 
populations at two adjacent nodes for an intermediate position.  
This completes the induction step and the proof of (1).

For part (2), assume every member of $W^{J}$ is fully commutative, and 
let $\lambda$ be any position in $C_{J}$.  Let $(i_{1},\ldots,i_{L})$ 
be a reduced word for $(w_{0})^{^{J}}$.  Suppose an intermediate 
position $s_{i_{k}}\cdots{s}_{i_{1}}.\lambda$ for the game sequence 
$(\gamma_{i_{1}},\ldots,\gamma_{i_{L}})$ 
has positive populations 
on adjacent nodes $\gamma_{x}$ and $\gamma_{y}$.  Then by Eriksson's 
Strong Convergence Theorem, there is a game sequence of length $L$ 
from $\lambda$ corresponding to a reduced word $\selt = 
(i_{1},\ldots,i_{k},\langle x,y 
\rangle_{m_{xy}},j_{k+m_{xy}+1},\ldots,j_{L})$ for $u = \psi(\selt)$.  
As in the proof of part (1), we can see that $u$ is just 
$(w_{0})^{^{J}}$.  So $(w_{0})^{^{J}}$ is fully commutative (by 
hypothesis) and has reduced word $\selt$, in violation of Proposition 
1.1 of \cite{StemFC}.  Therefore $\lambda$ is adjacency-free.\hfill\QED

In Theorem 5.1 of \cite{StemFC}, Stembridge classifies those $W^{J}$ 
for irreducible Coxeter groups $W$ such that every member of $W^{J}$ 
is fully commutative.  In view of \AdjacencyFreeForFinite\ and the 
classification of finite Coxeter groups, 
we may apply this result here to conclude that for finite irreducible $W$, 
the adjacency-free dominant positions are exactly those 
specified in the following theorem. Observe that a dominant 
position $\lambda$ is adjacency-free if and only if $r\lambda := 
(r\lambda_{i})_{i \in I_{n}}$ is adjacency-free for all positive real 
numbers $r$; call any such $r\lambda$ a {\em positive multiple} of 
$\lambda$.  

\noindent 
{\bf \AdjacencyFreeClassification}\ \ {\sl Suppose the E-GCM graph 
$(\Gamma,M)$ is connected.  If $W = W(\Gamma,M)$ is 
finite, then an adjacency-free dominant position is a positive 
multiple of a fundamental position.  
All fundamental positions for any E-Coxeter graph of type 
$\mathcal{A}_{n}$ are adjacency-free.  The 
adjacency-free 
fundamental positions for any graph of type 
$\mathcal{B}_{n}$, 
$\mathcal{D}_{n}$, or $\mathcal{I}_{2}(m)$ are 
precisely those corresponding to end nodes. The adjacency-free 
fundamental positions for any graph of type 
$\mathcal{E}_{6}$, $\mathcal{E}_{7}$, 
or $\mathcal{H}_{3}$ are precisely those 
corresponding to the nodes marked with asterisks in 
\ECoxeterGraphFigure.  Any graph of type $\mathcal{E}_{8}$,  
$\mathcal{F}_{4}$, or $\mathcal{H}_{4}$ has no 
adjacency-free fundamental positions.}\hfill\QED

\vspace{2ex} 

\noindent
{\Large \bf \RemarksNum.\ \ Comments}

\vspace{2ex} 
Continue to think of a Coxeter or Weyl group as a group associated 
to some E-GCM graph with index set $I_{n}$.  For $J \subseteq I_{n}$, 
let $(W^{J},\leq)$ denote the Bruhat (partial) 
order on $W^{J}$ (see for example 
\cite{BB} Ch.\ 2 for a definition).  
In Proposition 3.1 of \cite{PrEur}, 
Proctor shows that for finite irreducible Weyl groups $W$, those 
$(W^{J},\leq)$ for which the Bruhat order is a lattice have $|J| = 
1$ and correspond  
precisely to the adjacency-free fundamental positions for Dynkin 
diagrams identified in \AdjacencyFreeClassification; in Proposition 
3.2 of that paper, he shows 
that these lattices 
are, in fact, distributive. 
For finite irreducible Coxeter groups $W$, it is a consequence of 
Theorems 5.1 and 6.1 of \cite{StemFC} that $(W^{J},\leq)$ is a lattice 
if and only if $(W^{J},\leq)$ is a distributive lattice if and only if 
each element of $W^{J}$ is fully commutative, in which case $|J| = 
1$ and all such $J$'s correspond with the adjacency-free fundamental 
positions from \AdjacencyFreeClassification\ above.  
\AdjacencyFreeForFinite\ above adds to 
these equivalences the property that each element of $W^{J}$ is fully 
commutative if and only if for any associated E-GCM graph, 
any position in $C_{J}$ is adjacency-free.  The adjacency-free 
viewpoint is similar to Proctor's original viewpoint (cf.\ Lemma 3.2 
of \cite{PrEur}). 

In \cite{DonTwoColor} we use \MarsAttacksTheorem\ to 
show that if a finite ``edge-colored'' ranked 
poset meets a certain condition relative to an $n \times n$ matrix 
$M$, then $M$ must be a Cartan matrix, i.e.\ $M$ must be a GCM such 
that $(\Gamma,M)$ is a Dynkin diagram.  This so-called ``structure  
property'' is necessary for edge-colored ranked posets to carry certain 
information about semisimple Lie algebra representations.  Indeed, 
identifying such combinatorial properties is part of our program for 
obtaining combinatorial models for Lie algebra representations.  For 
example, in \cite{ADLMPW}, we introduced four families of finite 
distributive lattices whose elements are ``weighted'' by a simple 
combinatorial rule; their ``weight-generating functions'' are Weyl 
characters for the irreducible representations of the four rank two 
semisimple Lie algebras.  In \cite{DonTwoColor} we show that the 
posets of join irreducibles (cf.\ \cite{Stanley}) for these 
distributive lattices are characterized by 
a short list of combinatorial properties.  These are called 
``semistandard posets'' in \cite{ADLMPW}; the smallest of these posets 
are called ``fundamental posets.''  In \cite{DW}, we will say how 
these fundamental and semistandard posets can be constructed from 
information obtained by playing the numbers game on two-node Dynkin 
diagrams.  More generally, we will show how to construct fundamental 
posets for the adjacency-free fundamental positions for any given Dynkin 
diagram.  These fundamental posets can be combined to obtain 
semistandard posets whose corresponding distributive lattices are 
models (as above) for Weyl characters for certain irreducible 
representations of the corresponding semisimple Lie algebra.  When an 
adjacency-free fundamental position for a Dynkin diagram corresponds 
to a ``minuscule'' fundamental weight (see \cite{PrEur}, 
\cite{PrDComplete}, 
\cite{StemFC}), then our fundamental poset is a vertex-colored 
version of the corresponding ``wave'' poset 
of \cite{PrDComplete} and of the corresponding ``heap'' of 
\cite{StemFC}.  

Our work here leaves open some questions about the relationship 
between the numbers game and quasi-standard geometric representations of 
Coxeter groups.  For example, in the notation of Section 
\ConvergenceNum, if $J \subseteq I_{n}$ and $W_{J}$ is 
infinite, then for $\lambda \in C_{J}$ is it the case that 
$\mathfrak{P}(\lambda) = W\lambda$? 
In \AdjacencyFreeForFinite, can the assumption that $W$ is finite be 
relaxed?   Also, when $(\Gamma,M)$ is 
not unital OA-cyclic, what can be said about $U_{M} \cap (-U_{M})$? 
In a different direction, one could ask to what extent the proofs of 
\MainTheorems\ can be combinatorialized, avoiding Coxeter groups 
altogether.  For \MarsAttacksTheorem, one only needs a combinatorial 
proof of Eriksson's Comparison Theorem; such a proof is given for 
``unweighted'' GCM graphs in \cite{ErikssonLinear}. 

\renewcommand{\baselinestretch}{1}

\end{document}